\documentclass{elsart}

\usepackage[dvips]{color}

\usepackage{graphicx}
\usepackage{psfrag}

\usepackage{amssymb}
\usepackage{amsmath}
\usepackage{dsfont}
\usepackage{rotating}

\newcommand{\halb}{\frac{1}{2}}

\renewcommand{\u}{\mathbf{u}}
\newcommand{\w}{\mathbf{w}}
\newcommand{\q}{\mathbf{q}}
\newcommand{\be}{\begin{equation}}
\newcommand{\ee}{\end{equation}}
\newcommand{\bdm}{\begin{displaymath}}
\newcommand{\edm}{\end{displaymath}}
\newcommand{\bea}{\begin{eqnarray}}
\newcommand{\eea}{\end{eqnarray}}

\journal{Journal of Computational Physics}
\begin{document}
\begin{frontmatter}
\title{ADER-WENO Finite Volume Schemes with Space-Time Adaptive Mesh Refinement}

\author[UNITN]{Michael Dumbser}
\ead{michael.dumbser@ing.unitn.it}
\author[UNITN]{Olindo Zanotti}
\ead{olindo.zanotti@ing.unitn.it}
\author[UPM]{Arturo Hidalgo}
\ead{arturo.hidalgo@upm.es}
\author[ND]{Dinshaw S. Balsara}
\ead{dinshaw.balsara@nd.edu}
\address[UNITN]{Laboratory of Applied Mathematics, Department of Civil, Environmental and Mechanical Engineering, University of Trento, Via Mesiano 77, I-38123 Trento, Italy}
\address[UPM]{Departamento de Matem\'atica Aplicada y M\'etodos Inform\'aticos, Universidad Polit\'ecnica de Madrid, Calle R\'ios Rosas 21, E-28003 Madrid, Spain}
\address[ND]{Physics Department, University of Notre Dame du Lac, 225 Nieuwland Science Hall, Notre Dame, IN 46556, USA}

\begin{abstract}
We present the first high order one-step ADER-WENO finite volume scheme with Adaptive Mesh Refinement (AMR) in multiple space 
dimensions. High order spatial accuracy is obtained through a WENO reconstruction, while a high order one-step time discretization 
is achieved using a local space-time discontinuous Galerkin predictor method. Due to the one-step nature of the underlying scheme, 
the resulting algorithm is particularly well suited for an AMR strategy on space-time adaptive meshes, i.e.with time-accurate 
local time stepping. The AMR property has been implemented 'cell-by-cell', with a standard tree-type algorithm, while 
the scheme has been parallelized via the Message Passing Interface (MPI) paradigm. 
The new scheme has been tested over a wide range of examples for nonlinear systems of hyperbolic conservation laws, including the 
classical Euler equations of compressible gas dynamics and the equations of magnetohydrodynamics (MHD). High order in space and 
time have been confirmed via a numerical convergence study and a detailed analysis of the computational speed-up with respect to highly 
refined uniform meshes is also presented. We also show test problems where the presented high order AMR scheme behaves clearly better 
than traditional second order AMR methods. The proposed scheme that combines for the first time high order ADER methods with space--time 
adaptive grids in two and three space dimensions is likely to become a useful tool in several fields of computational physics, applied 
mathematics and mechanics. 
\end{abstract}

\end{frontmatter}

\section{Introduction}

The idea of using fully discrete one-step time update schemes, rather than multi-stage Runge-Kutta 
time integrators, dates back to van Leer~\cite{vanLeer1979} and Harten et al. ~\cite{eno}, who applied 
it to the case of MUSCL and ENO schemes, respectively.
The key feature of one-step time update schemes is that, after a piecewise polynomial approximation to the 
solution has been obtained at time $t^n$ through a prescribed reconstruction procedure, an element--local 
time evolution of such reconstructed polynomials is performed, thus allowing for a better than first order 
computation of the numerical fluxes between adjacent cells. 
The one-step time integrators can also be seen as a particular procedure to solve (approximately) the 
generalized Riemann problem at the element interfaces, where the initial data consists in piecewise polynomials 
instead of piecewise constant data, as it was the case in the original first order Godunov scheme \cite{godunov}.  
The generalized Riemann problem has been discussed in \cite{Raviart.GRP.1,Raviart.GRP.2,LeFloch:1991a,toro4,CastroToro,Montecinos} 
and has been used as a building block of numerical methods e.g. in \cite{Artzi,BenArtzi:2006a} as well as in the ADER approach 
\cite{titarevtoro,titarevtoro2,toro3,taube_jsc,dumbser_jsc,DumbserKaeser07,ADERNC,Balsara2009,Balsara2013}.    
In the original version of the ADER schemes of Titarev and Toro \cite{toro3,titarevtoro,schwartzkopff}, the state 
vector was first expanded in a Taylor series in time at the element interface and second these time derivatives have 
been replaced by spatial derivatives using the so--called Cauchy-Kovalewski procedure that is based on a repeated use 
of the governing conservation law in differential form. The leading state at the interface was computed by the exact 
Riemann solver applied to the boundary extrapolated left and right states. The higher order spatial  derivatives at the interface have then been defined by solving linearized Riemann problems for the  boundary extrapolated values of the 
derivatives, where the linearization was performed about the leading state. Though successful, the Cauchy-Kovalewski 
procedure, which is also equivalently called the Lax--Wendroff procedure \cite{lax}, turns out to be impracticable for 
complex systems of non-linear equations, and alternative methods need to be followed. Other explicit one--step schemes 
in time that are based on the Lax--Wendroff procedure can be found, for example, in 
\cite{qiu,QiuDumbserShu,stedg1,stedg2,chengshu2}. 

While the Cauchy-Kovalewski procedure is based on the \textit{strong} differential form of the PDE, an alternative has been 
proposed in \cite{DumbserEnauxToro,Dumbser2008}, where an element--local space--time Galerkin predictor is introduced that 
is based on the \textit{weak} integral form of the PDE. This approach is also able to account for stiff source terms. 
An overview of explicit one--step time discretization schemes can be found in \cite{GassnerDumbserMunz}. 
%
%
The high order one--step schemes proposed in \cite{DumbserEnauxToro,Dumbser2008} can be divided in three different steps. 
First, a high order WENO reconstruction \cite{shu_efficient_weno,titarevtoro3,HuShuTri,shi,friedrich,kaeserjcp,DumbserKaeser06b,DumbserKaeser07,ZhangShu3D,MixedWENO2D,MixedWENO3D,Pringuey} 
from the given cell averages ${\bf \bar u}$ is performed at time $t^n$, but also other nonlinear high  order reconstructions are possible, e.g. 
\cite{eno,sonar,abgrall_eno,MOOD}. Let the resulting piecewise polynomial solution be denoted by $\mathbf{w}_h=\mathbf{w}_h(\mathbf{x},t^n)$ in the following. Second, 
the local space-time DG predictor \cite{DumbserEnauxToro,DumbserZanotti} is applied for the time evolution of the reconstructed polynomials $\mathbf{w}_h$, 
the result of which are piecewise space--time polynomials denoted by $\mathbf{q}_h=\mathbf{q}_h(\mathbf{x},t)$. 
This allows a high order accurate computation of the numerical fluxes. Finally, the cell averages are updated in time with a one-step scheme according to 
the integral form of the conservation law. 
%
%
The numerical method just described has been successfully applied to a variety of physical systems, including stiff advection-diffusion-reaction problems~\cite{HidalgoDumbser}, 
compressible magnetohydrodynamics~\cite{Dumbser2008,Balsara2009} and magnetic reconnection in an astrophysical context~\cite{ZanottiDumbser2011}. 

Along with high order numerical schemes, another frontier of numerical research is represented by the implementation of efficient AMR algorithms, allowing for numerical simulations 
of very complex  and highly dynamical structures that require an adaptive refinement of the grids in specific flow regions. Starting from the pioneering AMR implementation by 
Berger et al. \cite{Berger-Oliger1984,Berger-Jameson1985,Berger-Colella1989}, who first introduced a patched-based block-structured AMR for finite difference methods, several codes 
have been developed over the years providing AMR infrastructures in combination with  different numerical techniques. We recall here, to list but a few, the AMRCLAW 
package~\cite{LeVequeCLAWPACK,Berger-Leveque1998,Bell1994} based on the second order accurate wave-propagation algorithm; AstroBEAR \cite{Cunningham2009,Carroll-Nellenback2011}, 
where a collection of TVD and piecewise parabolic methods are applied for the spatial reconstruction, while a MUSCL-Hancock predictor-corrector scheme is adopted for the temporal 
evolution; RAMSES \cite{Teyssier2002}, using a tree-based data structure with a second order Godunov method; NIRVANA \cite{Ziegler2008}, with a second order, directionally unsplit 
central-upwind scheme of Godunov type~\cite{Kurganov2005}. 
Further well-known AMR schemes can also be found in the following list of references, \cite{QuirkAMR,Balsara2001,Mulet1,KeppensAMR,FalleAMR1,FalleAMR2}, which does not  
pretend to be complete. 

In the last few years, moreover, significant progress has been obtained in combining high order numerical methods  with AMR techniques. 
It is worth mentioning the fourth order finite volume AMR scheme with Runge-Kutta time integrator by Colella et al.~\cite{Colella2009}; 
the PLUTO-AMR code by Mignone et al.~\cite{Mignone2012}, specifically devised for astrophysical applications,
which includes a Corner-Transport-Upwind scheme for time updating and a third-order WENO spatial reconstruction; 
the spectral WENO schemes by Burger et al.~\cite{Burger2012}, combined with a third-order TVD Runge-Kutta method, 
meant for the study of the sedimentation of a polydisperse suspension. Finite difference WENO schemes together
with AMR have been considered very recently in \cite{FDWENOAMR}. 

Common in most of the above AMR methods is the use of the method of lines (MOL) together with a higher order (TVD) Runge--Kutta scheme in time \cite{shuosher1,shuosher2,shu2}.  
In this paper we improve with respect to the approaches presented so far by providing the first tree-based 'cell-by-cell' AMR algorithm in combination with a finite volume 
ADER-WENO \textit{one-step} time update scheme. The use of high order one--step schemes in time allows directly for time--accurate local time stepping in a straight forward 
manner and has already been successfully applied in the context of high order Discontinuous Galerkin schemes with time--accurate local time stepping 
\cite{dumbserkaeser06d,TaubeMaxwell,stedg1,stedg2}, in multi--block domain decomposition methods \cite{DomainDecomp} and finite volume and DG schemes 
based on static mesh adaptation \cite{CastroLTS}. 

The structure of the paper is the following. In Sect.~\ref{sec:num-approach} we provide a description 
of the ADER approach. Sect.~\ref{sec:AMR} is specifically devoted to the description of the Adaptive 
Mesh  Refinement infrastructure. In Sect.~\ref{sec:num-tests} we present a variety of numerical tests 
for which two different systems of hyperbolic equations in conservative form have been considered, 
namely the classical Euler and magnetohydrodynamics equations. Finally, Sect.~\ref{sec:concl} contains 
the conclusions of our work. 

\section{Numerical Method}
\label{sec:num-approach}

In this section, the numerical method is first presented for purely regular Cartesian meshes. Due to the one--step nature of our high order finite volume schemes, all basic 
ingredients can later be transferred also to AMR grids with only few modifications. 

\subsection{The finite volume scheme}

We consider hyperbolic systems of balance laws in Cartesian coordinates 
\begin{equation}
\label{eq:governingPDE}
\frac{\partial{\bf u}}{\partial t}+\frac{\partial{\bf f}}{\partial x}+\frac{\partial{\bf g}}{\partial y}+\frac{\partial{\bf h}}{\partial z}={\bf S}({\bf u},\mathbf{x},t)\,,
\end{equation}
where $\bf {u}$ is the vector of conserved quantities, while ${\bf f}(\bf{u})$, ${\bf g}(\bf{u})$ and ${\bf h}(\bf{u})$ are the physical flux 
vectors along the $x$, $y$ and $z$ directions, respectively. As usual, we define the three-dimensional control volume on a regular Cartesian grid as  
$I_{ijk}=[x_{i-\halb};x_{i+\halb}]\times[y_{j-\halb};y_{j+\halb}]\times[z_{k-\halb};z_{k+\halb}]$, with 
$\Delta x_i=x_{i+\halb}-x_{i-\halb}$, $\Delta y_j=y_{j+\halb}-y_{j-\halb}$, $\Delta z_k=z_{k+\halb}-z_{k-\halb}$, and, in addition, the spacetime 
control volume ${\mathcal I}_{ijk}=I_{ijk}\times [t^n,t^n+\Delta t]$, where $\Delta t=t^{n+1}-t^n$. 
On adaptive meshes as used in this article, it is convenient to address each three--dimensional control volume with a unique mono-index 
$m$, like in the case of unstructured meshes, with $1<m<N_{{\rm Cells}}$ and $N_{{\rm Cells}}$ being the total number of cells at any given time 
(see Sect.~\ref{AMR-implementation}). Therefore, in the rest of the paper we will also use the symbol ${\mathcal C}_m$ to denote a control volume 
$I_{ijk}$. The numerical scheme, however, is more conveniently written using the indices $i$, $j$ and $k$. The integration over ${\mathcal I}_{ijk}$ 
provides the standard finite volume discretization 
\vspace{-6mm}  
\begin{eqnarray}
\label{eq:finite_vol}
{\bf \bar u}_{ijk}^{n+1}&=&{\bf \bar u}_{ijk}^{n}-\frac{\Delta t}{\Delta x_i}\left({\bf f}_{i+\halb,jk}-{\bf f}_{i-\halb,jk} \right)-\frac{\Delta t}{\Delta y_j}\left({\bf g}_{i,j+\halb,k}-{\bf g}_{i,j-\halb,k} \right)
\nonumber \\
&& \hspace{8mm} -\frac{\Delta t}{\Delta z_k}\left({\bf h}_{ij,k+\halb}-{\bf h}_{ij,k-\halb} \right)+ \Delta t {\bf \bar{S}}_{ijk}, 
\end{eqnarray}
where
\begin{equation}
{\bf \bar u}_{ijk}^{n}=\frac{1}{\Delta x_i}\frac{1}{\Delta y_j}\frac{1}{\Delta z_k}\int \limits_{x_{i-\halb}}^{x_{i+\halb}}\int \limits_{y_{j-\halb}}^{y_{j+\halb}}\int \limits_{z_{k-\halb}}^{z_{k+\halb}}{ \bf u}(x,y,z,t^n)dz\,\,dy\,\,dx
\end{equation}
is the spatial average of the solution in the element ${ I}_{ijk}$ at time $t^n$, while
\begin{equation}
\label{flux:F}
{\bf f}_{i+\halb,jk}= \frac{1}{\Delta t}\frac{1}{\Delta y_j}\frac{1}{\Delta z_k} \hspace{-1mm}  \int \limits_{t^n}^{t^{n+1}} \! \int \limits_{y_{j-\halb}}^{y_{j+\halb}}\int \limits_{z_{k-\halb}}^{z_{k+\halb}} \hspace{-1mm} 
{\bf \tilde f} \! \left({\bf q}_h^-(x_{i+\halb},y,z,t),{\bf q}_h^+(x_{i+\halb},y,z,t)\right) dz \, dy \, dt, 
\end{equation}
\begin{equation}
\label{flux:G}
{\bf g}_{i,j+\halb,k}=\frac{1}{\Delta t}\frac{1}{\Delta x_i}\frac{1}{\Delta z_k} \hspace{-1mm}  \int \limits_{t^n}^{t^{n+1}} \! \int \limits_{x_{i-\halb}}^{x_{i+\halb}}\int \limits_{z_{k-\halb}}^{z_{k+\halb}} \hspace{-1mm} 
{\bf \tilde g} \! \left({\bf q}_h^-(x,y_{j+\halb},z,t),{\bf q}_h^+(x,y_{j+\halb},z,t)\right) dz\,dx\,dt, \\
\end{equation}
\begin{equation}
\label{flux:H}
{\bf h}_{ij,k+\halb}=\frac{1}{\Delta t}\frac{1}{\Delta x_i}\frac{1}{\Delta y_j} \hspace{-1mm}  \int \limits_{t^n}^{t^{n+1}} \! \int \limits_{x_{i-\halb}}^{x_{i+\halb}}\int \limits_{y_{j-\halb}}^{y_{j+\halb}} \hspace{-1mm}  
{\bf \tilde h}\! \left({\bf q}_h^-(x,y,z_{k+\halb},t),{\bf q}_h^+(x,y,z_{k+\halb},t)\right) dy\,dx\,dt, 
\end{equation}
and
\begin{equation}
\label{source:S}
{\bf \bar{S}}_{ijk}=\frac{1}{\Delta t}\frac{1}{\Delta x_i}\frac{1}{\Delta y_j}\frac{1}{\Delta z_k}\int \limits_{t^n}^{t^{n+1}}\int \limits_{x_{i-\halb}}^{x_{i+\halb}}\int \limits_{y_{j-\halb}}^{y_{j+\halb}}\int \limits_{z_{k-\halb}}^{z_{k+\halb}}{ \bf S}\left(\mathbf{q}_h(x,y,z,t)\right) dz\,dy\,dx\,dt\,.
\end{equation}
are the space--time averaged numerical fluxes and sources, respectively. In the integrals above, the symbol $\mathbf{q}_h$ denotes the local space--time DG predictor solution illustrated in section \ref{sec:localDG}. 
In this article, we will use two alternative numerical fluxes, the classical Rusanov flux \cite{Rusanov:1961a,toro-book} or a new variant of the Osher--Solomon flux proposed in \cite{OsherUniversal}. The Rusanov 
flux, also called local Lax-Friedrichs flux in literature, reads 
\begin{equation}
  {\bf \tilde f}\left( \q_h^-, \q_h^+ \right) = \frac{1}{2}\left( \mathbf{f}(\q_h^-) + \mathbf{f}(\q_h^+) \right) - \frac{1}{2} |s_{\max}| \left( \q_h^+ - \q_h^- \right), 
\label{eqn.rusanov} 
\end{equation} 
where $|s_{\max}|$ denotes the maximum absolute value of the eigenvalues of the Jacobian matrix $\mathbf{A} = \partial \mathbf{f} / \partial \u$. 
The Osher--type flux proposed in \cite{OsherUniversal} reads 
\begin{equation}
  {\bf \tilde f}\left( \q_h^-, \q_h^+ \right) = \frac{1}{2}\left( \mathbf{f}(\q_h^-) + \mathbf{f}(\q_h^+) \right) - \frac{1}{2} \left( \int_0^1 \left| \mathbf{A}( \boldsymbol{\psi}(s) ) \right| ds \right) \left( \q_h^+ - \q_h^- \right), 
\label{eqn.osher} 
\end{equation} 
with the straight--line segment path 
\begin{equation}
  \boldsymbol{\psi}(s) = \q_h^- + s \left( \q_h^+ - \q_h^- \right), \qquad 0 \leq s \leq 1,   
\end{equation} 
which connects the left and right state with each other in phase space. The path integral in Eqn. \eqref{eqn.osher} is evaluated \textit{numerically} using Gauss--Legendre quadrature with at 
least two quadrature points \cite{OsherUniversal}. The numerical fluxes $\mathbf{\tilde g}$ and $\mathbf{\tilde h}$ are computed in the same way as the flux $\mathbf{\tilde f}$.

\subsection{ WENO reconstruction}
\label{sec:WENO_reconstruction}
Since the method used here differs to some extent from the original WENO scheme of Jiang and Shu \cite{shu_efficient_weno}, some more details 
are given in the following. For each spatial dimension, we consider a \textit{nodal} basis of polynomials of degree $M$ rescaled on the 
unit interval $I=[0;1]$. We remind that such a basis consists of $M+1$ Lagrange interpolating polynomials of maximum degree $M$, $\{\psi_l(\lambda)\}_{l=1}^{M+1}$, associated to the $M+1$ Gauss-Legendre nodes  $\{\lambda_k\}_{k=1}^{M+1}$ in the interval $[0;1]$ and with the standard property that
\begin{equation}
\psi_l(\lambda_k)=\delta_{lk}\hspace{1cm}l,k=1,2,\ldots, M+1\,.
\end{equation}
Here, $\delta_{lk}$ is the usual Kronecker symbol. This choice produces by construction an 
\textit{orthogonal} basis over $[0;1]$. Moreover, it has the additional advantage that the 
data is immediately available in the Gaussian quadrature points anytime it is  
necessary to compute an integral over $[0;1]$. The reconstruction is performed for each  
cell $I_{ijk}$ on a set of one--dimensional reconstruction stencils, which are given for each Cartesian direction by 
\begin{equation}
\label{eqn.stencildef}  
\mathcal{S}_{ijk}^{s,x} = \bigcup \limits_{e=i-L}^{i+R} {I_{ejk}}, \quad 
\mathcal{S}_{ijk}^{s,y} = \bigcup \limits_{e=j-L}^{j+R} {I_{iek}}, \quad 
\mathcal{S}_{ijk}^{s,z} = \bigcup \limits_{e=k-L}^{k+R} {I_{ije}},  
\end{equation}
where $L=L(M,s)$ and $R=R(M,s)$ are the order and stencil dependent spatial extension of the stencil to the left and to the right, 
respectively. 
Odd order schemes (even polynomial degrees $M$) always adopt three stencils, one central 
stencil ($s=1$, $L=R=M/2$), one fully left--sided stencil ($s=2$, $L=M$, $R=0$) and  one fully  right--sided stencil ($s=3$, $L=0$, $R=M$). 
Even order schemes (odd polynomial degree $M$) always adopt four stencils, two of which 
are central ($s=0$, $L=$floor$(M/2)+1$, $R=$floor$(M/2)$) and ($s=1$, $L=$floor$(M/2)$,  $R=$floor$(M/2)+1$), while the remaining two are again given by the fully left--sided and by the 
fully right--sided stencil, respectively, as defined before. The total amount of cells of the 
stencil is the same as that of the order of the scheme, namely $M+1$. 
After introducing a set of reference coordinates for each element $I_{ijk}$ given by
\begin{equation}
\label{eq:xi}
x = x_{i-\halb} + \xi   \Delta x_i, \quad 
y = x_{j-\halb} + \eta  \Delta y_j, \quad 
z = x_{k-\halb} + \zeta \Delta z_k, 
\end{equation} 
our reconstruction algorithm works according to a dimension by dimension fashion and can be described in the following steps:
%
\subsubsection{ Reconstruction in $x$ direction}
The reconstruction polynomial for each candidate stencil for element $I_{ijk}$ along the first coordinate direction $x$  
is written in terms of the basis function $\psi_l(\xi)$ as \footnote{Throughout this paper we use the Einstein summation convention, 
implying summation over indices appearing twice, although there is no need to distinguish among covariant and contra-variant indices.}  
\begin{equation}
\label{eqn.recpolydef.x} 
 \w^{s,x}_h(x,t^n) = \sum \limits_{p=0}^M \psi_p(\xi) \hat \w^{n,s}_{ijk,p} := \psi_p(\xi) \hat \w^{n,s}_{ijk,p}\,. 
\end{equation}
Integral conservation on all elements of the stencil then yields the following linear algebraic system for the unknown coefficients 
$\hat \w^{n,s}_{ijk,p}$ of the reconstruction polynomial $\w^{s,x}_h(x,t^n)$ in element $I_{ijk}$ 
\begin{equation}
 \frac{1}{\Delta x_e} \int _{x_{e-\halb}}^{x_{e+\halb}} \psi_p(\xi(x)) \hat \w^{n,s}_{ijk,p} \, dx = {\bf \bar u}^n_{ejk}, \qquad \forall {I}_{ejk} \in \mathcal{S}_{ijk}^{s,x}.      
 \label{eqn.rec.x} 
\end{equation}
For regular Cartesian meshes, the coefficients of the linear algebraic system above depend only on the choice of the basis functions, hence the system can be conveniently 
solved for the unknown coefficients $\hat \w^{n,s}_{ijk,p}$ by precomputing the inverse of the coefficient matrix, which can be done once and for all for each stencil on the 
reference element. 
%
%
Once the reconstruction has been performed for each of the stencils relative to the element $I_{ijk}$, we finally construct a data-dependent nonlinear combination of the 
polynomials obtained for each stencil, i.e. 
\begin{equation}
\label{eqn.weno} 
 \w_h^x(x,t^n) = \psi_p(\xi) \hat \w^{n}_{ijk,p}, \quad \textnormal{ with } \quad  
 \hat \w^{n}_{ijk,p} = \sum_{s=1}^{N_s} \omega_s \hat \w^{n,s}_{ijk,p},   
\end{equation}   
where the number of stencils is ${N_s}=3$ or ${N_s}=4$, depending on $M$ being even or odd, respectively. 
The nonlinear weights are given by the usual relations \cite{shu_efficient_weno} 
\begin{equation}
\omega_s = \frac{\tilde{\omega}_s}{\sum_k \tilde{\omega}_k}\,,  \qquad
\tilde{\omega}_s = \frac{\lambda_s}{\left(\sigma_s + \epsilon \right)^r}, 
\label{eqn.omegas}  
\end{equation} 
where the oscillation indicator $\sigma_s$ is
\begin{equation}
\sigma_s = \Sigma_{pm} \hat \w^{n,s}_{ijk,p} \hat \w^{n,s}_{ijk,m}\,,
\label{eqn.sigmas} 
\end{equation}
and requires the computation of the oscillation indicator matrix \cite{DumbserEnauxToro}
\begin{equation}
\Sigma_{pm} = \sum \limits_{\alpha=1}^M \int \limits_0^1 \frac{\partial^\alpha \psi_p(\xi)}{\partial \xi^\alpha} \cdot \frac{\partial^\alpha \psi_m(\xi)}{\partial \xi^\alpha} d\xi\,. 
\end{equation}
In our implementation we have adopted $\lambda_s=1$ for the one--sided stencils and $\lambda=10^5$ for the central stencils. Moreover, we use  $\epsilon=10^{-14}$ and $r=8$.

We stress that the resulting reconstruction polynomial $\w_h^x(x,t^n)$ is only a polynomial in $x$ direction, but still an average in the $y$ and $z$ direction, respectively. Hence, the reconstruction 
algorithm just described before can be applied again to the remaining two directions, as described below.  

\subsubsection{ Reconstruction in $y$ direction }
When the reconstruction algorithm is applied along the second direction $y$, the steps from (\ref{eqn.recpolydef.x})
to (\ref{eqn.weno}) are repeated for \textit{each} degree of freedom $\hat \w^{n}_{ijk,p}$. More precisely, we have 

\begin{equation}
\label{eqn.recpolydef.y} 
 \w^{s,y}_h(x,y,t^n) = \psi_p(\xi) \psi_q(\eta) \hat \w^{n,s}_{ijk,pq}\,. 
\end{equation}
Integral conservation is now applied for each degree of freedom in $x$ direction on all elements of the stencil $\mathcal{S}_{ijk}^{s,y}$ in $y$ direction, since the polynomial 
is still an average in the $y$ direction. This yields 
\begin{equation}
 \frac{1}{\Delta y_e} \int _{y_{e-\halb}}^{y_{e+\halb}} \psi_q(\eta(y)) \hat \w^{n,s}_{ijk,pq} \, dy = \hat \w^{n}_{iek,p}, \qquad \forall {I}_{iek} \in \mathcal{S}_{ijk}^{s,y}.      
 \label{eqn.rec.y} 
\end{equation}
Again, the essentially non--oscillatory property of the reconstruction polynomial is assured using the nonlinear weighting of the individual reconstruction polynomials as 
\begin{equation}
\label{eqn.weno.y} 
 \w_h^y(x,y,t^n) = \psi_p(\xi) \psi_q(\eta) \hat \w^{n}_{ijk,pq}, \quad \textnormal{ with } \quad  
 \hat \w^{n}_{ijk,pq} = \sum_{s=1}^{N_s} \omega_s \hat \w^{n,s}_{ijk,pq},   
\end{equation}   
with $\omega_s$ and $\sigma_s$ defined as above in Eqn. \eqref{eqn.omegas} and \eqref{eqn.sigmas}. 

\subsubsection{ Reconstruction in $z$ direction}
Finally, when the reconstruction algorithm is applied along the last direction $z$, the steps from (\ref{eqn.recpolydef.x}) 
to (\ref{eqn.weno}) are repeated for the $(M+1)^2$ degrees of freedom of the polynomial already reconstructed along $x$ and $y$. 
We therefore have 
\begin{equation}
\label{eqn.recpolydef.z} 
 \w^{s,z}_h(x,y,z,t^n) = \psi_p(\xi) \psi_q(\eta) \psi_r(\zeta) \hat \w^{n,s}_{ijk,pqr}\,. 
\end{equation}
with the integral conservation written as above,
\begin{equation}
 \frac{1}{\Delta z_e} \int _{z_{e-\halb}}^{z_{e+\halb}} \psi_r(\zeta(z)) \hat \w^{n,s}_{ijk,pqr} \, dz = \hat \w^{n}_{iek,pq}, \qquad \forall {I}_{ije} \in \mathcal{S}_{ijk}^{s,z}.      
 \label{eqn.rec.z} 
\end{equation}
The final three--dimensional WENO polynomial is then given by 
\begin{equation}
 \w_h(\mathbf{x},t^n) = \psi_p(\xi) \psi_q(\eta) \psi_r(\zeta) \hat \w^{n}_{ijk,pqr},  
\label{eqn.weno.z} 
\end{equation}   
with
\begin{equation}
 \hat \w^{n}_{ijk,pqr} = \sum_{s=1}^{N_s} \omega_s \hat \w^{n,s}_{ijk,pqr}.    
\end{equation}
We stress that the resulting polynomial WENO reconstruction produces \textit{entire polynomials}, 
although they are given under the form of a \textit{nodal basis}, but this is just a technical issue. 
Any other set of basis functions could have been chosen as well. This is the main difference with 
respect to the original optimal WENO scheme of Jiang \& Shu \cite{shu_efficient_weno}, which 
produces \textit{point values} at the element interfaces. The drawback of our present method is of 
course that the method is not optimal, in the sense that to obtain a given order of accuracy the
total stencil of the present method is larger than the one of the optimal WENO. 

\subsection{Local space--time DG predictor}
\label{sec:localDG}

Once a high order polynomial in space $\w_h$ has been reconstructed for each cell, it is necessary to evolve it in time in order 
to compute the fluxes and sources according to Eqs.~\eqref{flux:F}--\eqref{source:S}. Recall that the arguments of the numerical
fluxes \eqref{flux:F}-\eqref{flux:H} are the yet undefined space--time polynomials $\q_h$. The strategy of a high order time
evolution follows the spirit of the original MUSCL scheme of van Leer \cite{leer5} as well as the one of the ENO scheme \cite{eno}. 
It can furthermore also be interpreted as an approximate solution of the generalized Riemann problem at the cell interfaces as done
in the ADER approach \cite{titarevtoro,titarevtoro2}, see \cite{CastroToro,Montecinos}. 
However, the MUSCL scheme, the ENO method as well as the ADER approach obtain all the higher order in time through the use of Taylor 
series and a (repeated) use of the governing conservation law in strong differential form, which for higher than second order and for 
complicated PDE systems may quickly become very cumbersome or even unfeasible. Starting from the work by Dumbser et al.~\cite{DumbserEnauxToro,Dumbser2008} 
an alternative time evolution has been proposed, which is able to deal with general nonlinear conservation laws, thus providing a more 
flexible one--step scheme compared to the ones based on Taylor expansions. The new method relies on a \textit{weak} integral formulation of 
the governing PDE in space--time using an element--local space--time Galerkin predictor method. For this time evolution strategy, no algebraic 
manipulations of space and time derivatives are necessary, but all that is required is a point--wise evaluation of fluxes and source terms. 
The result of the local space--time Galerkin predictor are the high order space--time polynomials $\q_h$ that are needed for the evaluation of 
the numerical fluxes according to Eqs.~\eqref{flux:F}--\eqref{source:S}. In the following we briefly illustrate the method, addressing to 
\cite{DumbserEnauxToro,Dumbser2008,DumbserZanotti,HidalgoDumbser} for more details. 

In addition to the spatial reference coordinates $\boldsymbol{\xi}=(\xi,\eta,\zeta)$ already defined by 
Eqn.~\eqref{eq:xi}, we also introduce a reference time coordinate $\tau$ as $t = t^n + \tau \Delta t$, 
spanning the unit interval $[0,1]$. As a result, the governing PDE \eqref{eq:governingPDE} can be rewritten in compact form as 
\begin{equation}
\frac{\partial{\bf u}}{\partial \tau} + \frac{\partial \mathbf{f}^\ast}{\partial \xi} + \frac{\partial \mathbf{g}^\ast}{\partial \eta} + \frac{\partial \mathbf{h}^\ast}{\partial \zeta} ={\bf S}^\ast 
\label{eqn.pde.ref} 
\end{equation}
with
\begin{equation}
{\bf f}^\ast= \frac{\Delta t}{\Delta x_i} \, {\bf f}, \quad 
{\bf g}^\ast= \frac{\Delta t}{\Delta y_j} \, {\bf g}, \quad 
{\bf h}^\ast= \frac{\Delta t}{\Delta z_k} \, {\bf h}, \quad 
{\bf S}^\ast= \Delta t {\bf S}. 
\end{equation}
We then introduce the space--time basis functions $\theta_\mathfrak{p}(\boldsymbol{\xi},\tau)$, which are piecewise space--time polynomials
of degree $M$, using the multi--index $\mathfrak{p}=(p,q,r,s)$. The $\theta_\mathfrak{p}$ are given by a tensor--product of the basis functions 
$\psi_l$ already used before in the reconstruction procedure, hence 
\begin{equation}
  \theta_\mathfrak{p}(\boldsymbol{\xi},\tau) = \psi_p(\xi) \psi_q(\eta) \psi_r(\zeta) \psi_s(\tau).  
\end{equation} 
Multiplication of the PDE \eqref{eqn.pde.ref} with the space--time test functions $\theta_\mathfrak{q}$ and integrating over
the space--time reference control volume $[0;1]^4$ yields 
\begin{equation}
 \int \limits_{0}^{1} \int \limits_{0}^{1}  \int \limits_{0}^{1}   \int \limits_{0}^{1}   
\theta_\mathfrak{q} \left( 
  \frac{\partial{\bf u}}{\partial \tau} + \frac{\partial \mathbf{f}^\ast}{\partial \xi} + \frac{\partial \mathbf{g}^\ast}{\partial \eta} + \frac{\partial \mathbf{h}^\ast}{\partial \zeta} - {\bf S}^\ast \right) d\xi d\eta d\zeta d\tau = 0.  
\label{eqn.pde.weak1} 
\end{equation}
Integration of the first term that contains the time derivative by parts yields 
\vspace{-8mm}
\begin{eqnarray}
 && \int \limits_{0}^{1} \int \limits_{0}^{1}  \int \limits_{0}^{1} \theta_\mathfrak{q}(\boldsymbol{\xi},1) \u(\boldsymbol{\xi},1) d\xi d\eta d\zeta - 
  \int \limits_{0}^{1} \int \limits_{0}^{1}  \int \limits_{0}^{1}   \int \limits_{0}^{1} \left( \frac{\partial}{\partial \tau} \theta_\mathfrak{q} \right) \u d\xi d\eta d\zeta d\tau   \nonumber \\ 
 && + \int \limits_{0}^{1} \int \limits_{0}^{1}  \int \limits_{0}^{1}   \int \limits_{0}^{1} \left[   
    \theta_\mathfrak{q} \left( 
   \frac{\partial \mathbf{f}^\ast}{\partial \xi} + \frac{\partial \mathbf{g}^\ast}{\partial \eta} + \frac{\partial \mathbf{h}^\ast}{\partial \zeta} - {\bf S}^\ast \right) \right] d\xi d\eta d\zeta d\tau 
   \nonumber \\ 
 &&  = \int \limits_{0}^{1} \int \limits_{0}^{1}  \int \limits_{0}^{1} \theta_\mathfrak{q}(\boldsymbol{\xi},0) \w_h(\boldsymbol{\xi},t^n) d\xi d\eta d\zeta.  
\label{eqn.pde.weak2} 
\end{eqnarray}
In the following, we denote the discrete space--time solution by $\mathbf{q}_h$, for which we make the following ansatz 
\begin{equation}
 \mathbf{q}_h = \mathbf{q}_h(\boldsymbol{\xi},\tau) = \theta_\mathfrak{p}\left(\boldsymbol{\xi},\tau \right) \hat \q_\mathfrak{p},   
 \label{eqn.st.q} 
\end{equation}
with the yet unknown degrees of freedom of the space--time polynomial $\hat \q_\mathfrak{p} = \hat \q_{pqrs}$ 
We use the same representation for the fluxes and source terms, hence
\begin{equation}
 \mathbf{f}^{\ast}_h = \theta_\mathfrak{p} \hat{\mathbf{f}}^{\ast}_\mathfrak{p},   \qquad 
 \mathbf{g}^{\ast}_h = \theta_\mathfrak{p} \hat{\mathbf{g}}^{\ast}_\mathfrak{p},   \qquad 
 \mathbf{h}^{\ast}_h = \theta_\mathfrak{p} \hat{\mathbf{h}}^{\ast}_\mathfrak{p},   \qquad 
 \mathbf{S}^{\ast}_h = \theta_\mathfrak{p} \hat{\mathbf{S}}^{\ast}_\mathfrak{p}.  
 \label{eqn.st.fs} 
\end{equation}
Due to the nodal approach, the above degrees of freedom for the fluxes and source terms are simply the point--wise evaluation of the physical fluxes and source terms, hence 
\begin{equation}
 \hat{\mathbf{f}}^{\ast}_\mathfrak{p} = {\mathbf{f}}^{\ast}\left( \hat \q_\mathfrak{p} \right), \qquad 
 \hat{\mathbf{g}}^{\ast}_\mathfrak{p} = {\mathbf{g}}^{\ast}\left( \hat \q_\mathfrak{p} \right), \qquad 
 \hat{\mathbf{h}}^{\ast}_\mathfrak{p} = {\mathbf{h}}^{\ast}\left( \hat \q_\mathfrak{p} \right), \qquad 
 \hat{\mathbf{S}}^{\ast}_\mathfrak{p} = {\mathbf{S}}^{\ast}\left( \hat \q_\mathfrak{p} \right). 
 \label{eqn.nodal.eval} 
\end{equation}
Inserting Eqns. \eqref{eqn.st.q} and \eqref{eqn.st.fs} into \eqref{eqn.pde.weak2} yields 
\vspace{-9mm} 
\begin{eqnarray}
 && \int \limits_{0}^{1} \int \limits_{0}^{1}  \int \limits_{0}^{1} \theta_\mathfrak{q}(\boldsymbol{\xi},1) \theta_\mathfrak{p}(\boldsymbol{\xi},1) \hat \q_\mathfrak{p}
  \, d\xi d\eta d\zeta 
 - \int \limits_{0}^{1} \int \limits_{0}^{1}  \int \limits_{0}^{1}   \int \limits_{0}^{1} \left(\frac{\partial}{\partial \tau} \theta_\mathfrak{q} \right) \theta_\mathfrak{p} \hat \q_\mathfrak{p}  
 \, d\xi d\eta d\zeta d\tau 
   \nonumber \\ 
 && + \int \limits_{0}^{1} \int \limits_{0}^{1}  \int \limits_{0}^{1}   \int \limits_{0}^{1} \left[   
     \theta_\mathfrak{q} \left( 
     \frac{\partial}{\partial \xi}   \theta_\mathfrak{p} \hat{\mathbf{f}}^{\ast}_\mathfrak{p} 
   + \frac{\partial}{\partial \eta}  \theta_\mathfrak{p} \hat{\mathbf{g}}^{\ast}_\mathfrak{p} 
   + \frac{\partial}{\partial \zeta} \theta_\mathfrak{p} \hat{\mathbf{h}}^{\ast}_\mathfrak{p} 
   - \theta_\mathfrak{p} \hat{\mathbf{S}}^{\ast}_\mathfrak{p}  \right) \right] \, d\xi d\eta d\zeta d\tau 
   \nonumber \\ 
   && = \int \limits_{0}^{1} \int \limits_{0}^{1}  \int \limits_{0}^{1} \theta_\mathfrak{q}(\boldsymbol{\xi},0) \w_h(\boldsymbol{\xi},t^n) \, d\xi d\eta d\zeta.  
\label{eqn.pde.weak3} 
\end{eqnarray}
The above weak form is an element--local nonlinear algebraic equation system for the unknown coefficients $\hat \q_\mathfrak{p}$. The initial condition is 
included in a weak sense by the integral on the right hand side, where the reconstructed solution $\w_h(\boldsymbol{\xi},t^n)$ is given by Eqn. 
\eqref{eqn.weno.z}. 
After introducing the integrals 
\begin{equation}
 \mathbf{K}^1_{\mathfrak{q} \mathfrak{p}} = \int \limits_{0}^{1} \int \limits_{0}^{1}  \int \limits_{0}^{1} \theta_\mathfrak{q}(\boldsymbol{\xi},1) \theta_\mathfrak{p}(\boldsymbol{\xi},1) d \boldsymbol{\xi} - 
                \int \limits_{0}^{1} \int \limits_{0}^{1}  \int \limits_{0}^{1}   \int \limits_{0}^{1} \left(\frac{\partial}{\partial \tau} \theta_\mathfrak{q} \right) \theta_\mathfrak{p} d \boldsymbol{\xi} d\tau,  
\end{equation} 
\begin{equation}
 \mathbf{K}^{\boldsymbol{\xi}}_{\mathfrak{q} \mathfrak{p}} = \left( \mathbf{K}^{\xi}_{\mathfrak{q} \mathfrak{p}} , \mathbf{K}^{\eta}_{\mathfrak{q} \mathfrak{p}}, \mathbf{K}^{\zeta}_{\mathfrak{q} \mathfrak{p}} \right) =  \int \limits_{0}^{1} \int \limits_{0}^{1}  \int \limits_{0}^{1} \int \limits_{0}^{1} \theta_\mathfrak{q}  
     \frac{\partial}{\partial \boldsymbol{\xi}}  \theta_\mathfrak{p} d \boldsymbol{\xi} d\tau,  
\end{equation} 
\begin{equation}
 \mathbf{M}_{\mathfrak{q} \mathfrak{p}} =  \int \limits_{0}^{1} \int \limits_{0}^{1}  \int \limits_{0}^{1} \int \limits_{0}^{1} \theta_\mathfrak{q} \theta_\mathfrak{p} d \boldsymbol{\xi} d\tau,  
\end{equation} 
and
\begin{equation}
 \mathbf{F}^0_{\mathfrak{q} \mathfrak{p}} = \int \limits_{0}^{1} \int \limits_{0}^{1}  \int \limits_{0}^{1} \theta_\mathfrak{q}(\boldsymbol{\xi},0) \psi_\mathfrak{m}(\boldsymbol{\xi}) d \boldsymbol{\xi},  
\end{equation} 
where $d \boldsymbol{\xi} = d\xi d\eta d\zeta$ one can rewrite the above system in compact matrix--vector form as 
\begin{equation}
 \mathbf{K}^1_{\mathfrak{q} \mathfrak{p}}   \hat \q_\mathfrak{p} + 
 \mathbf{K}^\xi  _{\mathfrak{q} \mathfrak{p}} \cdot \hat{\mathbf{f}}^\ast_\mathfrak{p} + 
 \mathbf{K}^\eta _{\mathfrak{q} \mathfrak{p}} \hat{\mathbf{g}}^\ast_\mathfrak{p} + 
 \mathbf{K}^\zeta_{\mathfrak{q} \mathfrak{p}} \hat{\mathbf{h}}^\ast_\mathfrak{p} = 
 \mathbf{M}_{\mathfrak{q} \mathfrak{p}} \hat{\mathbf{S}}^\ast_\mathfrak{p} + \mathbf{F}^0_{\mathfrak{q} \mathfrak{m}} \hat{\mathbf{w}}_{\mathfrak{m}}^n, 
\label{eqn.pde.weak4} 
\end{equation} 
with the spatial multi--index $\mathfrak{m}=(k,l,m)$ and $\psi_\mathfrak{m}(\boldsymbol{\xi}) = \psi_k(\xi) \psi_l(\eta) \psi_m(\zeta)$.   
Due to the tensor--product nature of the basis functions and the control volumes, the above matrices are all very sparse block--matrices, where all sub--blocks 
contain purely one--dimensional integrals. Hence, the product of the matrices with the vectors of degrees of freedom can be efficiently implemented in a 
dimension--by--dimension manner. 
Equation \eqref{eqn.pde.weak4} is conveniently solved by the following iterative scheme, introduced in \cite{Dumbser2008,DumbserZanotti}, 
\begin{equation}
 \mathbf{K}^1_{\mathfrak{q} \mathfrak{p}}   \hat \q_\mathfrak{p}^{k+1} - \mathbf{M}_{\mathfrak{q} \mathfrak{p}} \hat{\mathbf{S}}^{\ast,k+1}_\mathfrak{p} = 
  \mathbf{F}^0_{\mathfrak{q} \mathfrak{m}} \hat{\mathbf{w}}_{\mathfrak{m}}^n  
 -\mathbf{K}^\xi  _{\mathfrak{q} \mathfrak{p}} \cdot \hat{\mathbf{f}}^{\ast,k}_\mathfrak{p}   
 -\mathbf{K}^\eta _{\mathfrak{q} \mathfrak{p}} \hat{\mathbf{g}}^{\ast,k}_\mathfrak{p}  
 -\mathbf{K}^\zeta_{\mathfrak{q} \mathfrak{p}} \hat{\mathbf{h}}^{\ast,k}_\mathfrak{p}    
\label{eqn.stdg.final} 
\end{equation} 
for which an efficient second--order MUSCL--type initial guess has been proposed in \cite{HidalgoDumbser}. Moreover, the matrix 
$\mathbf{K}^1_{\mathfrak{q} \mathfrak{p}}$ can be easily inverted once and for all on the reference element by inverting the 
one--dimensional sub--blocks in time that are associated with each Gaussian quadrature point in space. 

\section{Adaptive Mesh Refinement}
\label{sec:AMR}

There are two major strategies for implementing an AMR algorithm. The first one adopts nested arrays of logically rectangular grid patches, according to the 
original Berger-Colella-Oliger approach~\cite{Berger-Oliger1984,Berger-Jameson1985,Berger-Colella1989}. The second strategy is referred to as the 'cell-by-cell'  
refinement, the implementation of which is somehow similar to the one of finite volume schemes on unstructured meshes~\cite{Khokhlov1998}. In our work we have 
followed the second choice, since it can be easily implemented using a tree-type data structure and is slightly more general than the first one. 

As already mentioned at the beginning of the previous section, the three ingredients of our high order one--step finite volume schemes (reconstruction, time evolution, 
finite volume update) can be transferred in a rather straightforward manner to space--time adaptive grids. 
In particular, the element--local space--time DG predictor on regular Cartesian grids is identical to the one on AMR grids, since for the predictor there is no need to exchange 
information with neighbor elements. Hence, even if two adjacent cells are on different levels of grid refinement this will not alter the local space--time DG predictor scheme
at all. The reconstruction procedure obviously needs cell averages from neighbor elements. Therefore, in our AMR implementation the case of adjacent cells with different levels 
of refinement is handled in such a way that a real cell of grid level $\ell$ is \textit{always} surrounded by a sufficiently thick layer of virtual (ghost) cells on the same layer. 
The cell averages of the virtual cells are obtained from the real cells on the same 
location either by averaging or projection, depending whether the virtual cell is on 
a coarser or finer level of refinement. 
Since for high order WENO schemes the total effective stencil needed for obtaining the reconstruction polynomial within a spatial control volume $I_{ijk}$ grows with the degree
of the reconstruction polynomial $M$, the layer of ghost cells must always have at least a thickness of $M$ ghost cells. In order to keep the information needed for reconstruction \textit{local} on the 
\textit{coarser} grid level, the grid refinement factor $\mathfrak{r}$ adopted in the scheme [see definition \eqref{refine-factor} below]
must 
satisfy $\mathfrak{r} \ge M$. The fact that active cells are surrounded by a layer of ghost 
cells can of course also be interpreted as use of \textit{micro--patches}. 
The virtual ghost cells are logically connected neighbors of the real cells in order to allow the uniform Cartesian reconstruction procedure to be performed exactly as described 
in the previous section. The only real difference between uniform Cartesian grid and AMR grid can be found in the computation of the fluxes at edges with two adjacent cells of different level 
of refinement. To keep the algorithmic complexity as low as possible, in our approach the 
level of refinement of two adjacent grid cells may differ by at most one. 

\subsection{AMR implementation}
\label{AMR-implementation}
%
We have developed a cell-by-cell AMR technique in which the computational domain is 
discretized with a uniform Cartesian grid at the coarsest level. 
We use ${\mathcal L}_0$ to denote this initial grid on the coarsest level of refinement  ($\ell=0$), while ${\mathcal L}_{\ell}$ indicates the union of all elements 
up to level $\ell$. 
Already at time $t=0$, the refinement criterion (see Sect.~\ref{Ref-criterion}) is applied to  the initial condition, thus producing a hierarchy of refinement levels up to a prescribed  maximum level of refinement $\ell_{\rm max}$, i.e. $0\leq \ell \leq \ell_{\rm max}$. In the 
rest of this section we will use the following terminology: 
With \textit{children} we intend the cells on the next refinement level $\ell+1$ contained 
in a cell ${\mathcal C}_m$ of level $\ell$ after its refinement. For the children cells the  original cell ${\mathcal C}_m$ is denoted as their \textit{mother} cell. The \textit{Neumann
neighbors} $\mathcal{N}_m$ of a cell ${\mathcal C}_m$ are the neighbor cells that share a 
common face with cell ${\mathcal C}_m$. There are $2d$ Neumann neighbors in $d$ space  dimensions. For example, the numerical fluxes \eqref{flux:F}-\eqref{flux:H} are evaluated between Neumann neighbors. 
The \textit{Voronoi neighbors} $\mathcal{V}_m$ of a cell ${\mathcal C}_m$ are cells which 
share common nodes, hence the Neumann neighbors are a subset of the Voronoi neighbors. There 
are $3^d-1$ Voronoi neighbors in $d$ space dimensions. 

\begin{figure}
\begin{center}
\includegraphics[angle=0,width=0.40\textwidth]{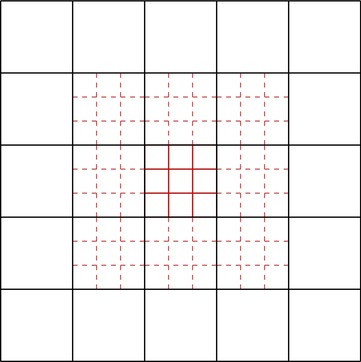}
\hspace{5mm}
\includegraphics[angle=0,width=0.40\textwidth]{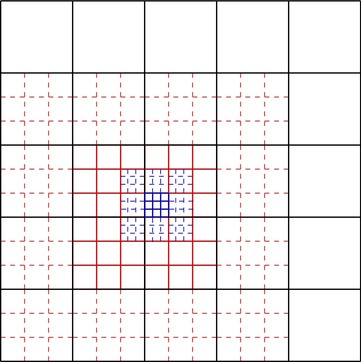}
\caption{Refinement mechanism involving two levels of refinement ($\ell \leq 2$) in a 2D geometry. Solid lines denote active cells, while dashed lines denote virtual cells.  Left panel: when the central cell ${\mathcal C}_{j}$ at level $\ell = 0$ is 
refined (solid-dark-red line), all its Voronoi neighbors $\mathcal{V}_j$ are virtually 
refined (dashed-light-red lines). 
Right panel: if the cell ${\mathcal C}_{k}$ at level $\ell=1$ on the left corner of ${\mathcal C}_{j}$ requires further refinement to level $\ell=2$ (solid-dark-blue lines), the virtual  children of the lower, left and lower left neighbors of ${\mathcal C}_{j}$ must be activated  
(to avoid violation of rule (5)), hence all their Voronoi neighbors without children must be virtually refined to level $\ell=1$. Furthermore, the neighbors of ${\mathcal C}_{k}$ must be virtually refined to level $\ell=2$  (dashed-light-blue lines). 
}
\label{fig_three_levels}
\end{center}
\end{figure}

All over the simulation, the following general rules have to be met: 
\begin{enumerate}

\item 
Whenever a cell of the level $\ell$ is refined, it is subdivided into an integer number $\mathfrak{r}$ of finer cells along {\em each direction}, such that 
\begin{equation}
\label{refine-factor}
\Delta x_{\ell} = \mathfrak{r} \Delta x_{\ell+1}\, \quad \Delta y_{\ell} = \mathfrak{r} \Delta y_{\ell+1} \, \quad 
\Delta z_{\ell} = \mathfrak{r} \Delta z_{\ell+1}, 
\end{equation}
and also the time steps are chosen \textit{locally} on each level so that 
\begin{equation}
\Delta t_{\ell} = \mathfrak{r} \Delta t_{\ell+1}.  
\end{equation}
As a result, each {\em mother} cell generates $\mathfrak{r}^d$ {\em children} cells in $d$ 
space dimensions.
%
Since the high order WENO reconstruction needs information from more cells than just the direct 
neighbors, the following condition must always hold to keep reconstruction \textit{local} 
on the coarser grid level: 
\begin{equation}
 \mathfrak{r} \ge M. 
\end{equation}

\item 
Each cell ${\mathcal C}_m$, at any level of refinement, has one among three possible 
\textit{status} flags. The status flag is denoted by $\sigma$ in the following. A cell is 
either an \textit{active cell} ($\sigma=0$), and therefore has to be updated 
through the finite volume scheme; or, it is a \textit{virtual child cell} (with status $\sigma=1$) and is updated by projection of the mother's high order space--time polynomial $\q_h$; or, finally, a \textit{virtual mother cell} ($\sigma=-1$), updated by recursively  averaging over all children from higher refinement levels.  
A virtual child cell has always an active mother cell, while a virtual mother cell has 
children with status $\sigma \leq 0$. A virtual child cell cannot be further refined, 
unless it is first activated.  

\item Any cell ${\mathcal C}_m$, at any level of refinement and with any status, is identified with a unique positive\footnote{ See Sect.~\ref{sec:AMR-parallelization} for 
the possibility of negative integer numbers assigned to those cells 
belonging to the ghost zone at the MPI border between two processors.}
integer number $m$, with $1\leq m \leq N_{{\rm Cells}}$ and 
$N_{{\rm Cells}}$ being the total number of cells at any given time. $N_{{\rm Cells}}$ is of course a time-dependent quantity, increasing for any refinement operation, and decreasing for any recoarsening operation.

\item
A refined cell with children of status $\sigma \leq 0$ cannot have any Voronoi neighbor  
without children. Therefore, as soon as a cell ${\mathcal C}_m$ is refined and the 
generated children have status $\sigma=0$, or as soon as the existing virtual children of cell 
${\mathcal C}_m$ are activated, all Voronoi neighbors without children 
are virtually refined, i.e. they generate virtual children with status $\sigma=1$. 
The status of the regularly refined cell ${\mathcal C}_m$ changes to $\sigma = -1$. 

\item
The levels of refinement of two cells that are Voronoi neighbors of each other can only 
differ by at most unity. Violation of this rule is avoided through appropriate 
activation of virtual children. Such a situation is schematically depicted in  Fig.~\ref{fig_three_levels}. 

\item
The maximum level of refinement is a prescribed level $\ell_{\max}$. Cells are not 
allowed to be refined beyond this level. 

\item When a cell ${\mathcal C}_m$ is recoarsened, its children cells are destroyed. 
However, this is not allowed if the cell ${\mathcal C}_m$ contains children that contain themselves children of any status. Furthermore, if the cell ${\mathcal C}_m$ is to be 
recoarsened and has a neighbor cell with active children, then the cell ${\mathcal C}_m$ 
can only deactivate its children cells, changing them to virtual, but cannot destroy them. 

\end{enumerate}
For each cell we store the indices of the Voronoi neighbors, as well as pointers to the 
mother and the first child of a cell. The indices of all $\mathfrak{r}^d$ children are
obtained by the convention that all children of a cell have consecutive numbers. 
We furthermore store the status $\sigma$ of all elements. 

\begin{figure}
\begin{center}
\vspace{5mm} 
\includegraphics[angle=0,width=0.9\textwidth]{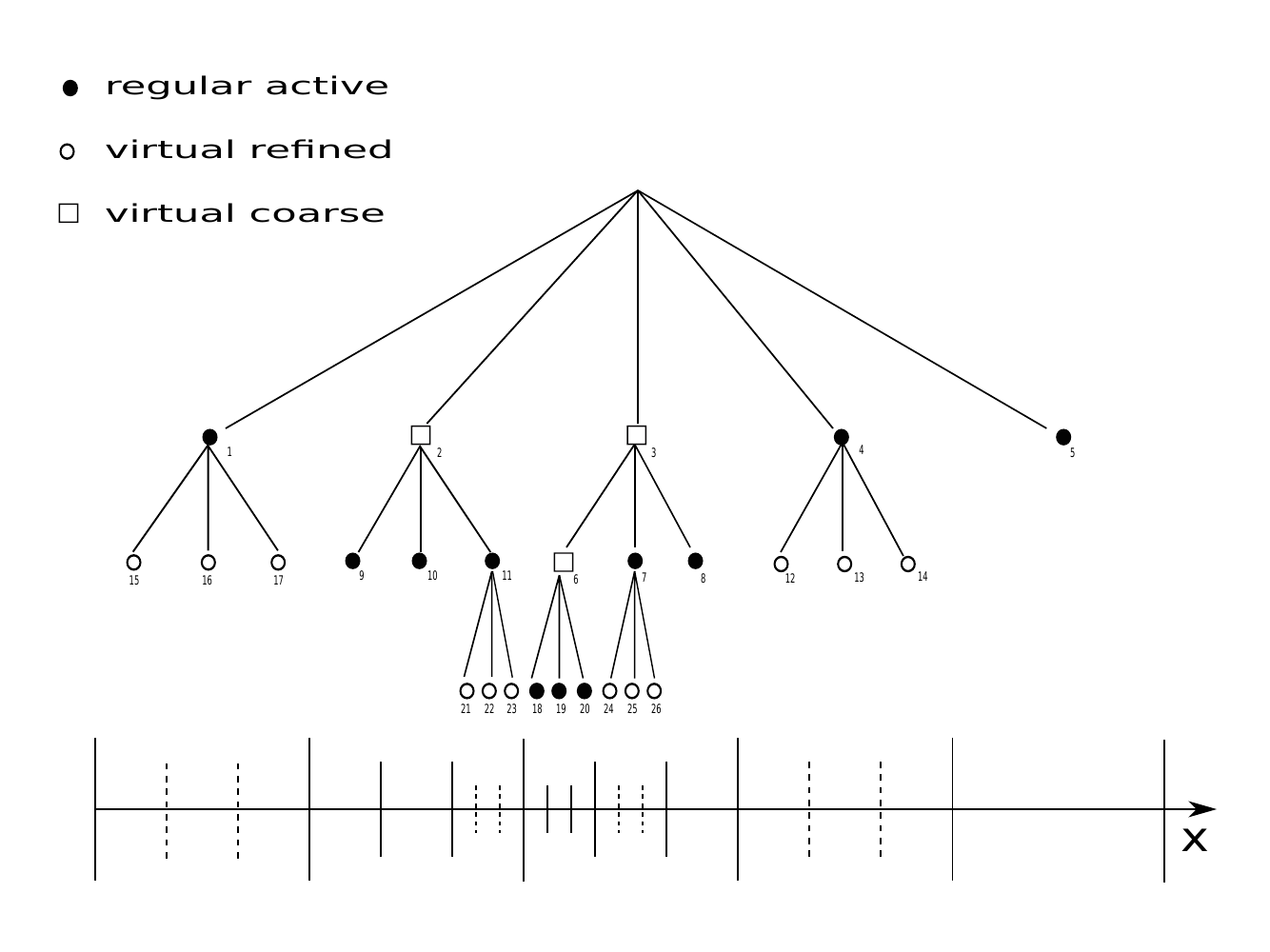}
\caption{
One--dimensional example for the tree structure used in the present cell-by-cell 
implementation of AMR with $\mathfrak{r}=3$ and $\ell \leq 2$. 
}
\label{fig.tree}
\end{center}
\end{figure}

A simple 1D example of the tree--structure used in our algorithm together with an 
illustration of the use of the status flag $\sigma$ is depicted in Fig. \ref{fig.tree}. 

\subsubsection{Refinement criterion}
\label{Ref-criterion}

As a refinement criterion, we have adopted the same strategy described in 
\cite{Lohner1987} and later adopted by \cite{Fryxell2000} and \cite{Mignone2012}. 
Such a criterion has two very desirable properties. The first one is that it is entirely local, 
thus avoiding any global computation. The second one is that it is based on the calculation of 
a second derivative error, thus avoiding unnecessary refinement in smooth regions, such as the wake of a rarefaction wave.
In practice, a cell ${\mathcal C}_m$ is marked for regular refinement (generation of 
children with status $\sigma=0$) if $\chi_m>\chi_{\rm ref}$, while it is marked for 
recoarsening if $\chi_m<\chi_{\rm rec}$, where
\begin{equation}
\chi_m=\sqrt{\frac{\sum_{k,l} (\partial^2 \Phi/\partial x_k \partial x_l)^2 }{\sum_{k,l}[(|\partial \Phi/\partial x_k|_{i+1}+|\partial \Phi/\partial x_k|_i)/\Delta x_l+\varepsilon|(\partial^2 /\partial x_k \partial x_l )||\Phi|]^2} }\,.
\label{eqn.indicator}
\end{equation}
The summation $\sum_{k,l}$ is taken over the number of space dimension of the problem in order to include the cross term derivatives. We stress that $\Phi=\Phi(\u)$ can be any suitable indicator 
function of the conservative variables $\u$. We have observed that the threshold values   $\chi_{\rm ref}$ and $\chi_{\rm rec}$ can be slightly model dependent. In most of our tests 
we have chosen $\chi_{\rm ref}$ in the range $\sim[0.2,0.25]$ and $\chi_{\rm rec}$ in the 
range $\sim[0.05,0.15]$. Finally, the parameter $\varepsilon$ acts as a filter preventing  refinement in regions of small  ripples and is given the value $\varepsilon = 0.01$.

\subsubsection{Time accurate local timestepping}
\label{sec:time}
Each of the refinement levels is advanced in time with its own \textit{local} time-step 
$\Delta t_\ell = \mathfrak{r} \Delta t_{\ell+1}$, as already defined above. The use of
time steps that are integer multiples of each other for each level is very convenient, 
but not strictly necessary, see e.g. the local time stepping schemes presented in 
\cite{dumbserkaeser06d,TaubeMaxwell,stedg1}, where a different local time step is 
allowed for each element. 
Let us denote by $t^n_\ell$ and $t^{n+1}_\ell$ the current and future times of the level $\ell$. 
Then, the level scheduled for updating is the largest value of $\ell$ that satisfies the  \textit{update criterion} \cite{dumbserkaeser06d}  
\begin{equation} 
 t^{n+1}_{\ell} \leq t^{n+1}_{\ell-1},  \qquad  0 \leq \ell \leq \ell_{\max}, 
 \label{eqn.update.criterion}  
\end{equation} 
where we define $t^{n+1}_{-1}:=t^{n+1}_{0}$ for convenience, so that also the scheduling of 
level $\ell=0$ is included in Eqn. \eqref{eqn.update.criterion}.  

In other words, starting from the common initial time $t=0$, the finest level of refinement  $\ell_{\max}$ is evolved first and performs a number of $\mathfrak{r}$ sub-timesteps before 
the next coarser level $\ell_{\max}-1$ performs its first time update. This procedure is then applied recursively and it implies a total amount of $\mathfrak{r}^\ell$ sub-timesteps on 
each level to be performed in order to reach the time $t_0^{n+1}$ of the coarsest level. 

The computation of numerical fluxes between two adjacent cells on different levels of 
refinement is rather straightforward thanks to the use of the local space--time predictor, 
which computes the predictor solution $\q_h$ for each element after reconstruction and 
which is valid from time $t^n_\ell$ to time $t^{n+1}_\ell$. In this way, and with the update 
criterion \eqref{eqn.update.criterion}, a high order space--time polynomial is always 
available on both sides of an element interface for flux computation.  

\begin{figure}[!htbp] 
\begin{center}
\includegraphics[angle=0,width=0.95\textwidth]{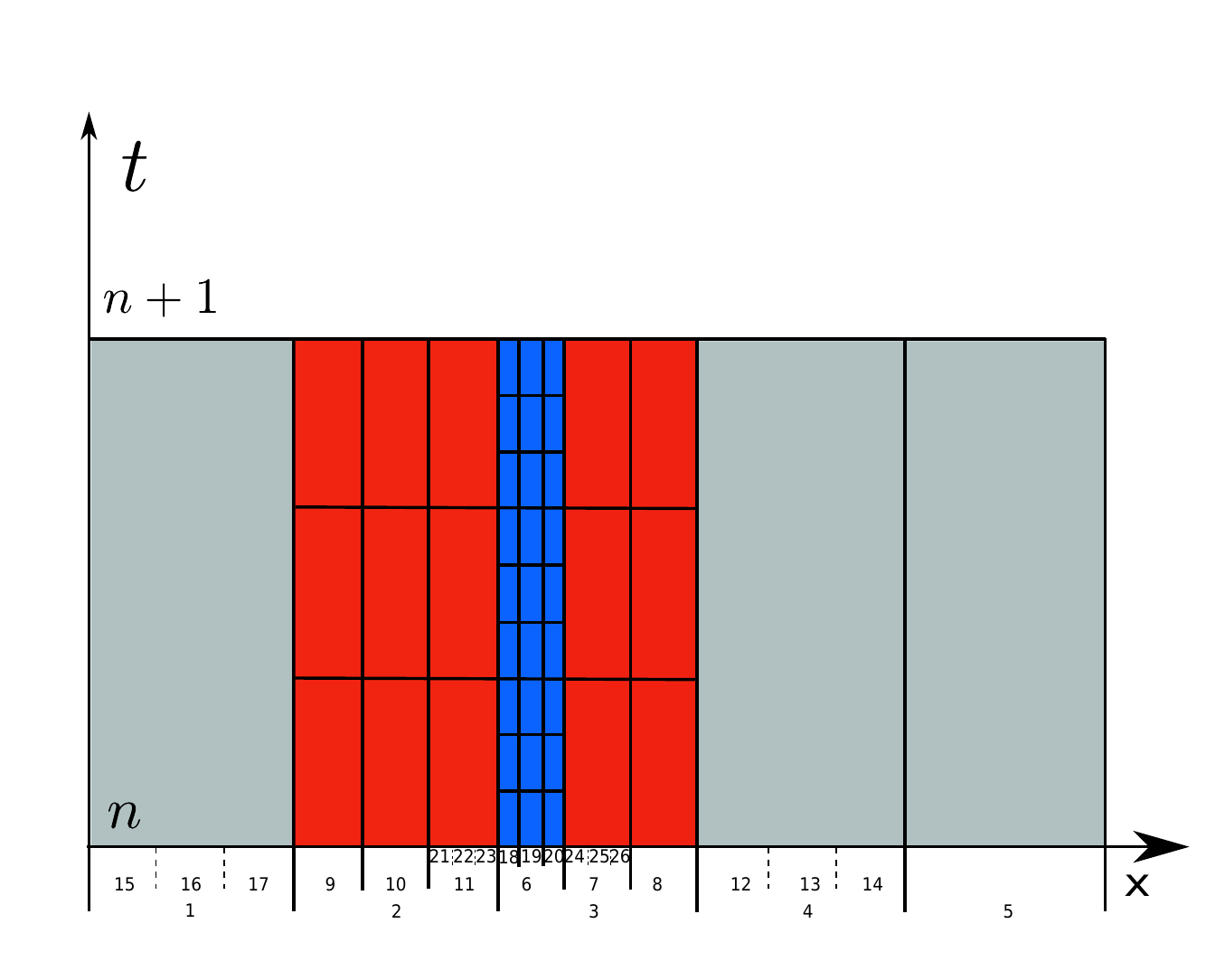}
\caption{
One--dimensional example for the local timestepping with $\mathfrak{r}=3$ and $\ell \leq 2$. 
The grid and the element numbers correspond to the ones depicted in Fig. \ref{fig.tree}. 
}
\label{fig.time.schedule}
\end{center}
\end{figure}

To illustrate the  procedure, we use the one--dimensional example depicted in Fig.  \ref{fig.time.schedule}, which corresponds to the case already shown in Fig. \ref{fig.tree} before. 
All elements start from a common time level $t^n_0$ and the space--time predictor solution 
$\q_h$ has been computed in all elements. Then, the first level that satisfies the update  criterion \eqref{eqn.update.criterion} is $\ell=2$. Computing the fluxes for cell  $\mathcal{C}_{19}$ is exactly 
as on a uniform Cartesian mesh, since $\mathcal{C}_{19}$ has two active neighbors on the same grid level.  The situation is different for $\mathcal{C}_{18}$ because its left real (active) neighbor is cell  $\mathcal{C}_{11}$, which is on a coarser grid level. Nevertheless, the numerical flux between  $\mathcal{C}_{18}$ and $\mathcal{C}_{11}$ can be computed since the predictor solution $\q_h$ is available 
for the necessary time interval $[t^n_2;t^n_2 + \Delta t_2]$ in both cells. To make the approach  \textit{conservative},  the computed numerical flux between $\mathcal{C}_{11}$ and $\mathcal{C}_{18}$ is 
used to update $\mathcal{C}_{18}$, but it is at the same time also stored in a memory variable of the  virtual neighbor cell $\mathcal{C}_{23}$ and will be used later to update the real element  $\mathcal{C}_{11}$. For convenience, the virtual neighbor cell $\mathcal{C}_{23}$ can also be used 
to technically handle the hanging node in time (and in space for $d>1$) by scaling the space--time polynomial $\q_h$ of the coarse element $\mathcal{C}_{11}$ down to the finer grid level $\ell=2$. 
A similar situation occurs for element $\mathcal{C}_{20}$ and its right active neighbor 
$\mathcal{C}_{7}$. After the first time step of level $\ell=2$, the elements  
$\mathcal{C}_{18}-\mathcal{C}_{20}$ are at time $t^n_0 + \Delta t_2$, while all other elements
are still at time $t^n_0$. To perform reconstruction for cells $\mathcal{C}_{18}-\mathcal{C}_{20}$, 
the predictor solution $\q_h$ is \textit{projected} from elements $\mathcal{C}_{11}$ and $\mathcal{C}_{7}$ 
into the cell averages $\bar \u$ of their virtual children at time $t^n_0 + \Delta t_2$. 
Now, reconstruction for cells $\mathcal{C}_{18}-\mathcal{C}_{20}$ can be performed 
exactly as for the uniform Cartesian case and subsequently the local space--time predictor can 
be carried out with the new initial data $\w_h$. For the time update in the next time interval $[t_0^n+\Delta t_2;t_0^n+2\Delta t_2]$ the predictor solution in elements $\mathcal{C}_{11}$ and 
$\mathcal{C}_{7}$ is still valid, hence the numerical flux can again be directly computed between 
cells as described before. This procedure is carried out until all cells of level $\ell=2$ reach the 
time $t^n_0 + \mathfrak{r} \Delta t_2$. Then, their future time is $t^n_0 + (\mathfrak{r}+1) \Delta t_2$, 
hence the update criterion \eqref{eqn.update.criterion} is no longer fulfilled for level $\ell=2$ and 
the next coarser level $\ell=1$ can be updated. Cell $\mathcal{C}_{10}$ is only surrounded by 
active neighbors on the same level, hence the standard finite volume scheme on uniform Cartesian 
mesh can be used. For cells $\mathcal{C}_{9}$ and $\mathcal{C}_{8}$ the previously described procedure 
of flux computation between fine and coarse cell applies. Cell $\mathcal{C}_{11}$ now illustrates 
the last special case, namely the computation of the numerical flux between the coarse cell $\mathcal{C}_{11}$ and fine cell $\mathcal{C}_{18}$. Actually, the solution is particularly simple. 
Due to the requirement of \textit{conservation}, the coarse cell $\mathcal{C}_{11}$ must \textit{not}  compute any new flux at the interface with element $\mathcal{C}_{18}$, since all the necessary fluxes 
have already been computed before on the finer level. It is just sufficient for cell $\mathcal{C}_{11}$ 
to sum up the fluxes stored in the memory variable of its virtual child $\mathcal{C}_{23}$. In this way, the sum of time integrals of the 
fluxes on the left border of element $\mathcal{C}_{18}$ is equal to the time integral of the fluxes 
on the right border of $\mathcal{C}_{11}$ in the time interval $[t^n_0;t^n_0+\Delta t_1]$, which makes 
the method conservative. In the most general 3D case, the numerical flux computed at an interface  containing elements of level $\ell$ and $\ell+1$ reads 
\begin{equation}
  \mathbf{f}_{i+\halb,jk} = \frac{1}{\Delta t_\ell} \frac{1}{\Delta y_\ell} \frac{1}{\Delta z_\ell}  
  \sum \limits_{ii=1}^{\mathfrak{r}} \sum \limits_{jj=1}^{\mathfrak{r}} \sum \limits_{kk=1}^{\mathfrak{r}} 
  \, \, 
  \int \limits_{\mathcal{T}_{ii}}
  \int \limits_{\mathcal{Y}_{jj}}  
  \int \limits_{\mathcal{Z}_{kk}}  
  \mathbf{\tilde f}(\q_h^-,\q_h^+) \, dz \, dy \, dt, 
\label{eqn.memory.sum} 
\end{equation} 
with the integration intervals above defined as 
\vspace{-5mm} 
\begin{eqnarray} 
\mathcal{T}_{ii}&=&[t^n_\ell    + (ii-1) \Delta t_{\ell+1};t^n_\ell    + ii \Delta t_{\ell+1}], \nonumber \\  
\mathcal{Y}_{jj}&=&[y_{j-\halb} + (jj-1) \Delta y_{\ell+1};y_{j-\halb} + jj \Delta y_{\ell+1}], \nonumber \\ 
\mathcal{Z}_{kk}&=&[z_{k-\halb} + (kk-1) \Delta z_{\ell+1};z_{j-\halb} + kk \Delta z_{\ell+1}]. 
\end{eqnarray} 
The flux on the left hand side of \eqref{eqn.memory.sum} corresponds to the final averaged flux for 
the coarse grid cell, while the integrals on the right hand side of \eqref{eqn.memory.sum} are the
integrated fluxes for each time step of each fine grid cell adjacent to the coarse cell. 
The practical implementation of Eqn. \eqref{eqn.memory.sum} is conveniently achieved by the sum over
the memory variables, as described before. 
 
This completes the description of all cases to be treated by the AMR algorithm
concerning flux computation across elements on different refinement levels.  

\subsubsection{Projection}
\label{sec:projection}

Projection is the typical AMR operation, sometimes called "coarse-to-fine prolongation", by which 
values an active mother assigns values to the virtual children ($\sigma=1$) at intermediate times 
via standard $L_2$ projection. For this purpose, the space--time polynomials $\q_h$ can be 
conveniently evaluated at any time. This operation is needed for performing the reconstruction on 
the finer grid level at intermediate times. The projection operator for a cell $\mathcal{C}_m$ on
level $\ell$ is simply given by evaluating the space--time polynomial $\q_h$ of its \textit{mother}  
at any given time $t^n_\ell$ as follows: 
\begin{equation}
  \bar \u_m(t^n_\ell) = \frac{1}{\Delta x_\ell} \frac{1}{\Delta y_\ell}  \frac{1}{\Delta z_\ell}  \int \limits_{\mathcal{C}_m} \q_h(\mathbf{x},t^n_\ell) d \mathbf{x}.  
\end{equation} 

\subsubsection{Averaging}
\label{sec:averaging}
Averaging is another typical AMR operation by which a virtual mother cell ($\sigma=-1$) obtains 
its cell average by averaging recursively over the cell averages of all its children and their
children at higher refinement levels. 
Let us denote the set of children of a cell $\mathcal{C}_m$ by $\mathcal{B}_m$,
then the averaging operator is given by 
\begin{equation}
  \bar \u_m = \frac{1}{\mathfrak{r}^d} \sum \limits_{\mathcal{C}_k \in \mathcal{B}_m} \bar \u_k. 
\label{eqn.average} 
\end{equation}  

\subsection{Overall efficiency and MPI parallelization of higher order AMR}
\label{sec:AMR-parallelization}

In the following we study quantitatively the overhead introduced by the high order
one-step ADER-WENO AMR method proposed in this article. For this purpose, we report 
the CPU times needed for the simulation of the two-dimensional explosion problem 
discussed in more detail in section \ref{sec:num-tests} for uniform and AMR grids
for second to fourth order ADER-WENO schemes. The detailed CPU time results are 
summarized for all cases in Table \ref{tab.efficiency} and are normalized with 
respect to the standard second order scheme on uniform mesh. The data refer to the 
average CPU time needed for the update of one single real element, which has been 
computed by dividing the total wallclock time needed for the simulation by the number 
of time updates of the active elements ($\sigma=0$) contained in the domain. Hence, the results 
reported in Table  \ref{tab.efficiency} include the entire overhead necessary for the 
update, averaging and projection of the virtual ghost cells needed in the AMR approach. 
In the table, we also report separately the total overhead introduced by the AMR 
approach in percent for convenience. The CPU times have been obtained on one single 
core of an Intel i7-2600 CPU with 3.4 GHz clock speed and 12 GB RAM.

\begin{figure}
\begin{center}
\includegraphics[angle=0,width=0.45\textwidth]{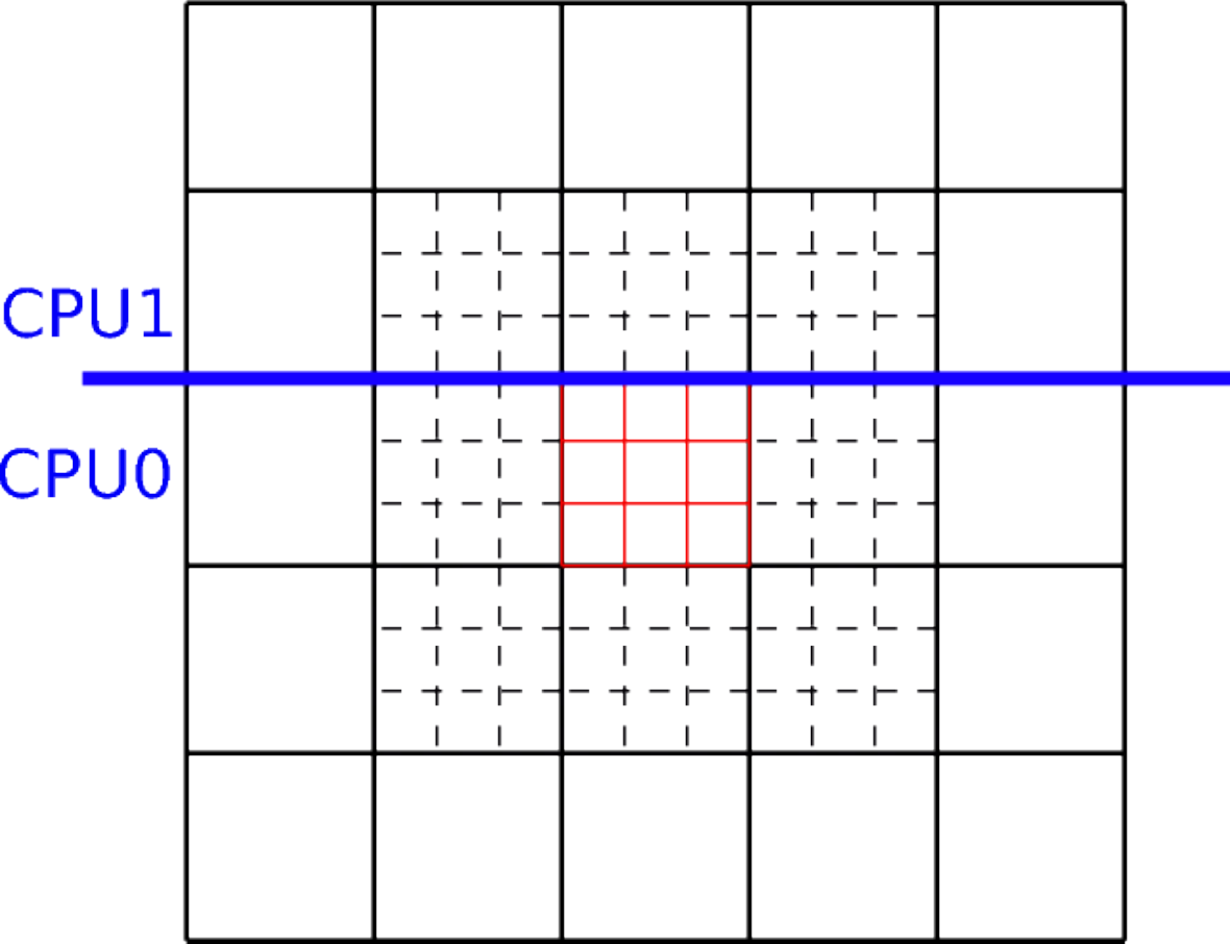}
\hspace{1cm}
\includegraphics[angle=0,width=0.45\textwidth]{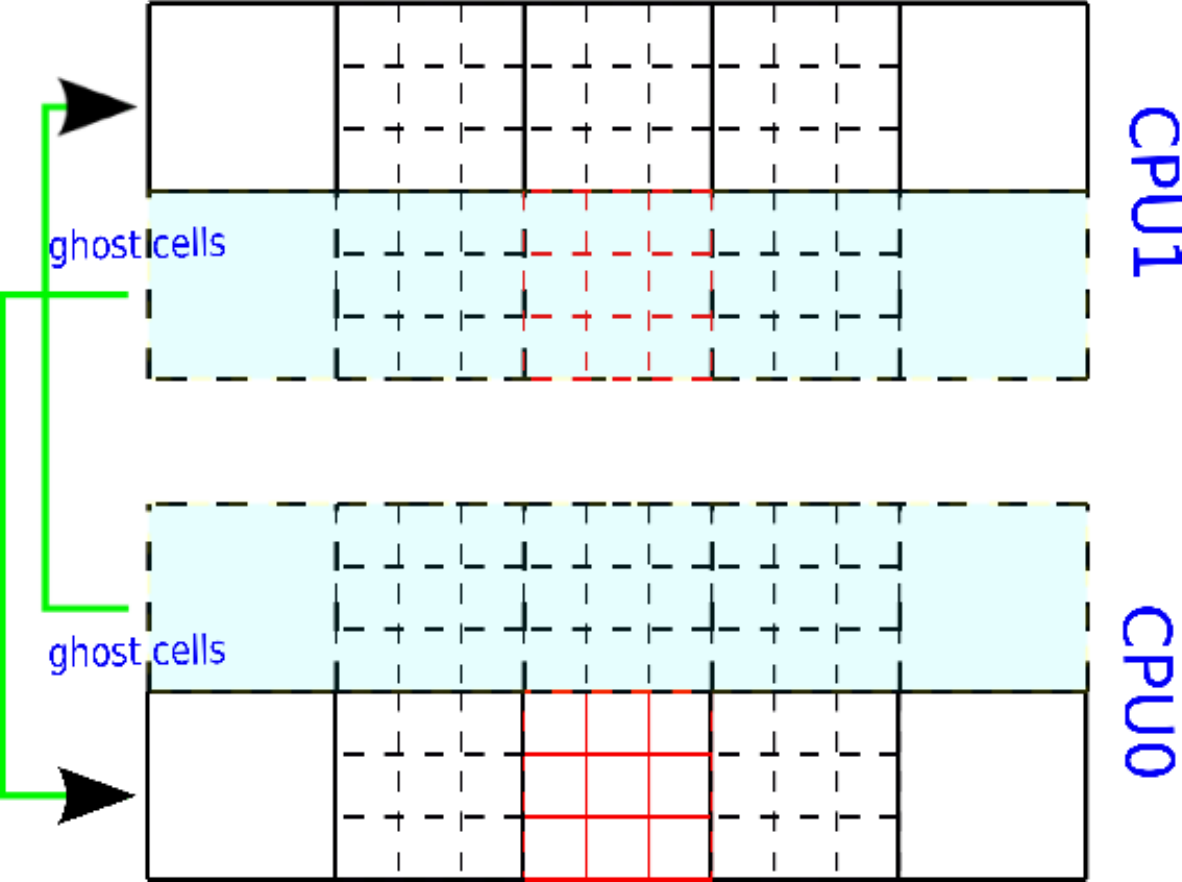}
\caption{Cell refinement at the MPI-border between two processors. 
Left panel: processors CPU0 and CPU1 must exchange information about which cells are refined (solid red line) or virtually 
refined (dashed black line).
Right panel: Each processor has a MPI--ghost zone of cells that are a copy 
of the true cells managed by the adjacent processor.}
\label{fig_refine-at-mpi-border}
\end{center}
\end{figure}

\begin{table}[!b]   
\caption{Assessment of the overall efficiency of high order one--step ADER-WENO 
schemes on space-time adaptive AMR grids ($\mathfrak{r}=4, \ell=2$). Normalized 
average CPU time per real element update with respect to the second order scheme 
on uniform grid.}
\begin{center} 
\renewcommand{\arraystretch}{1.0}
\begin{tabular}{cccc} 
\hline
 Scheme order & Uniform grid & AMR grid & Total AMR overhead \\ 
\hline 
  $\mathcal{O}2$ & 1.00  & 1.15  & 15 \% \\ 
  $\mathcal{O}3$ & 3.18  & 3.82  & 20 \% \\ 
  $\mathcal{O}4$ & 8.64  & 10.82 & 25 \% \\ 
\hline 
\end{tabular} 
\end{center}
\label{tab.efficiency}
\end{table} 
We furthermore have parallelized the three dimensional ADER-WENO code through the standard Message Passing Interface (MPI). Any AMR implementation poses additional challenging problems 
to the parallelization task, which become manifest when refinement of a cell ${\mathcal C}_m$  occurs at the MPI border between two processors (see left panel of  Fig.~\ref{fig_refine-at-mpi-border}). In this case, in fact, proper communication among the processors must be established in order to spread the knowledge about which cells must be 
either virtually refined or activated.  
For this purpose, each processor stores in its memory also MPI ghost--cells that are a copy of  the true cells, managed by the adjacent processor. In the practical implementation, we have 
found convenient to assign a negative integer number to each cell in the MPI ghost--zone,  
thus making the distinction with respect to real cells very transparent. When a cell in 
the domain of the processor CPU0, at the border with the domain of the processor CPU1, 
is refined (see right panel of Fig.~\ref{fig_refine-at-mpi-border}), CPU0 informs CPU1 that 
(i) a number of real cells belonging to CPU1 must be virtually refined, and (ii) that one 
cell (the one at the border) must be virtually refined in the MPI--ghost zone of CPU1. 
This information is used by CPU1 to (virtually) refine its corresponding cells in its 
MPI--ghost zone. Before doing that, CPU1 must also check whether such cells have already  received an instruction of virtual refinement internal to CPU1. The link between the true 
cells of CPU0 and those belonging to the MPI--ghost zone of CPU1 is obtained via so--called  \textit{exchange lists}. The exchange lists are used for the MPI communication that is 
necessary during the adaptive mesh refinement procedure, as well as to exchange the 
information about the cell averages $\bar \u$ and the space--time polynomials $\q_h$ between 
the processors. We stress that our present MPI-AMR implementation does not yet provide 
dynamic load-balancing among  processors, which is a rather complex topic that will be  considered in the future. 

%

\section{Numerical Tests}
\label{sec:num-tests}

In all the following numerical test problems we have used the density as indicator function in \eqref{eqn.indicator}, hence $\phi(\u)=\rho$. 

\subsection{Euler equations}

The first session of tests considers a sequence of applications for which the classical Euler equations 
of compressible gas dynamics are solved. In three space dimensions the vectors of the conserved variables 
$\bf{u}$ and of the fluxes $\bf{f}$, $\bf{g}$ and $\bf{h}$ are given respectively by
\begin{equation}
{\bf{u}}=\left(\begin{array}{c}
\rho \\ \rho v_x \\ \rho v_y \\ \rho v_z \\ E 
\end{array}\right) \!\!,  \, 
{\bf f}=
\left(\begin{array}{c}
\rho v_x \\
\rho v_x^2 + p \\
\rho v_xv_y \\
\rho v_xv_z \\
v_x(E+p)
\end{array}\right)\!\!,  \, 
{\bf g}=
\left(\begin{array}{c}
\rho v_y \\
\rho v_xv_y \\
\rho v_y^2 + p \\
\rho v_yv_z \\
v_y(E+p)
\end{array}\right)\!\!,  \, 
{\bf h}=
\left(\begin{array}{c}
\rho v_z \\
\rho v_xv_z \\
\rho v_yv_z \\
\rho v_z^2 + p \\
v_z(E+p)
\end{array}\right),
\label{eq:Euler-system}
\end{equation}
where $v_x$, $v_y$ and $v_z$ are the velocity components, $p$ is the pressure, $\rho$ is the mass density,
$E=p/(\gamma-1)+\rho (v_x^2+v_y^2+v_z^2)/2$ is the total energy density, while $\gamma$ is the adiabatic index. 

\paragraph*{2D isentropic vortex.}
%
The first test considered is a two-dimensional convected isentropic vortex, see e.g. \cite{HuShuTri}.  
The computational domain is $\Omega=[0;10]\times[0;10]$ and the initial conditions are given by a 
perturbation added to a uniform mean flow   
\begin{equation}
\left( \rho,v_x,v_y,v_z,p \right) =(1+\delta\rho, 1+\delta v_x, 1+\delta v_y, 0, 1+\delta p)\,,
\end{equation} 
with 
\begin{equation}
\left(\begin{array}{c}
\delta \rho \\ \delta v_x \\ \delta v_y \\ \delta p 
\end{array}\right)
=
\left(\begin{array}{c}
(1+\delta T)^{1/(\gamma-1)}-1 \\
-(y-5)\epsilon/2\pi \exp{[0.5(1-r^2)]} \\
\phantom{-}(x-5)\epsilon/2\pi \exp{[0.5(1-r^2)]} \\
(1+\delta T)^{\gamma/(\gamma-1)}-1
\end{array}\right).~~~
\label{eq:pert}
\end{equation}
The perturbation $\delta T$ in the temperature is 
\begin{equation}
\delta T=-\frac{\epsilon^2(\gamma-1)}{8\gamma\pi^2}~\exp{(1-r^2)}\,,
\end{equation}
with $r^2=(x-5)^2+(y-5)^2$, vortex strength $\epsilon=5$ and adiabatic index $\gamma=1.4$. The refinement factor adopted is $\mathfrak{r}=3$. 
In Table~\ref{tab.conv2} we have reported the results of the convergence tests, where we have used the third and fourth order version of the method. 
The convergence rates have been computed with respect to an initially uniform mesh, as proposed by Berger and Oliger in \cite{Berger-Oliger1984}. 

\begin{table}[!t]   
\caption{Numerical convergence results for the isentropic vortex test using the third and fourth order version of the one--step ADER-WENO finite volume scheme
presented in this article. The error norms refer to the variable $\rho$ (density) at the final time $t_f=10$. The asterisk $^\ast$ refers to a uniform grid. 
}
\begin{center} 
\renewcommand{\arraystretch}{1.0}
\begin{tabular}{lccllc} 
\hline
  \hline
  $\ell_{\rm max}=1$ &   & & & &   \\ 
  \hline
  $N_G\times N_G$  & $\epsilon_{L_2}$ & $\mathcal{O}(L_2)$ &  $N_G\times N_G$ & $\epsilon_{L_2}$ & $\mathcal{O}(L_2)$  \\ 
\hline
                     & & {$\mathcal{O}3$} & & & {$\mathcal{O}4$}  \\
\hline
  12$\times$12$^\ast$   & 5.0130E-01 &      & 10$\times$10$^\ast$  & 5.1496E-01 &       \\ 
  24$\times$24          & 4.9128E-02 & 3.35 & 15$\times$15         & 8.9093E-02 & 4.33   \\ 
  36$\times$36          & 1.6922E-02 & 3.08 & 21$\times$21         & 2.7906E-02 & 3.93  \\ 
  48$\times$48          & 7.5867E-03 & 3.02 & 28$\times$28         & 8.3878E-03 & 4.00   \\  
  72$\times$72          & 2.7106E-03 & 2.91 & 42$\times$42         & 1.5780E-03 & 4.03   \\
  108$\times$108        & 1.0579E-04 & 2.80 & 63$\times$63         & 4.0931E-04 & 3.88   \\
\hline
  \hline
  $\ell_{\rm max}=2$ &   & & & &   \\ 
  \hline
  $N_G\times N_G$  & $\epsilon_{L_2}$ & $\mathcal{O}(L_2)$ & $N_G\times N_G$ & $\epsilon_{L_2}$ & $\mathcal{O}(L_2)$  \\ 
\hline
                     & & {$\mathcal{O}3$} & & & {$\mathcal{O}4$}  \\
\hline                     
 12$\times$12$^\ast$   & 5.0131E-01 &      & 10$\times$10$^\ast$  & 5.1496E-01 &       \\ 
 24$\times$24          & 1.5223E-02 & 5.04 & 15$\times$15         & 3.2990E-02 & 6.78   \\ 
 36$\times$36          & 5.6974E-03 & 4.08 & 21$\times$21         & 1.2157E-02 & 5.05   \\ 
 48$\times$48          & 2.3935E-03 & 3.86 & 28$\times$28         & 4.5922E-03 & 4.58   \\  
 72$\times$72          & 7.5147E-04 & 3.63 & 42$\times$42         & 1.0334E-03 & 4.33   \\
 108$\times$108        & 5.4038E-04 & 3.11 & 63$\times$63         & 2.4593E-04 & 4.15   \\
\hline 
\end{tabular} 
\end{center}
\label{tab.conv2}
\end{table} 
\begin{figure}
\begin{center}
\includegraphics[angle=0,width=6.0cm]{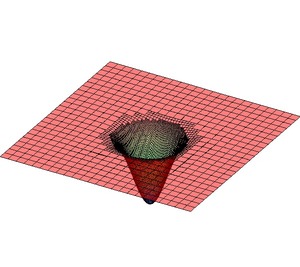}
\includegraphics[angle=0,width=6.0cm]{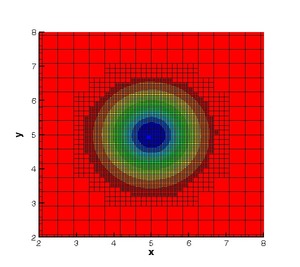}
\caption{  
{Isentropic-Vortex test at the final time $t=10$}. Left panel: Contour plot of the mass density. 
Right panel: zoom into the AMR grid. Two levels of refinement have been adopted ($\ell_{\max}=2$). 
}
\label{fig_shu_vortex}
\end{center}
\end{figure}
%

\begin{figure}
\begin{center}
\begin{tabular}{lr}
\includegraphics[angle=0,width=0.4\textwidth]{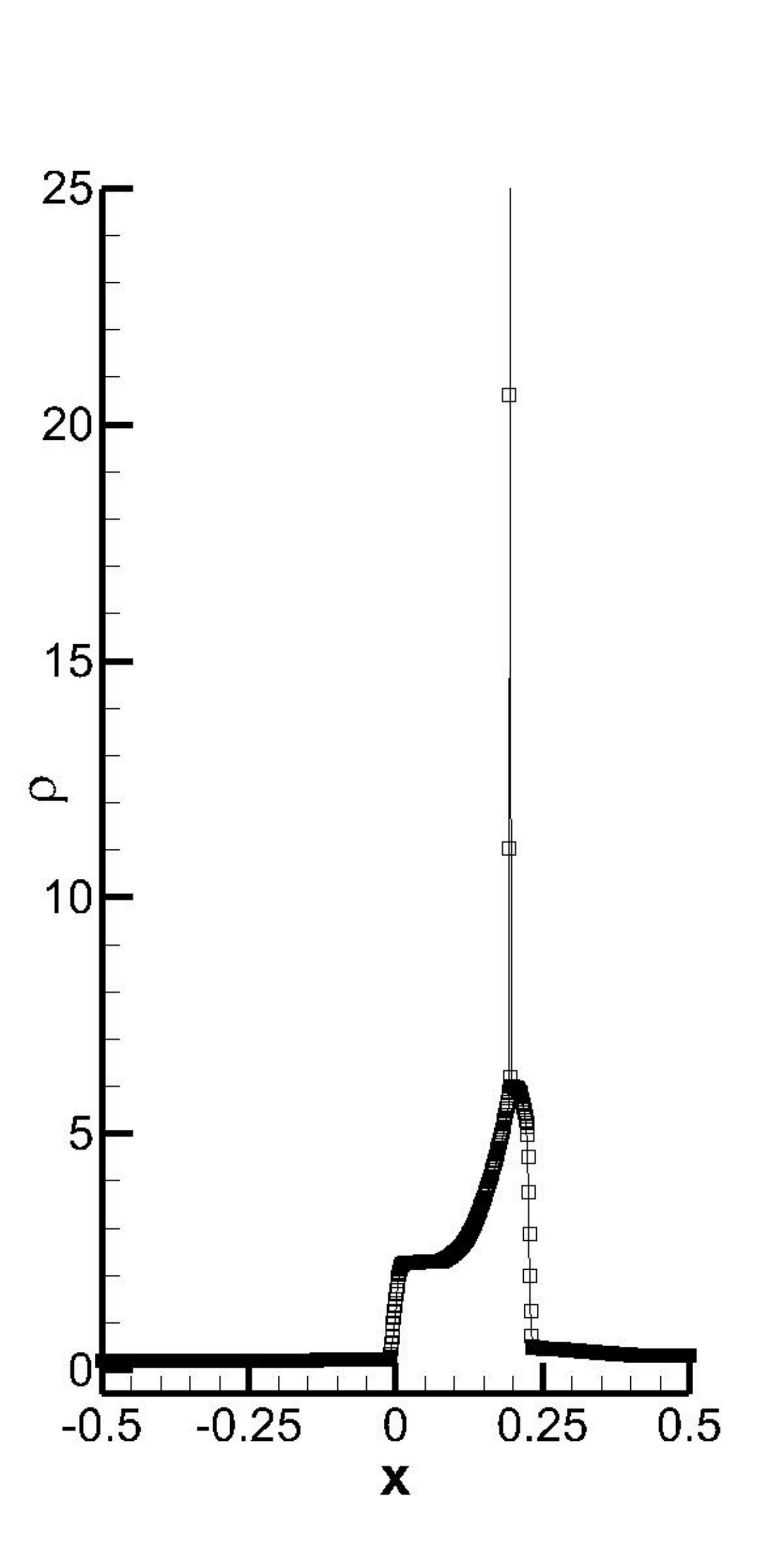}
\includegraphics[angle=0,width=0.4\textwidth]{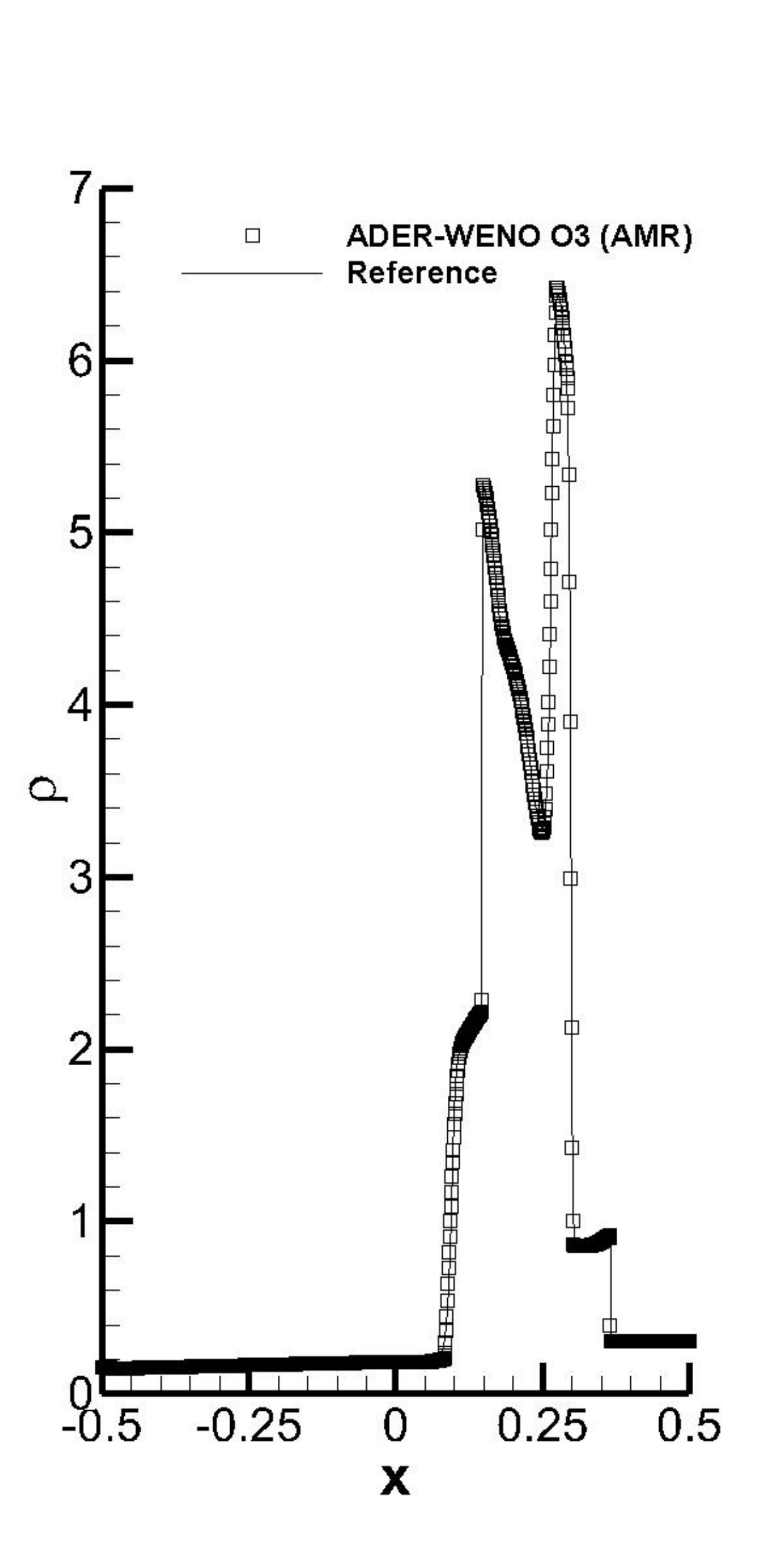}
\end{tabular}
\caption{
{Interacting Blast Waves Test}. Profile of the mass density at time $t=0.028$ (left panel)
and at the final time $t=0.038$ (right panel). 
Two levels of refinement have
been adopted from an initial grid $300\times10$. 
}
\label{fig_blast_wave}
\end{center}
\end{figure}

\paragraph*{Interacting blast waves.}
%
This test, originally proposed by \cite{woodwardcol84}, has by now become a classical problem 
in computational fluid dynamics
and consists in the interaction of blast waves with
initial conditions given by
\begin{equation}
\label{blast-wave}
(\rho,v_x,p)= \left\{
\begin{array}{llll}
(1.0,0.0,10^3) &   {\rm if} & -0.5 < x < -0.4 \,, \\
 
(1.0,0.0,10^{-2}) &   {\rm if} & -0.4 < x < 0.4 \,, \\
 
(1.0,0.0,10^2) &   {\rm if} & \phantom{-} 0.4 < x < 0.5 \,. 
\end{array} \right.
\end{equation}
Although one dimensional, we have evolved this problem in two spatial dimensions over the domain $[-0.5,0.5]\times[-0.5,0.5]$,
using reflecting boundary conditions in $x$ direction and periodic boundary conditions along the $y$ direction. The adiabatic 
index has been chosen as $\gamma=1.4$. 
Fig.~\ref{fig_blast_wave} shows the results of our test, where we have used two levels of refinement from an original uniform grid 
with $300\times 10$ cells, and adopting a third order ADER-WENO scheme. 
The left panel shows the solution at the time when the two waves hit each other from opposite directions, producing a very strong 
density peak. The right panel, on the other hand, shows the solution at the final time after the waves have crossed each other. 
A reference solution is also reported, obtained with a traditional finite difference TVD method
using $3600$ grid-points.
%
%
\paragraph*{Explosion problems in two and three space dimensions.}
%

\begin{figure}
\begin{center}
\begin{tabular}{lr}
\includegraphics[angle=0,width=0.45\textwidth]{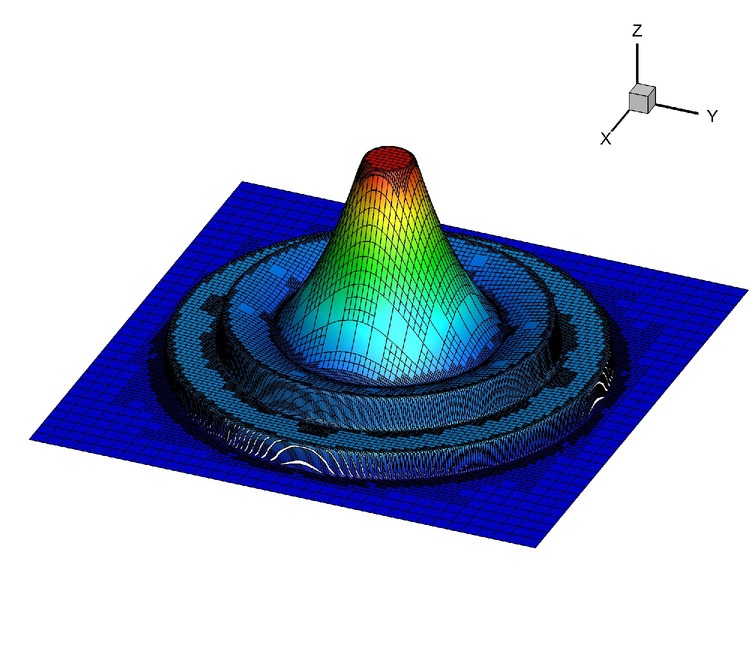} 
\includegraphics[angle=0,width=0.45\textwidth]{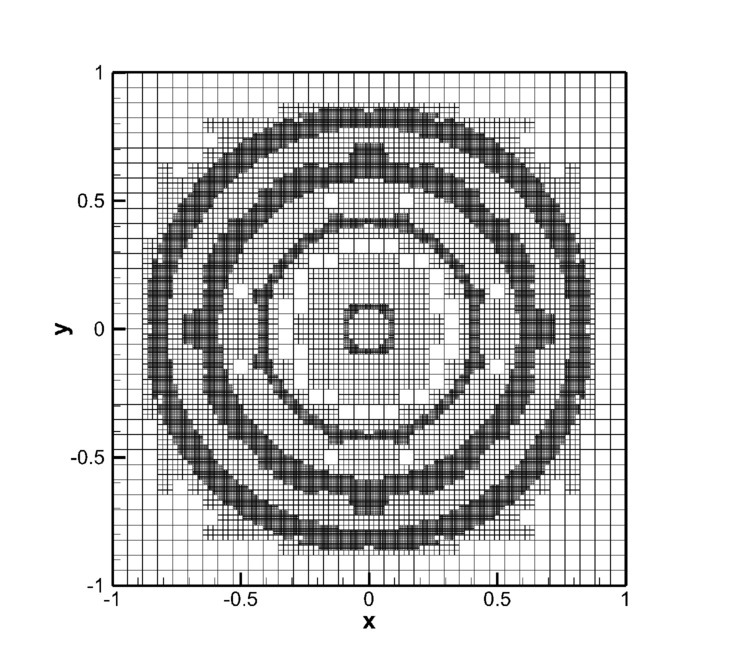}
\end{tabular}
\caption{{Explosion test in two space dimensions.} Density distribution at time $t=0.25$ obtained with a fourth order ADER-WENO scheme (left) and corresponding 
final AMR grid configuration (right). }
\label{fig_explosion2D}
\end{center}
\end{figure}

\begin{figure}
\begin{center}
\begin{tabular}{lr}
\includegraphics[angle=0,width=0.45\textwidth]{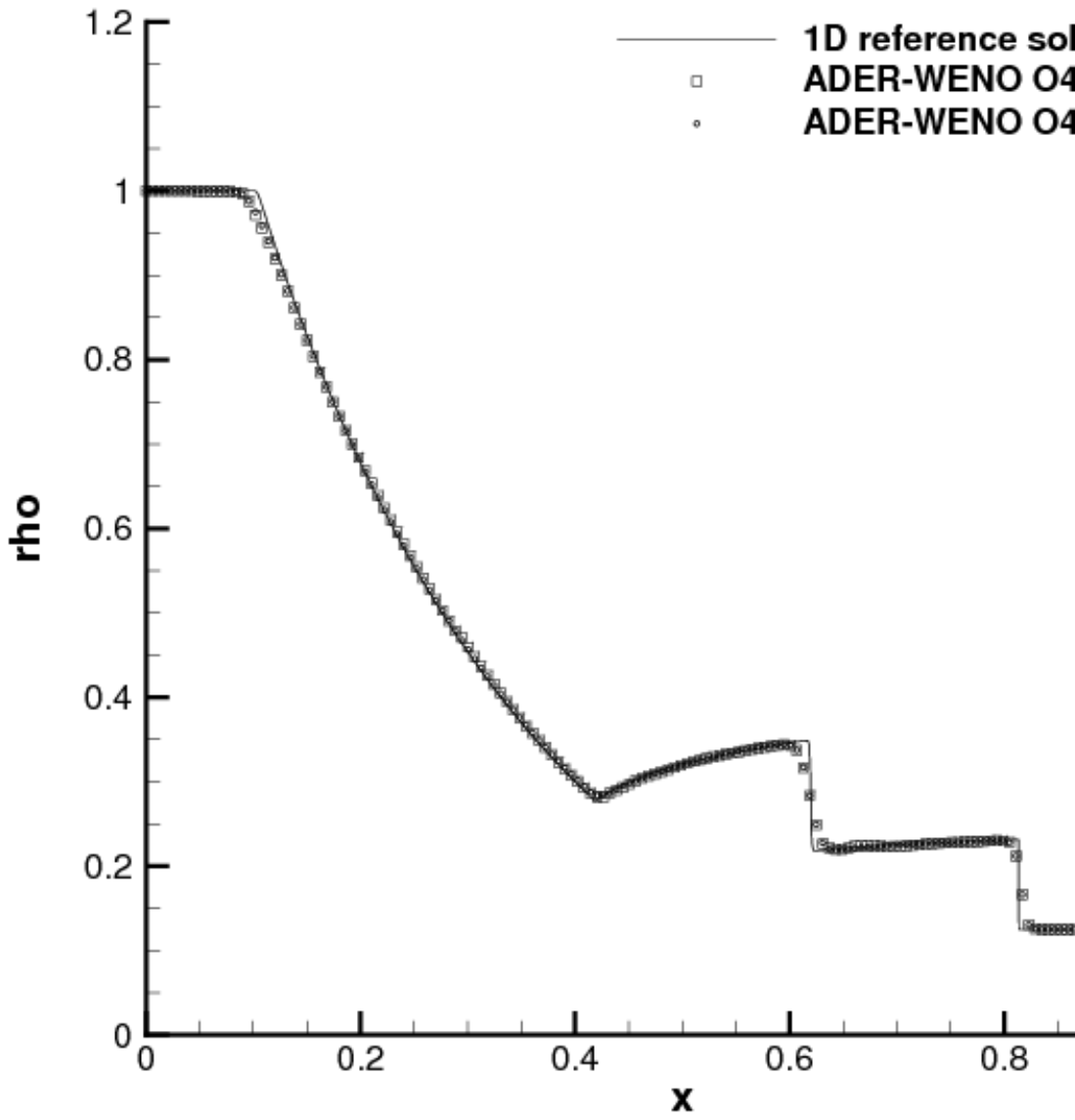} & 
\includegraphics[angle=0,width=0.45\textwidth]{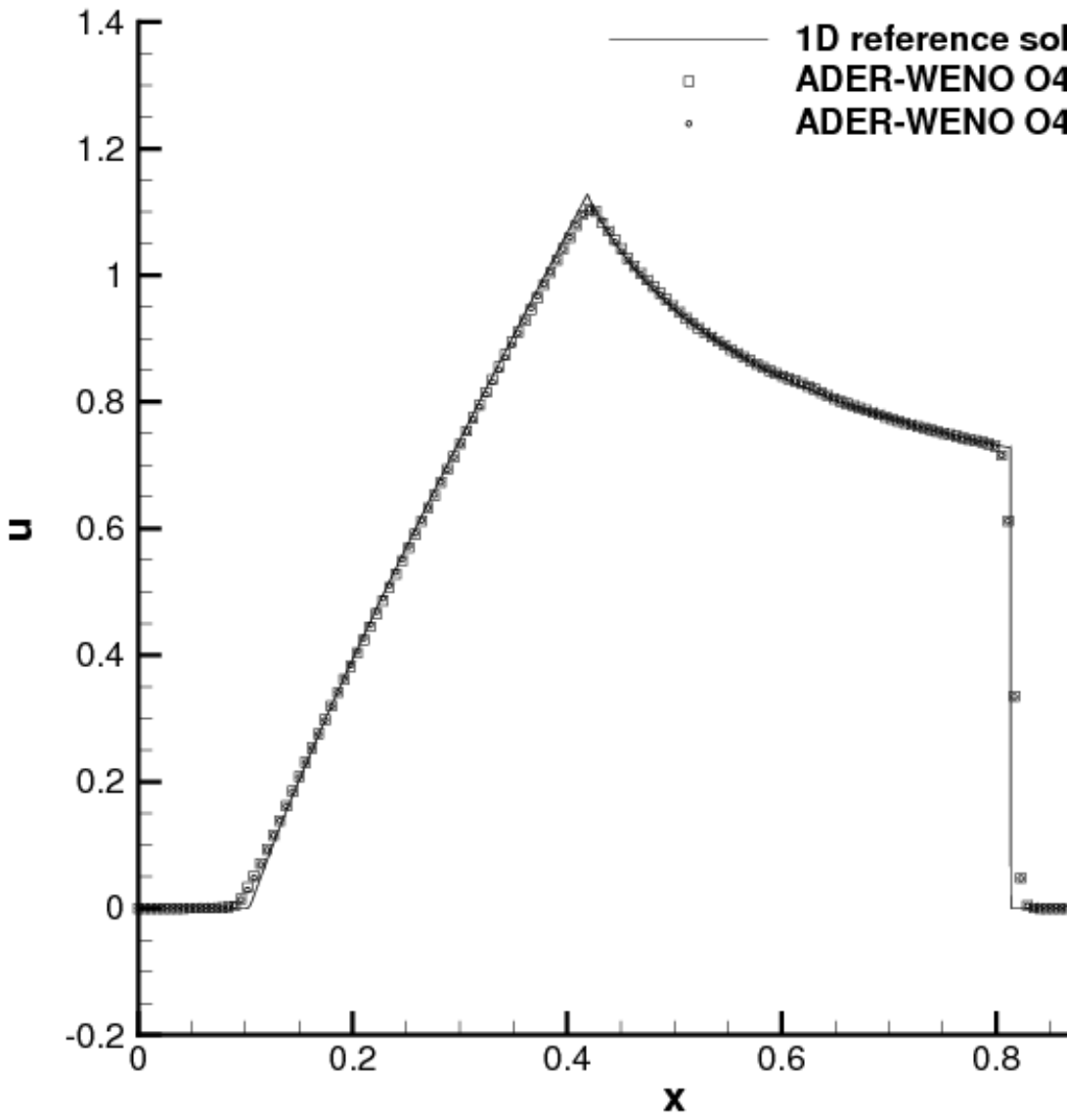}   
\end{tabular} 
\caption{{Explosion test in two space dimensions.} One-dimensional cut along the positive $x$-axis through the fourth order ADER-WENO solution obtained on the AMR 
grid for density (left) and velocity (right). The solution computed on a uniform fine mesh corresponding to the finest AMR grid level is also shown. }
\label{fig_explosion2D.cut}
\end{center}
\end{figure}

In this test, proposed in \cite{titarevtoro,toro-book}, we solve the Euler equations on the computational domain $\Omega=[-1; 1]^d$, where $d$ 
denotes the number of space dimensions. The initial flow variables take constant values for $r\leq R$ and for $r\geq R$, separated by a cylindrical 
or spherical discontinuity, respectively. Therefore the initial condition is given by 
\begin{equation}
  \mathbf{u}(\mathbf x,0) = \left\{ \begin{array}{ccc} \mathbf{u}_i & \textnormal{ if } & r \leq R, \\ 
                                                       \mathbf{u}_o & \textnormal{ if } & r > R.
  \end{array} \right. 
\end{equation} 
Here, $R=0.4$ denotes the radius of the initial discontinuity, $\mathbf{x}$ is the vector of spatial coordinates with the radial coordinate $r = \sqrt{\mathbf{x}^2}$.  
$\mathbf{u}_i$ and $\mathbf{u}_o$ are the inner and outer states, respectively, listed in detail in Table \ref{tab.3d.explos.ic}. The 
adiabatic index of the ideal-gas equation of state has been set to $\gamma=1.4$. Due to the symmetry of the problem, which is cylindrical in 
the two-dimensional case and spherical in the three-dimensional case, the solution can be compared with an equivalent one dimensional problem 
in radial direction $r$, see \cite{toro-book}: 
\begin{table}[!b]
\caption{Inner and outer initial states for the multidimensional explosion test problems. 
  The last column reports the final simulation time $t_e$.} 
\vspace{0.5cm}
\renewcommand{\arraystretch}{1.0}
\begin{center}
\begin{tabular}{ccccccc}
\hline
 Case  & $\rho$ & $p$ & $v_x$  & $v_y$ & $v_z$ & $t_e$ \\
\hline
Inner  & 1.0    & 1.0  & 0.0 & 0.0 & 0.0 &  0.25 \\
Outer  & 0.125  & 0.1  & 0.0 & 0.0 & 0.0 &       \\
\hline
\end{tabular}
\end{center}
\label{tab.3d.explos.ic}
\end{table}
\begin{equation}
\frac{\partial}{\partial t} \left( \begin{array}{c} \rho \\ \rho u \\ E \end{array} \right) +  
\frac{\partial}{\partial r} \left( \begin{array}{c} \rho u \\ \rho u^2 + p \\ u (E + p) \end{array} \right) = 
-\frac{d-1}{r}\left(\begin{array}{c} 
\rho u \\ \rho u^2 \\ u (E+p)  
\end{array}\right),    
\label{eq.1D}
\end{equation}
where $u$ is the radial velocity component. The 1D reference solution has been computed by a classical second order TVD finite volume scheme on a very fine 
mesh composed of 10000 grid zones and using the Osher-type flux proposed in \cite{OsherUniversal}. The two-dimensional AMR simulations have been carried out 
with a fourth order ADER-WENO scheme on a level zero grid with $34 \times 34$ control volumes, using the Osher flux \eqref{eqn.osher} and $\mathfrak{r}=3$ and $\ell_{\max}=2$. This leads to 
an equivalent resolution on a uniform fine grid of $306 \times 306 = 93636$ points. Fig.~\ref{fig_explosion2D} shows a 3D plot of the density distribution 
obtained for the cylindrical explosion case, as well as the AMR grid configuration at the final time $t=0.25$.   
Fig.~\ref{fig_explosion2D.cut} shows the results obtained on a one-dimensional cut along the $x$-axis, together with the 1D reference solution according to 
\eqref{eq.1D}. The cut is performed on equidistant points, evaluating locally the reconstruction polynomials $\mathbf{w}_h$. For comparison, we also show the results 
obtained on the uniform fine grid of $306 \times 306$ grid points that corresponds to the finest AMR grid level. Both simulations agree very well with the 1D 
reference solution. Furthermore one can note only very little differences in the numerical results obtained with AMR and without AMR, i.e. on the uniform 
grid. However, the simulation on the AMR grid took only 4 minutes on 4 cores of an Intel i7-2600 CPU with 3.4 GHz clock speed and 12 GB of RAM, while the fine 
uniform grid computation needed 14 minutes on 4 cores on the same machine, hence it took \textit{3.5 times longer}. This clearly confirms that even though the use of AMR adds a 
certain overhead of about 25 $\%$ to the fourth order finite volume scheme, according to Table \ref{tab.efficiency}, the use of space--time adaptive meshes 
can significantly speed up multidimensional computations also for higher order schemes. The total number of AMR grid cells present at the final time was 
28036, compared to the 93636 cells of the uniform grid.  


\begin{figure}
\begin{center}
\includegraphics[angle=0,width=0.75\textwidth]{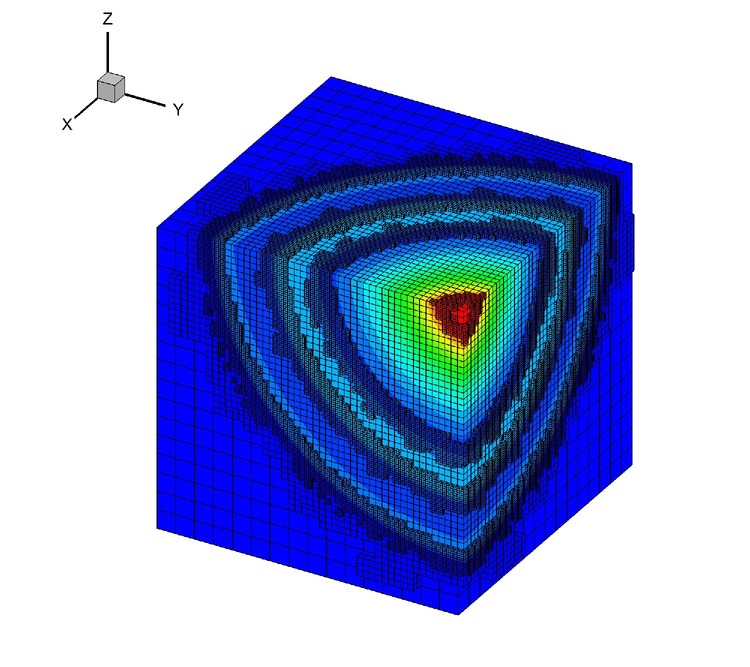}
\caption{{Explosion test in three space dimensions.} AMR grid structure at time $t=0.25$ and density contour colors. }
\label{fig_explosion3D.grid}
\end{center}
\end{figure}

\begin{figure}
\begin{center}
\begin{tabular}{lr} 
\includegraphics[angle=0,width=0.45\textwidth]{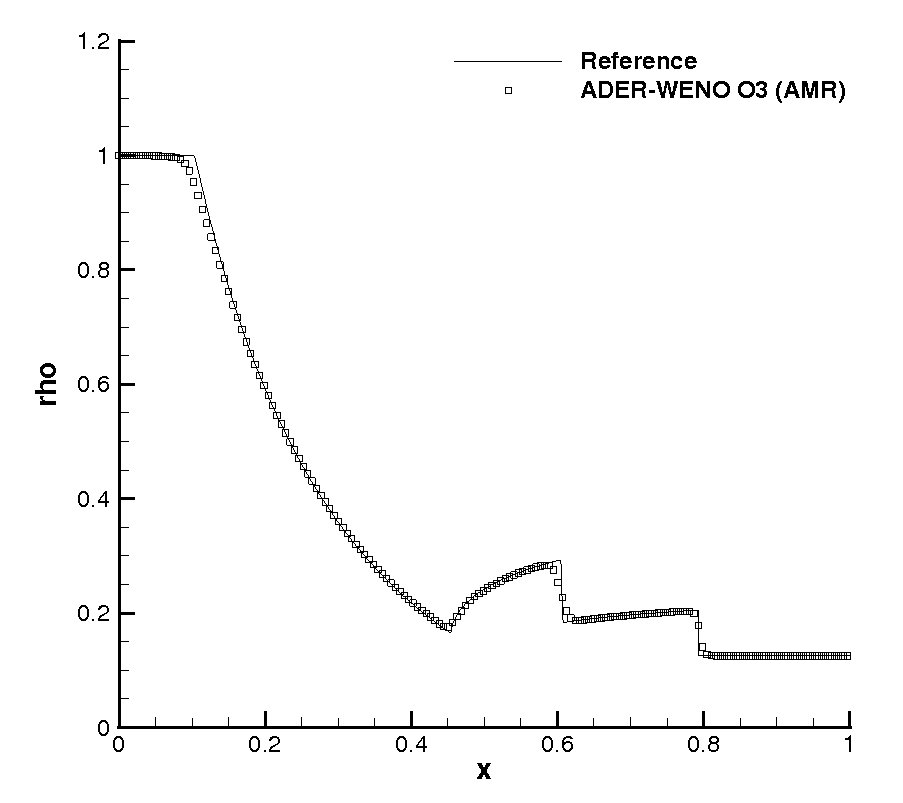}      & 
\includegraphics[angle=0,width=0.45\textwidth]{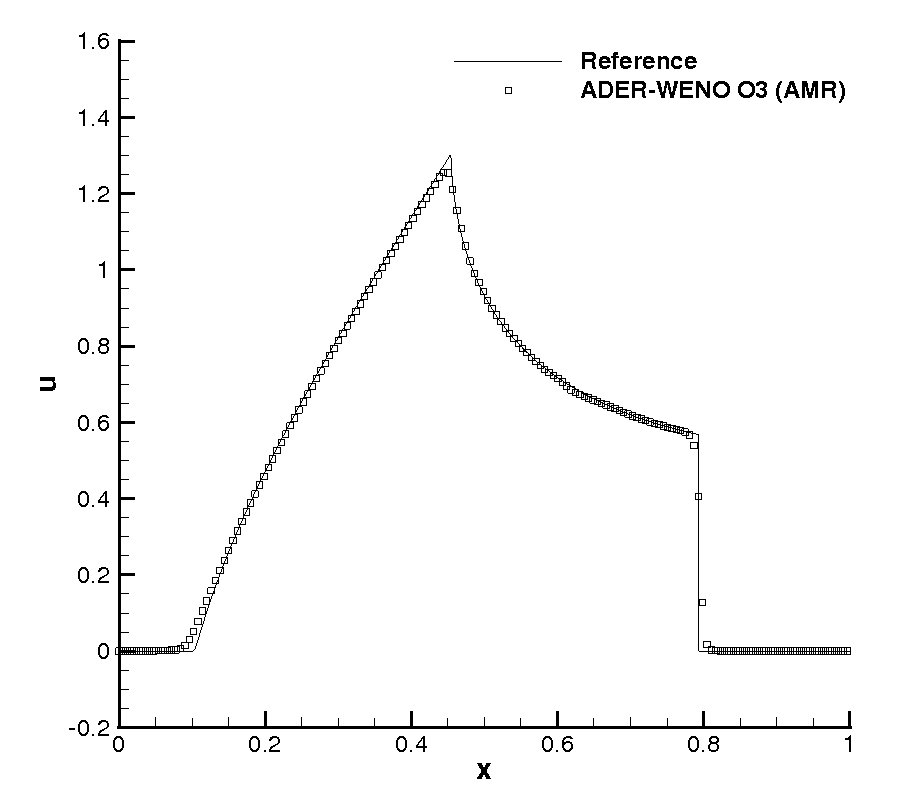}        \\ 
\includegraphics[angle=0,width=0.45\textwidth]{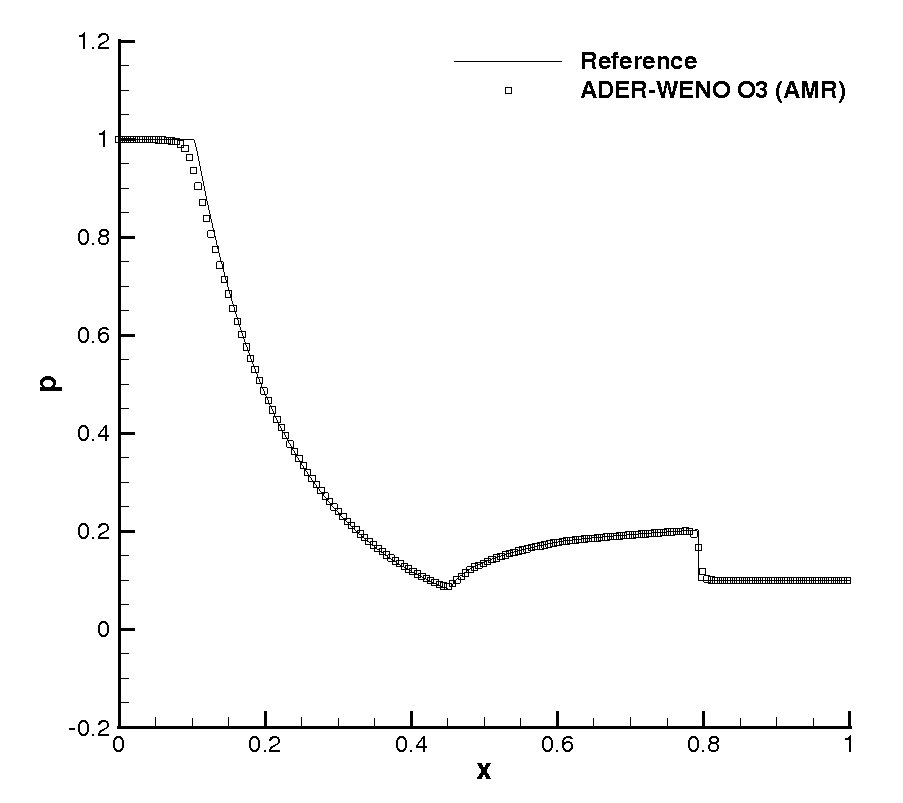}        & 
\includegraphics[angle=0,width=0.45\textwidth]{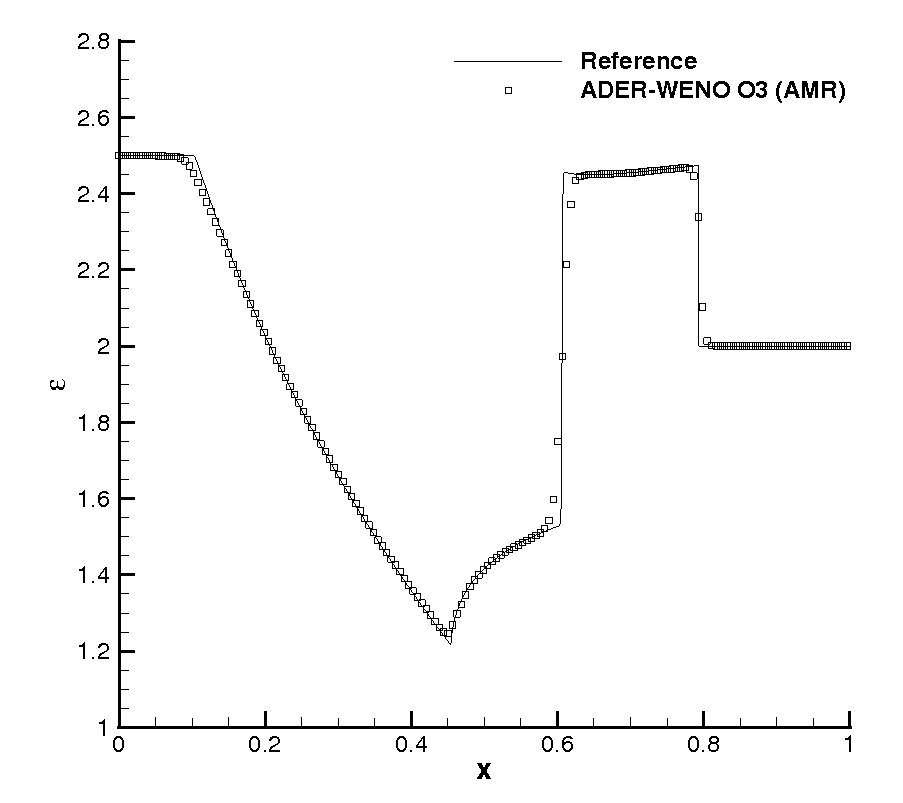}     
\end{tabular}  
\caption{{Explosion test in three space dimensions.} From top left to bottom right: comparison of the 1D reference solution with the numerical solution obtained 
with a third order ADER-WENO scheme on AMR mesh at time $t=0.25$. One-dimensional cuts along the $x$-axis are shown for density, velocity, pressure and internal energy.
}
\label{fig_explosion3D}
\end{center}
\end{figure}

Similarly, Fig.~\ref{fig_explosion3D} illustrates the solution for the three-dimensional problem, for which again a level zero mesh with $34\times34\times34$ cells has been
used, together with $\mathfrak{r}=3$ and $\ell_{\max}=2$. This corresponds to an equivalent fine grid resolution of $306^3 = 28,652,616$ cells. In the three-dimensional case, 
a third order ADER--WENO scheme has been employed based on the Osher flux \eqref{eqn.osher}. The final AMR grid at time $t=0.25$ contains 9,079,984 elements and is depicted in Fig. \ref{fig_explosion3D.grid}. 
The simulation took 7.5 hours on 8 CPU cores of an AMD Opteron 6272 cluster with 2.1 GHz clock speed and 256 GB of RAM. One--dimensional cuts (along the positive $x$-axis) 
through the reconstruction polynomials $\mathbf{w}_h$ are shown on equidistant points in Fig. \ref{fig_explosion3D} together with the 1D reference solution according to \eqref{eq.1D}. 
An excellent agreement is observed. 

\paragraph*{Double Mach reflection problem.}
%
A classical test problem that contains simultaneously strong shock waves, contact waves and shear waves is represented by the double Mach reflection 
test \cite{woodwardcol84}, whose initial conditions are given by a right-moving shock wave at shock Mach number $M=10$, intersecting the $x-$ axis at  
$x=1/6$ with an inclination angle of $\alpha=60^{\circ}$. The computational domain is $\Omega=[0;3] \times [0;1]$. On the left, top and right boundary,  
the exact solution is prescribed, while a reflective wall boundary is imposed on the bottom. This test problem is frequently used in the literature 
on high order WENO and Discontinuous Galerkin schemes, see e.g. 
\cite{shu_efficient_weno,eno,balsarashu,SebastianShu,cbs4,WENOComplexFlows,titarevtoro,QiuDumbserShu,DGWENOLimiters,DumbserKaeser07}, where many 
reference solutions for this problem can be found. 

The initial condition for this problem is given by the Rankine-Hugoniot conditions as follows: 
\begin{equation}
\label{eqn.dmr.ic}
 \left( \rho, v_x, v_y, v_z, p \right)(x,y,0) =  \left\{ \begin{array}{cl} 
 \left(  8.0,  8.25 \cos(\alpha), 8.25 \sin(\alpha), 0.,  116.5  \right) 
  & \textnormal{ if } x' < 0.0, \\ 
 \left(  1.4,  0.0, 0.0., 0.,  1.0   \right) 
 & \textnormal{ if } x' \geq 0.0, \end{array} \right. 
\end{equation}
with $x' = (x - 1/6) \cos(\alpha) - y \sin(\alpha)$. 

The ratio of specific heats is chosen as $\gamma=1.4$. The problem is solved with a third order ADER-WENO scheme using the Rusanov flux \eqref{eqn.rusanov}. We use 
a mesh on the coarsest level consisting of only $150 \times 50$ elements, together with a refine factor of $\mathfrak{r}=4$ and a maximum 
refinement level of $\ell_{\max} = 2$. On the finest level, this corresponds to an effective resolution of $2400 \times 800$ control volumes. 

\begin{figure}
\begin{center}
\includegraphics[width=0.95\textwidth]{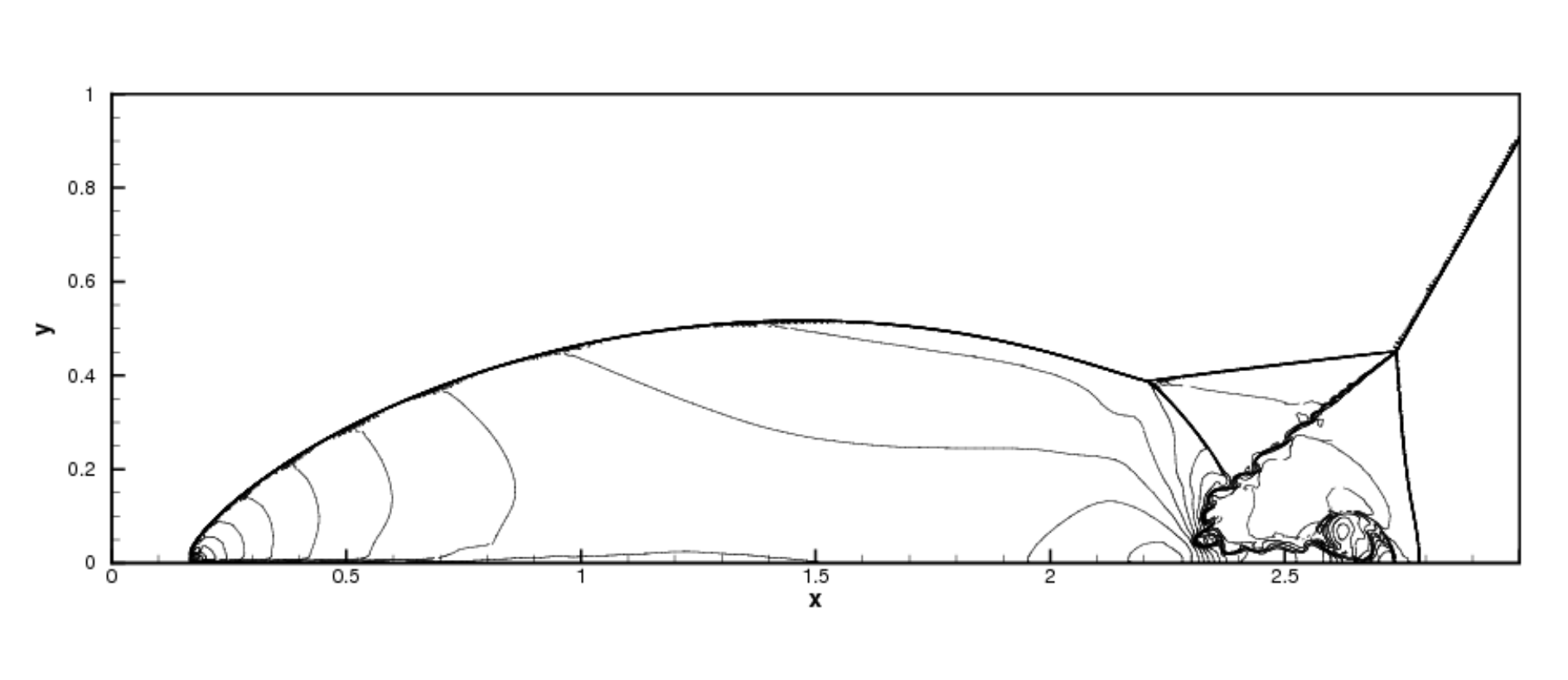} \\ 
\includegraphics[width=0.95\textwidth]{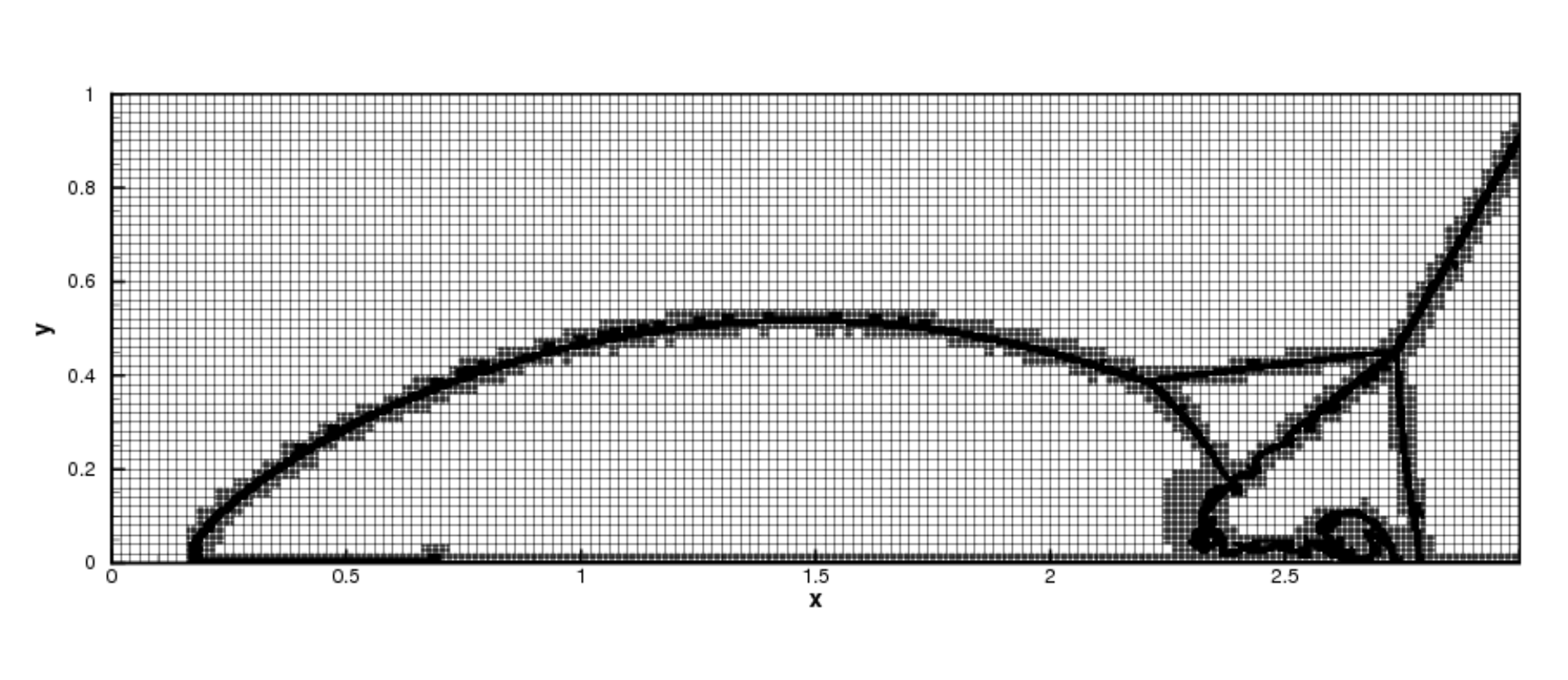}
\caption{ Double Mach reflection problem at time $t=0.2$. Top: equidistant density contour lines (contour spacing $\Delta \rho = 0.5$). 
Bottom: AMR grid with two levels of grid refinement.    
}
\label{fig.dmr}
\end{center}
\end{figure}

\begin{figure}
\begin{center}
\begin{tabular}{lr}
\includegraphics[width=0.48\textwidth]{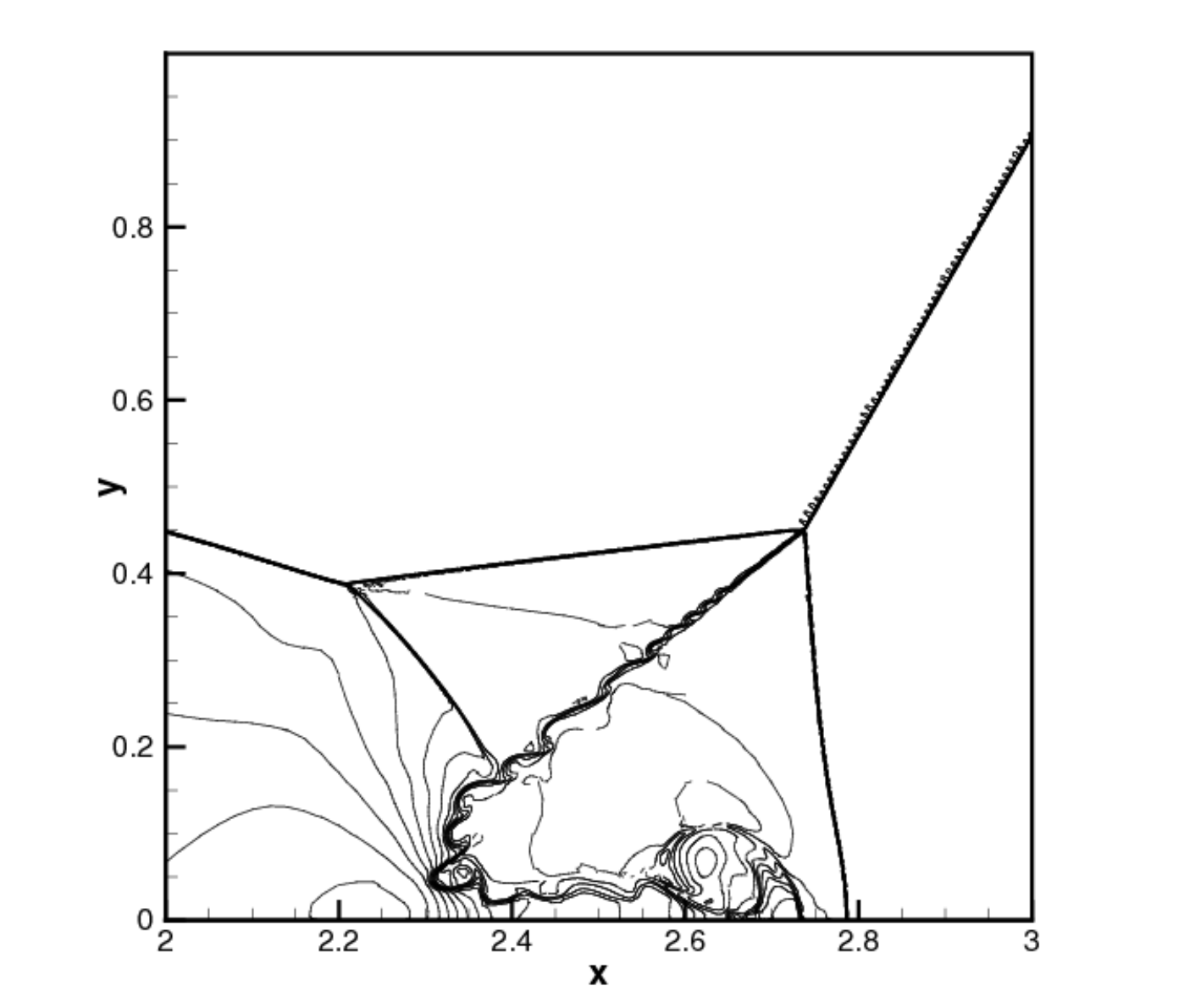} &  
\includegraphics[width=0.48\textwidth]{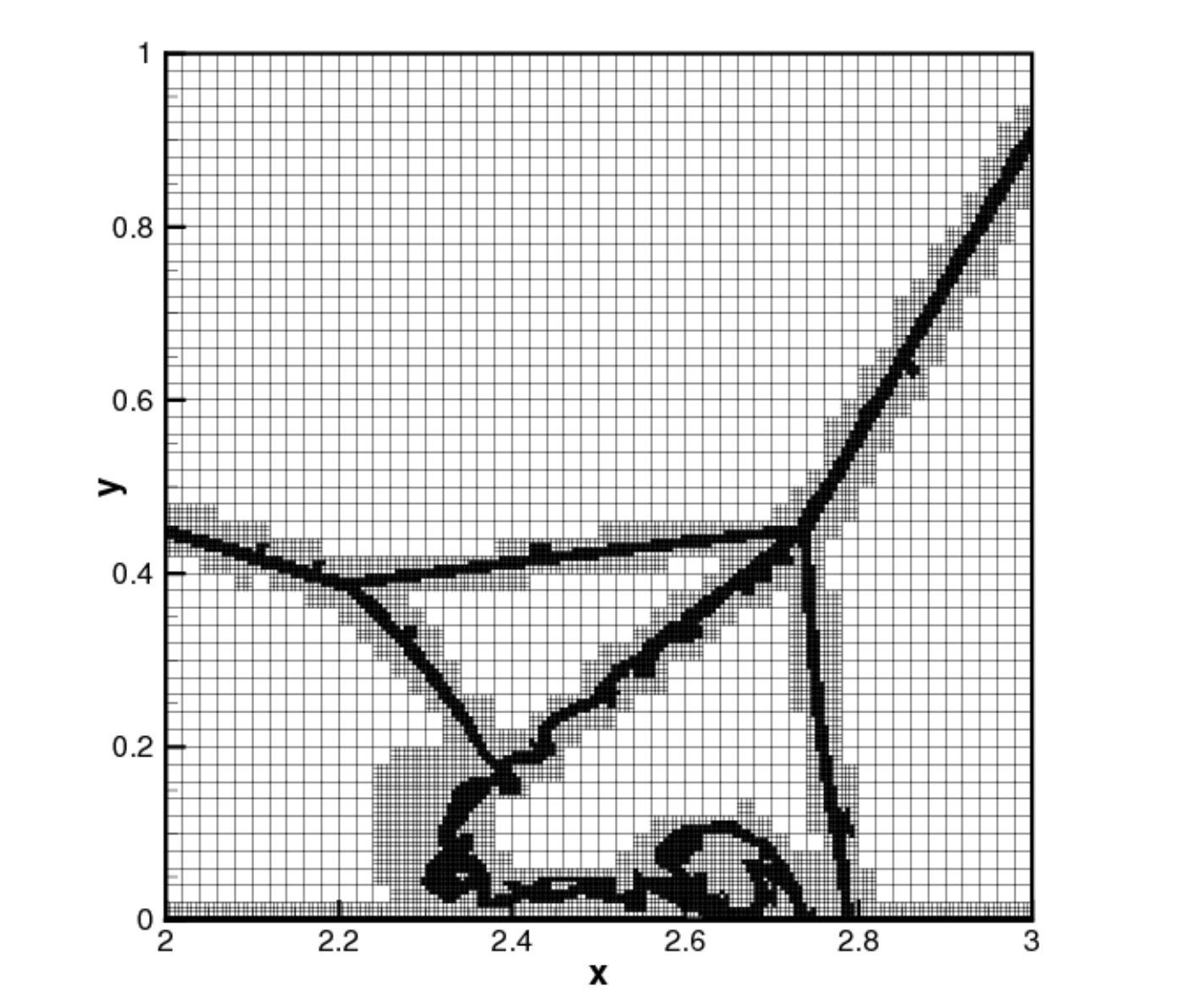}
\end{tabular} 
\caption{ Zoom into the double Mach reflection problem at time $t=0.2$. 
Left: equidistant density contour lines (contour spacing $\Delta \rho = 0.5$). 
Right: AMR grid with two levels of grid refinement.    
}
\label{fig.dmrz}
\end{center}
\end{figure}

The results for the density (31 equidistant contour levels from 1.5 to 22.5) are depicted at time $t=0.2$ in Figure \ref{fig.dmr} together with 
the final AMR grid. A zoom of the solution and the mesh is shown in Fig. \ref{fig.dmrz}. 

The shear waves present in this test are subject to the classical Kelvin--Helmholtz instability and therefore tend to roll up. Since there is no 
physical viscosity in the compressible Euler equations solved here, the developed small-scale flow features are purely governed by numerical 
viscosity. However, the amount of roll-up is a good qualitative indicator of the amount of numerical viscosity since more roll up indicates less 
numerical viscosity introduced by the scheme. The results obtained with the present ADER-WENO scheme on space-time adaptive Cartesian meshes is in 
good qualitative agreement with other published results for this test problem. 

\begin{figure}
\begin{center}
\includegraphics[width=0.95\textwidth]{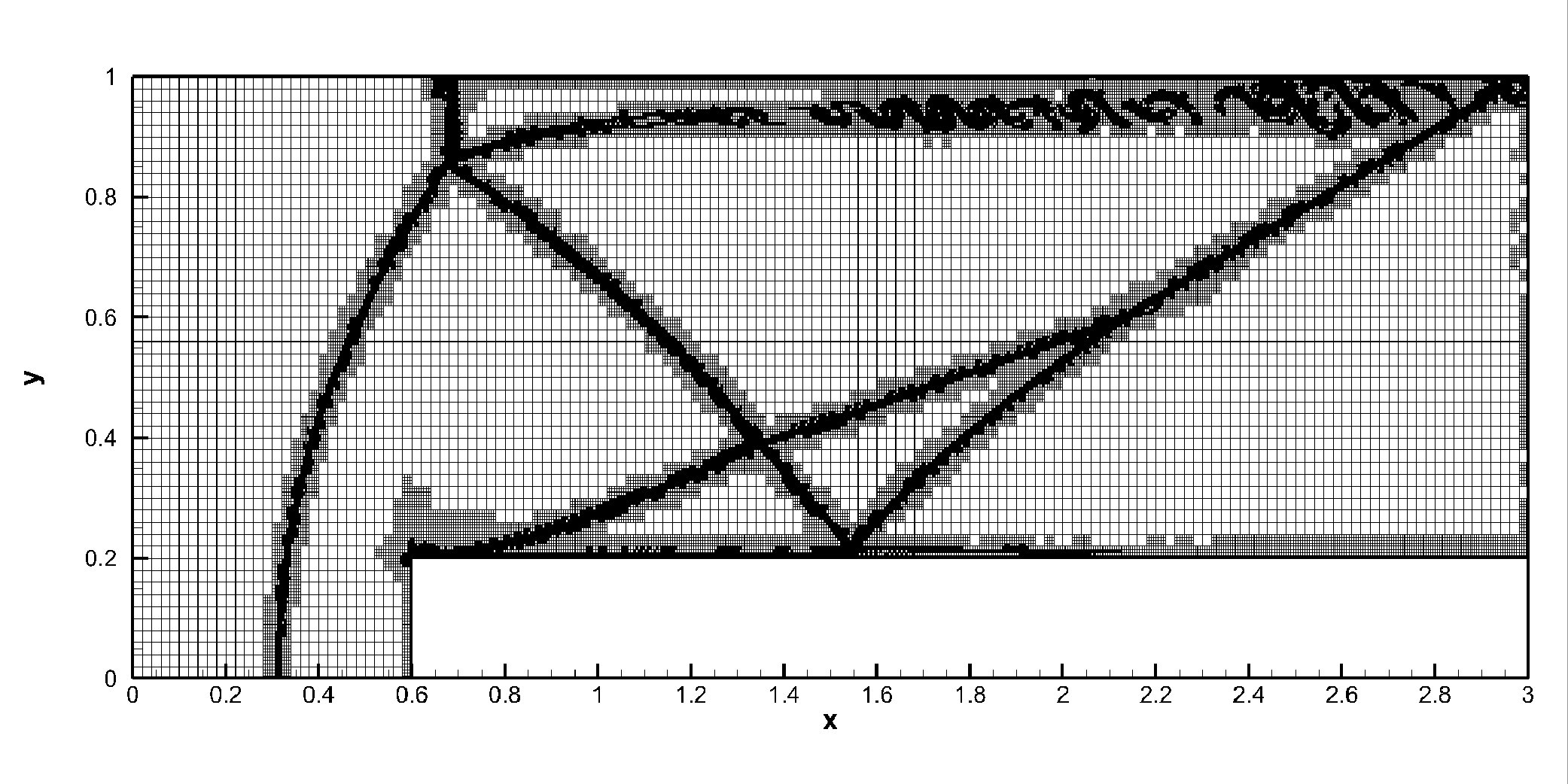} \\ 
\includegraphics[width=0.95\textwidth]{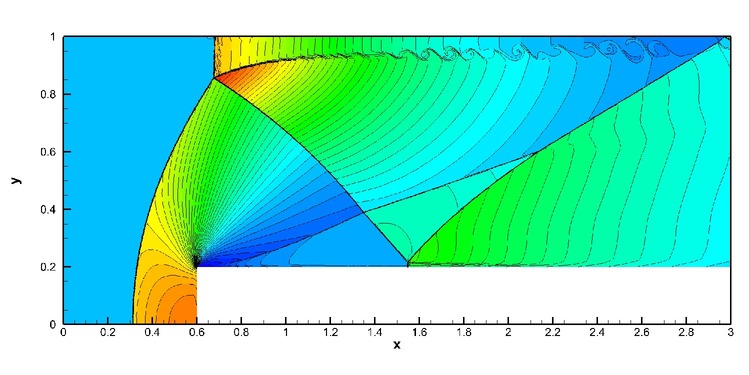} \\  
\includegraphics[width=0.95\textwidth]{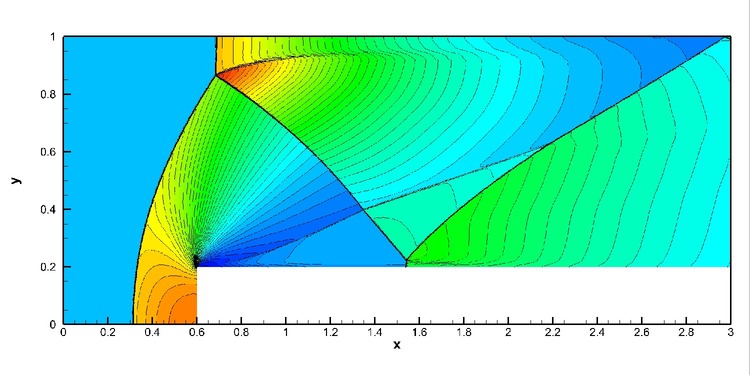}    
\caption{ Forward facing step problem. AMR grid (top) and density contours (center) obtained with the third order ADER-WENO scheme 
at time $t=2.5$. One can clearly observe the roll up of the slip lines. For comparison, also a second order solution is shown (bottom). 
}
\label{fig.ffs}
\end{center}
\end{figure}

\paragraph*{Forward facing step.} 
%
Another classical test problem for high resolution shock--capturing finite volume scheme consists in the forward facing step problem,
also called the Mach 3 wind tunnel test.   
It has also been proposed originally in \cite{woodwardcol84}. The computational domain is given by 
$\Omega = [0;3] \times [0;1] \backslash [0.6;3] \times [0;0.2]$ and the initial condition is a uniform flow at Mach number $M=3$ 
moving to the right. In particular, we use $\rho(x,y,0) = 1$, $p(x,y,0) = 1 / \gamma$, $v_x(x,y,0) = 3$ and $v_y=v_z=0$. The ratio 
of specific heats is set to $\gamma=1.4$. Simulations are carried out until $t=2.5$. Reflective boundary conditions are applied 
on the upper and lower boundary of the domain and inflow/outflow boundary conditions are applied at the entrance/exit. At the 
corner of the step, there is a singularity, which is properly resolved with the third order ADER-WENO scheme using adaptive mesh 
refinement. The mesh on the coarsest level contains $150 \times 50$ control volumes. We use $\mathfrak{r}=4$ and $\ell_{\max}=2$, 
hence on the finest level this corresponds to an equivalent resolution of $2400 \times 800$. 
The computational results obtained with the third order ADER-WENO method as well as a sketch of the final AMR mesh are depicted in Fig. 
\ref{fig.ffs}. For comparison, also a second order simulation is shown. One can clearly observe that the third order scheme provides 
a much better resolution of the physical instability and roll up of the contact line compared to the standard second order scheme. This 
indicates that even in the context of space--time adaptive mesh refinement, the use of higher order schemes may be appropriate to enhance
resolution and to reduce numerical viscosity for small scale turbulent structures.    
 
\begin{figure}
\begin{center}
\begin{tabular}{lr}
\includegraphics[width=0.4\textwidth]{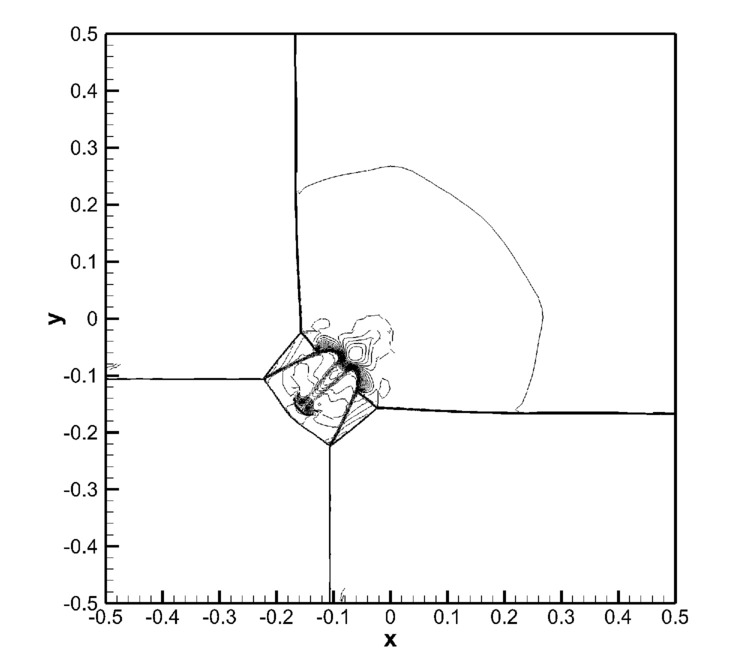}     &  
\includegraphics[width=0.4\textwidth]{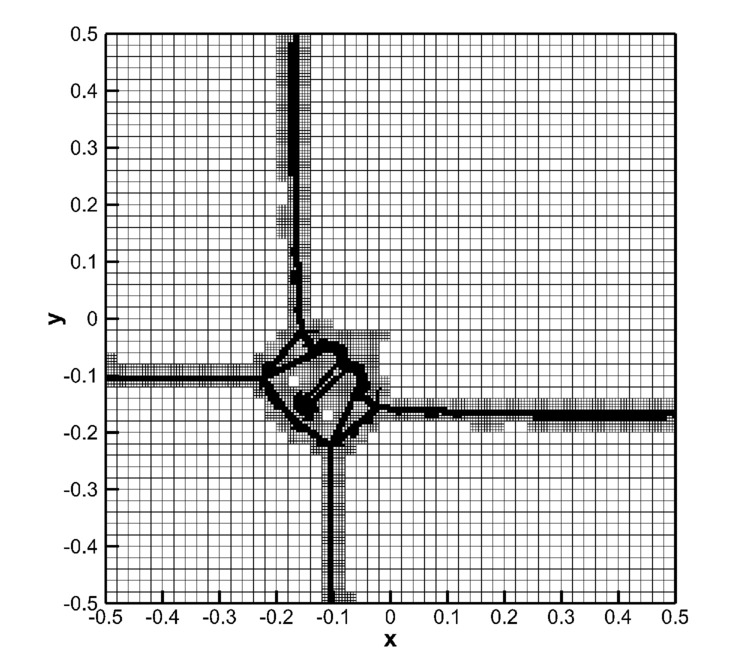} \\
\includegraphics[width=0.4\textwidth]{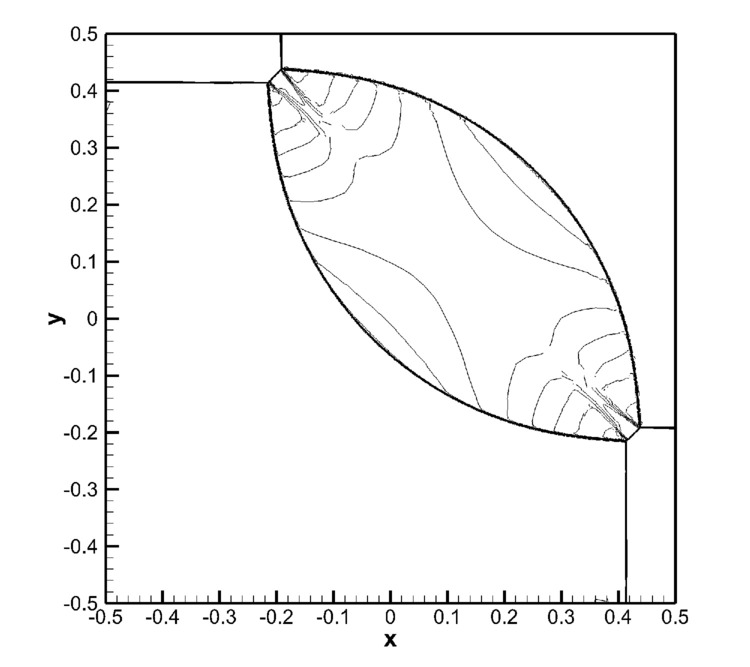}     &  
\includegraphics[width=0.4\textwidth]{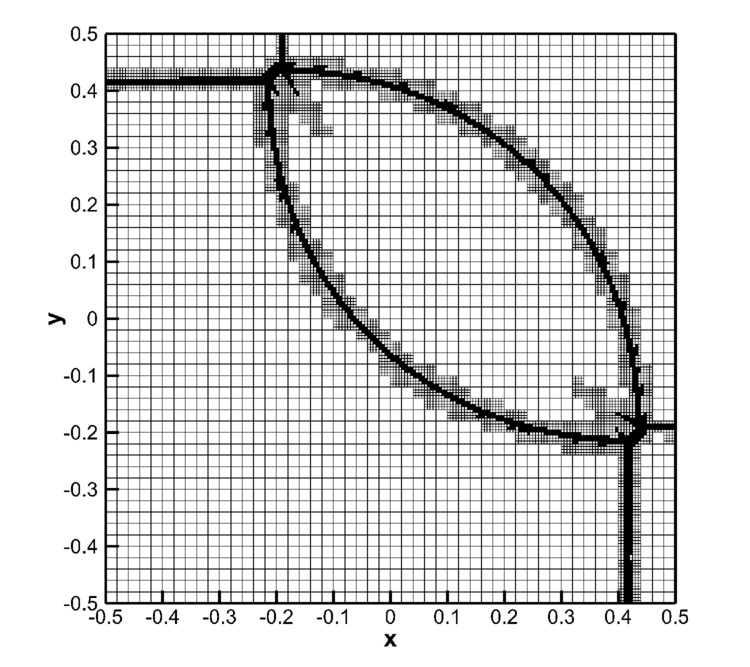} \\
\includegraphics[width=0.4\textwidth]{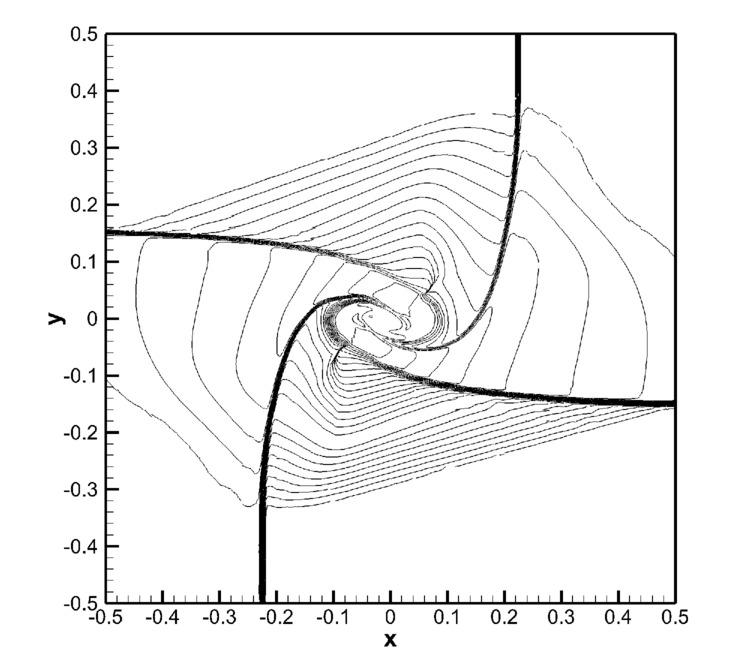}     &  
\includegraphics[width=0.4\textwidth]{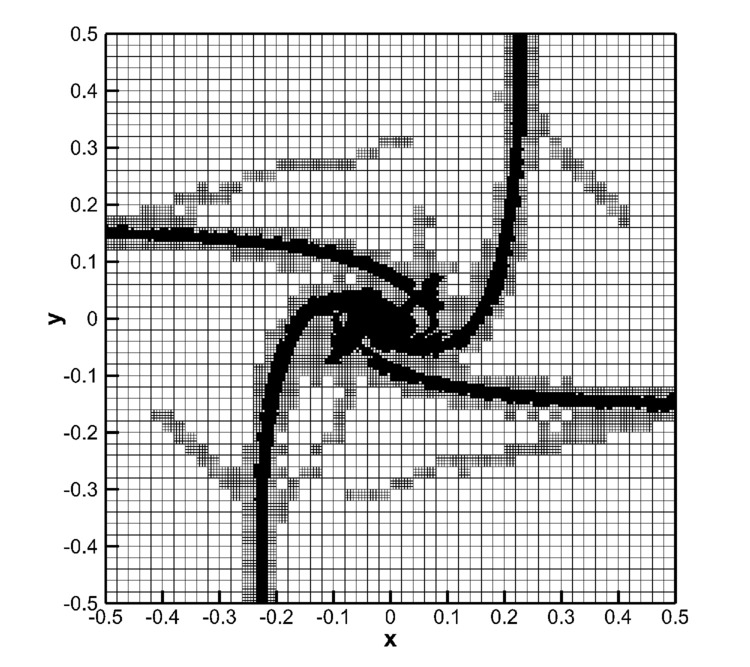} \\
\includegraphics[width=0.4\textwidth]{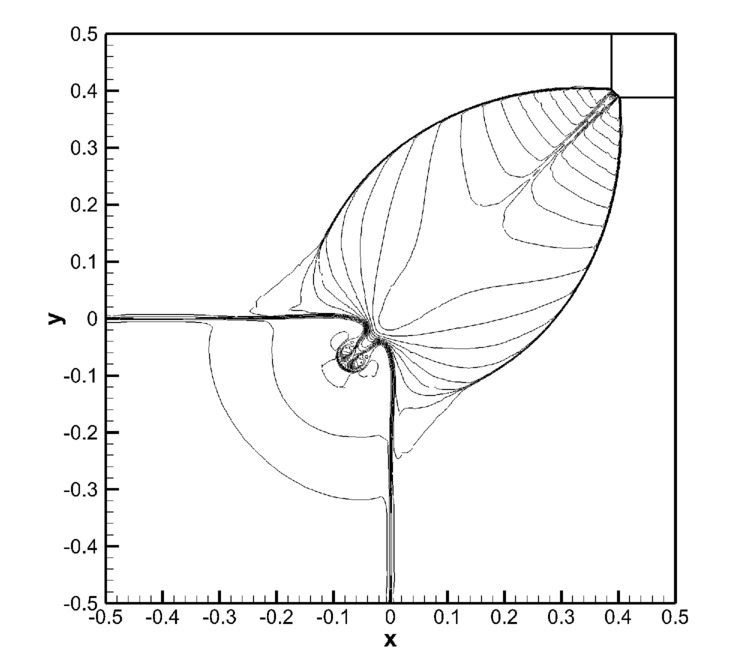}     &  
\includegraphics[width=0.4\textwidth]{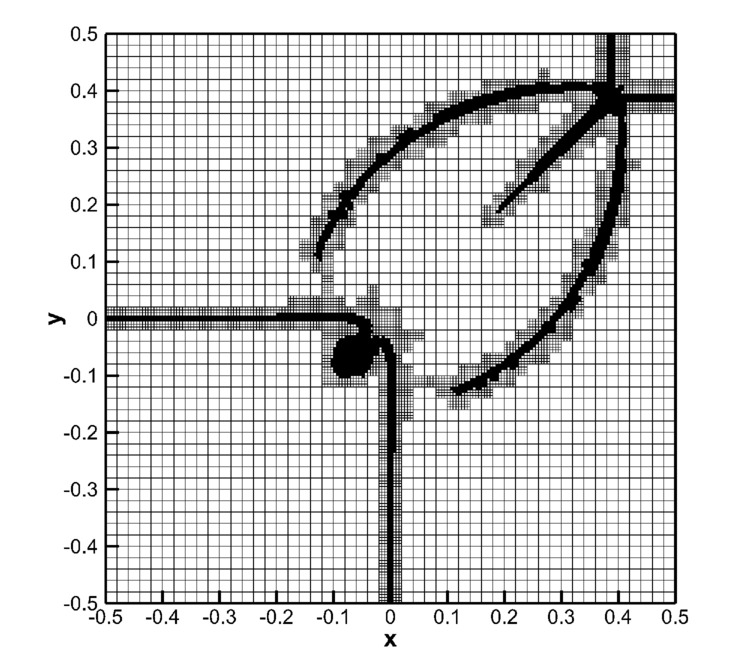}  
\end{tabular} 
\caption{ Two-dimensional Riemann problems solved with third order ADER-WENO schemes. 
Density contour lines (left) and AMR grid at the final time (right). 
}
\label{fig.rp2d}
\end{center}
\end{figure}

\paragraph*{2D Riemann problems.} 
%
A large set of two--dimensional Riemann problems has been cataloged in \cite{kurganovtadmor}. The computational domain is $\Omega = [-0.5;0.5] \times [-0.5;0.5]$ 
and the initial conditions are given by 
\begin{equation}
 \mathbf{u}(x,y,0) = \left\{ \begin{array}{ccc} 
 \mathbf{u}_1 & \textnormal{ if } & x > 0 \wedge y > 0,    \\ 
 \mathbf{u}_2 & \textnormal{ if } & x \leq 0 \wedge y > 0, \\ 
 \mathbf{u}_3 & \textnormal{ if } & x \leq 0 \wedge y \leq 0, \\ 
 \mathbf{u}_4 & \textnormal{ if } & x > 0 \wedge y \leq 0.    
  \end{array} \right. 
\end{equation} 
The initial conditions and the final simulation time $t_f$ for the four configurations presented in this article are listed in Table \ref{tab.rp2d.ic}. In 
all cases $\gamma=1.4$. The simulations are carried out with a third order one--step ADER WENO scheme using a level zero grid of $50 \times 50$ elements. 
The computational results together with the final AMR grids are depicted in Fig. \ref{fig.rp2d}. We can note a good agreement 
with the reference solution published in \cite{kurganovtadmor}. 

\begin{table}[!b]   
\caption{Initial conditions for the two--dimensional Riemann problems.} 
\begin{center} 
\renewcommand{\arraystretch}{1.0}
\begin{tabular}{c|cccc|cccc} 
\hline
\multicolumn{9}{c}{Problem RP1 (Configuration 3 in \cite{kurganovtadmor}), $t_f = 0.25$} \\ 
\hline
     & \multicolumn{4}{c|}{$x \leq 0$} & \multicolumn{4}{|c}{$x>0$} \\
\hline
     & $\rho$ & $u$ & $v$  & p & $\rho$ & $u$ & $v$  & p  \\ 
\hline
$y > 0$    & 0.5323 & 1.206   & 0.0     & 0.3   & 1.5    & 0.0 &  0.0    & 1.5 \\ 
$y \leq 0$ & 0.138  & 1.206   & 1.206   & 0.029 & 0.5323 & 0.0 &  1.206  & 0.3 \\ 
\hline
\multicolumn{9}{c}{Problem RP2 (Configuration 4 in \cite{kurganovtadmor}), $t_f = 0.25$} \\
\hline
     & \multicolumn{4}{c|}{$x \leq 0$} & \multicolumn{4}{|c}{$x>0$} \\
\hline
     & $\rho$ & $u$ & $v$  & p & $\rho$ & $u$ & $v$  & p  \\ 
\hline
$y > 0$    & 0.5065 &  0.8939 & 0.0     & 0.35 & 1.1    & 0.0 &  0.0    & 1.1  \\ 
$y \leq 0$ & 1.1    &  0.8939 & 0.8939  & 1.1  & 0.5065 & 0.0 &  0.8939 & 0.35 \\ 
\hline
\multicolumn{9}{c}{Problem RP3 (Configuration 6 in \cite{kurganovtadmor}), $t_f = 0.30$} \\
\hline
     & \multicolumn{4}{c|}{$x \leq 0$} & \multicolumn{4}{|c}{$x>0$} \\
\hline
     & $\rho$ & $u$ & $v$  & p & $\rho$ & $u$ & $v$  & p  \\ 
\hline
$y > 0$    & 2.0    &  0.75  & 0.5   & 1.0  & 1.0  &  0.75 &  -0.5  & 1.0  \\ 
$y \leq 0$ & 1.0    & -0.75  & 0.5   & 1.0  & 3.0  & -0.75 &  -0.5  & 1.0  \\ 
\hline
\multicolumn{9}{c}{Problem RP4 (Configuration 12 in \cite{kurganovtadmor}), $t_f = 0.25$} \\
\hline
     & \multicolumn{4}{c|}{$x \leq 0$} & \multicolumn{4}{|c}{$x>0$} \\
\hline
     & $\rho$ & $u$ & $v$  & p & $\rho$ & $u$ & $v$  & p  \\ 
\hline
$y > 0$    & 1.0    &  0.7276 & 0.0     & 1.0  & 0.5313 & 0.0 &  0.0    & 0.4  \\ 
$y \leq 0$ & 0.8    &  0.0    & 0.0     & 1.0  & 1.0    & 0.0 &  0.7276 & 1.0  \\ 
\hline 
\end{tabular} 
\end{center}
\label{tab.rp2d.ic}
\end{table} 

%
\paragraph*{A co-rotating vortex pair.} 
%
This multi-scale problem from aeroacoustics solved here differs from the previous ones in several aspects. First, it 
is a low Mach number problem without shock waves. Second, it is of true multi-scale nature. The problem is 
an academic example of sound generation mechanisms in aeroacoustics and is taken from 
\cite{Mitchelletal,leekoo,dahl,Diss-Thomas,DomainDecomp,MPVEIFLEE}. It consists of two isentropic  vortices with 
characteristic size $r_c$ (vortex core radius) that rotate around each other, thus generating sound waves with a 
wavelength that is about two orders of magnitude larger than the size of the vortices themselves. While the sound
waves are very smooth, the vortices contain strong gradients in velocity and pressure as one approaches the vortex 
core. 
The accurate propagation of sound waves with little amplitude and phase errors in large domains is 
the major challenge in computational aeroacoustics. For such problems, typically special high order 
schemes are employed, e.g. \cite{lele,TamDRP,tam,schwartzkopff-dumbser-munz,euromech03}, since the 
conventional high resolution TVD schemes that are typically used in classical AMR codes are too diffusive 
and lose accuracy at local  extrema. Since such extrema regularly occur in acoustic wave propagation 
problems the use of high order WENO schemes, as the ones used in this paper, that are able to handle strong 
gradients without degenerating at local extrema, seems to be appropriate. 

The initial condition for the velocity field is given by the superposition of the velocity fields induced by two 
potential vortices. The complex potential $w$ of the rotating vortex pair is given by 
\begin{equation}
  w(z,t) = \frac{\Gamma}{2\pi i} \ln \left( z^2 - b^2 \right), 
  \label{eqn.potential} 
\end{equation} 
with $z=x + i y$, $b = r_0 e^{i \omega t}$ and $ i^2 = -1$. The circulation of each vortex is denoted by $\Gamma$, 
the angular rotation frequency of the vortex pair is given by $\omega = \Gamma / (4 \pi r_0^2)$ and the rotation 
Mach number is $M = \Gamma/(4 \pi r_0 c_0)$, with the usual definition of the sound speed as 
$c_0^2 = \gamma p_0 / \rho_0$. The ambient reference density and pressure are denoted by $\rho_0$ and 
$p_0$, respectively. 
From \eqref{eqn.potential} one obtains the Cartesian velocity components $u_x$ and $u_y$ as  
\begin{equation}
 \frac{\partial w}{\partial z} = \frac{\Gamma}{\pi i} \frac{z}{z^2 - b^2}  = v_x - i v_y.  
\end{equation} 
The hydrodynamic pressure associated with the vortex pair is given by the unsteady Bernoulli equation as 
\begin{equation} 
 p = p_0 - \rho_0 \left( \textnormal{Re}\left(\frac{\partial w}{\partial t}\right) + \frac{1}{2}(v_x^2 + v_y^2) \right).
\end{equation} 
The initial density is defined by $\rho = p^{1/\gamma}$. 
Inside the vortex core radius $r_c$, i.e. for $r < r_c$, a Gaussian-type vorticity distribution is imposed according
to \cite{Mitchelletal,Diss-Thomas}, in order to avoid the singularity of the potential flow in the vortex center. With 
a matched asymptotic expansion technique \cite{leekoo,Mitchelletal}, the far field sound pressure produced by the ideal 
incompressible pair of potential vortices can be obtained in polar coordinates $r^2=x^2+y^2$ and $\tan(\theta) = x/y$ 
as 
\begin{equation} 
 p' = \frac{\rho_0 \Gamma^4}{64 \pi^3 r_0^4 c_0^2} 
       \left(   \textnormal{J}_2(kr) \sin\left( 2(\omega t - \theta) \right) 
              - \textnormal{Y}_2(kr) \cos\left( 2(\omega t - \theta) \right) \right), 
 \label{eqn.spl} 
\end{equation}
with $k=2\omega/c_0$. $\textnormal{J}_2(kr)$ and  $\textnormal{Y}_2(kr)$ are the second order Bessel functions of the
first and second kind, respectively, and $p'$ denotes the fluctuation of the sound pressure about the unperturbed ambient 
mean pressure $p_0$. 

For the numerical simulations, the two-dimensional computational domain of this problem is chosen as $\Omega = [-500;500] \times [-500;500]$ 
and the entire problem is solved with the same fourth order ADER-WENO scheme using a level zero grid of $250 \times 250$ elements, together 
with $\mathfrak{r}=4$ and $\ell_{\max}=3$. On a uniform fine grid this would  correspond to an effective resolution of $16000 \times 16000$ 
mesh points, hence the use of an AMR technique with local time stepping or at least a suitable domain decomposition with local time stepping 
such as the one presented in  \cite{DomainDecomp} is mandatory. For comparison, also a second order AMR simulation with $\mathfrak{r}=4$ and 
$\ell_{\max}=3$ is run with $500 \times 500$ elements on the level zero grid. 

%

\begin{figure}
\begin{center}
\includegraphics[width=0.75\textwidth]{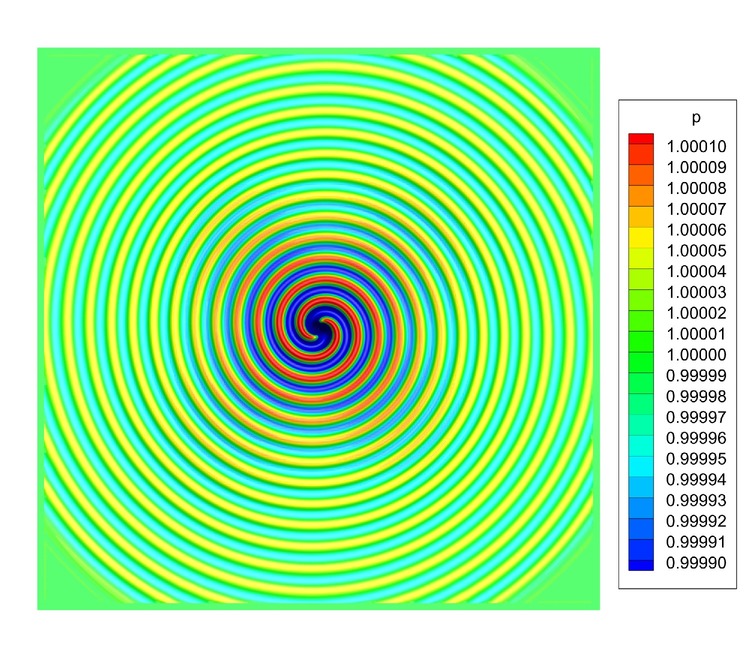}        
\caption{ Sound pressure field generated by the co-rotating vortex pair at time $t=500$. }
\label{fig.vortex.p}
\end{center}
\end{figure}

\begin{figure}
\begin{center}
\begin{tabular}{lr} 
\includegraphics[width=0.45\textwidth]{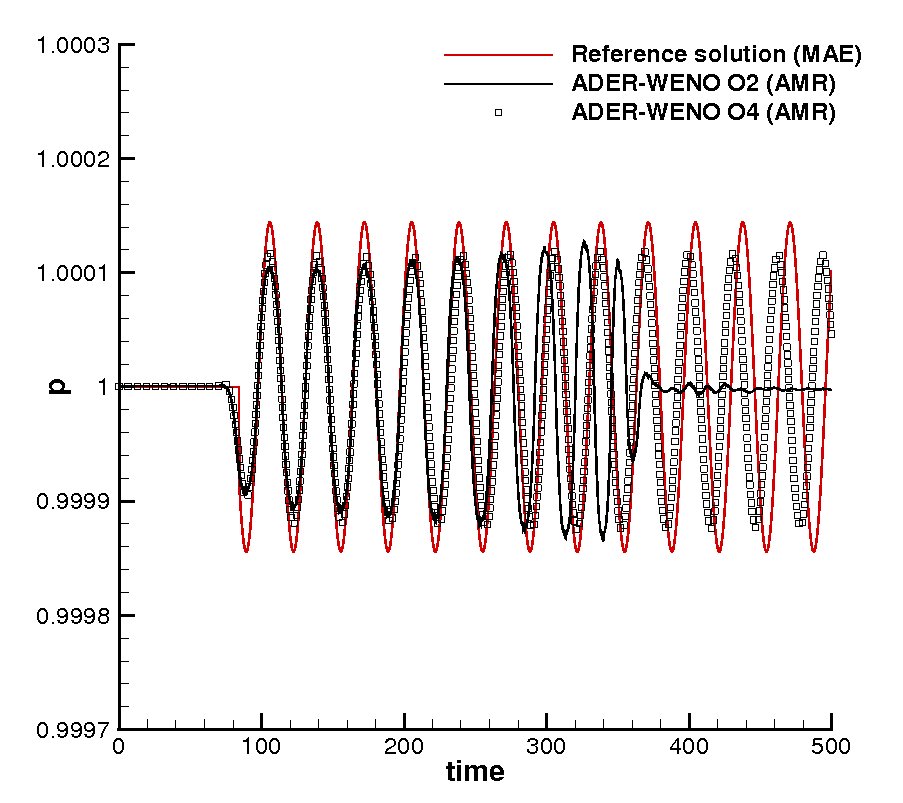}     &  
\includegraphics[width=0.45\textwidth]{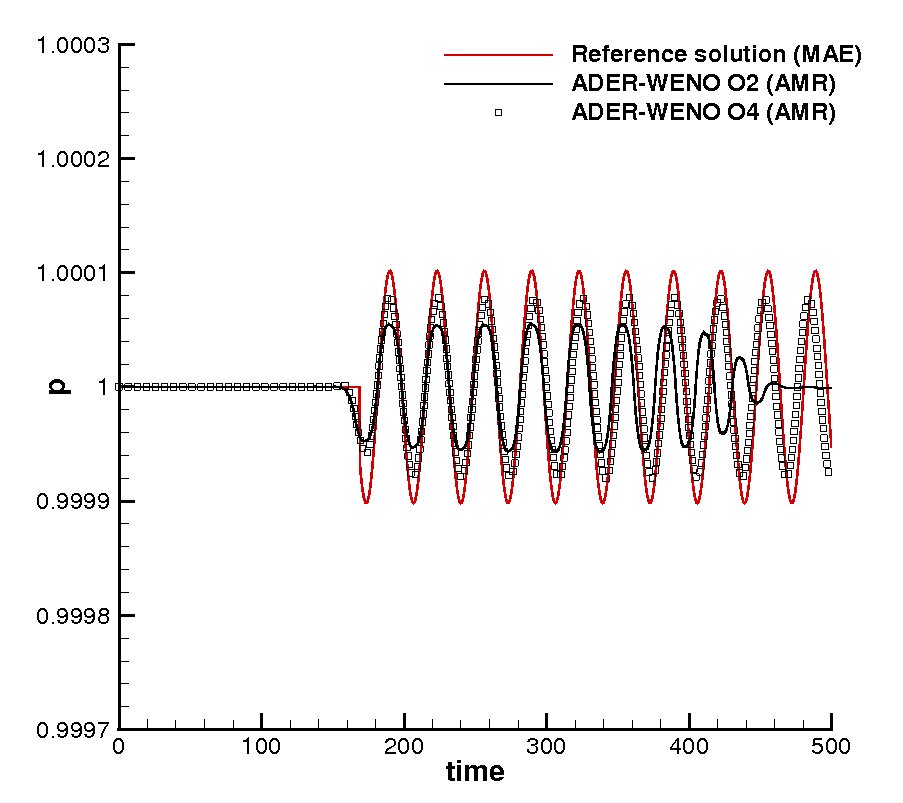}        
\end{tabular} 
\caption{ Temporal evolution of the sound pressure in the points $\mathbf{x}_1=(100,0)$ (left) and $\mathbf{x}_2=(200,0)$ (right). 
Comparison of the second and fourth order ADER-WENO AMR results with the matched asymptotic expansion (MAE) solution for the 
far field sound pressure generated by an ideal incompressible co-rotating potential vortex pair. The acoustic signal of the 
second order scheme ceases due to a spurious unphysical vortex merging caused by excessive numerical diffusion.} 
\label{fig.vortex.compare}
\end{center}
\end{figure}

\begin{figure}
\begin{center}
\begin{tabular}{lcr} 
\includegraphics[width=0.32\textwidth]{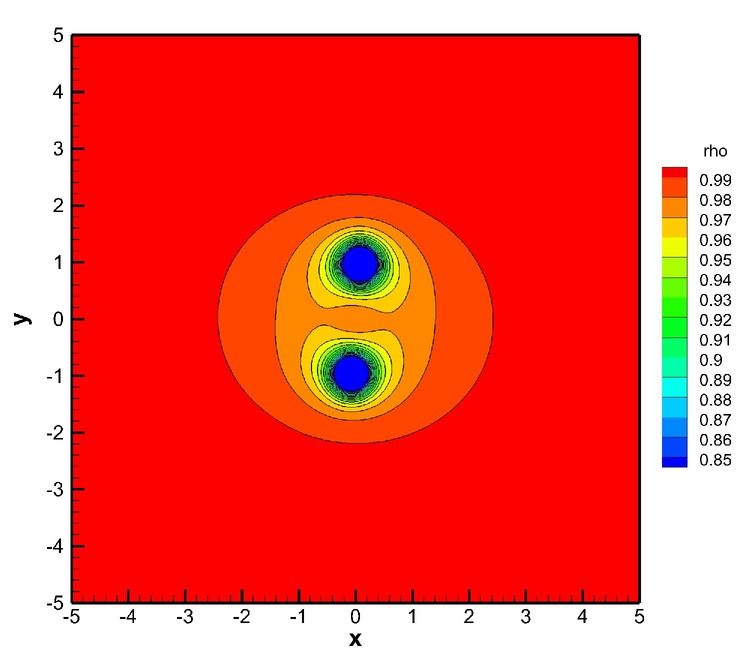}     &  
\includegraphics[width=0.32\textwidth]{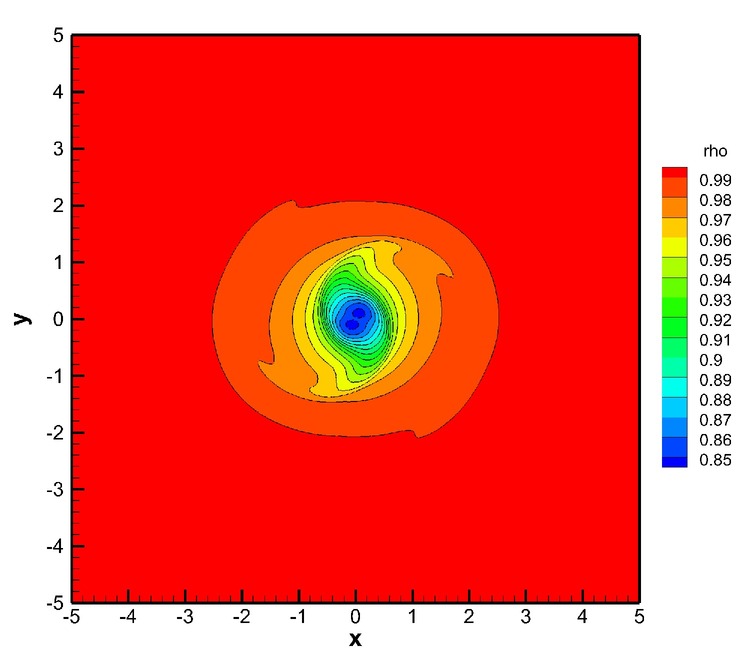}         & 
\includegraphics[width=0.32\textwidth]{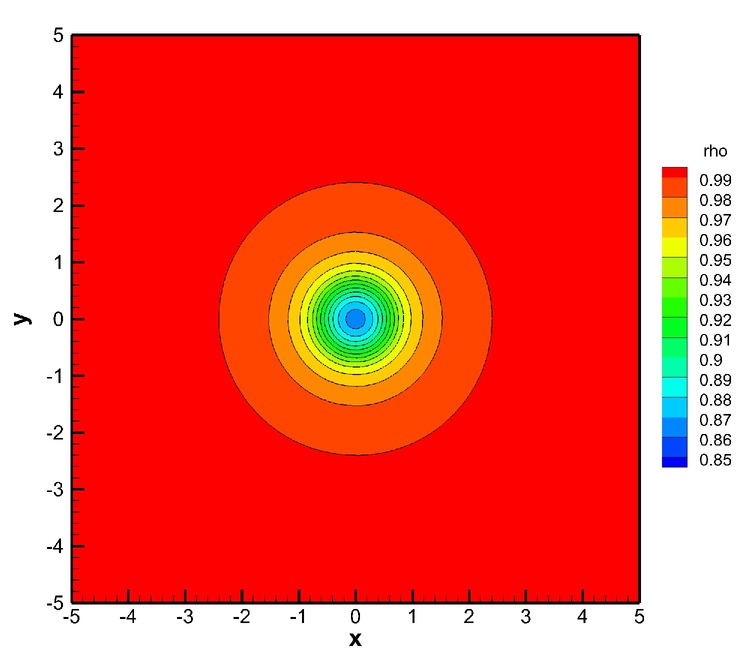}      
\end{tabular} 
\caption{ Density contours of the co-rotating vortex pair at time $t=500$ for the fourth order scheme (left) 
and for the second order scheme (right). The spurious vortex merging obtained with the second order scheme 
is also depicted for $t=281$ (center). } 
\label{fig.vortex.merge}
\end{center}
\end{figure}

The parameters used for this simulation are $r_c=0.2$, $\gamma=1.4$, $p_0=\rho_0=1$, 
$c_0 = \sqrt{\gamma p_0 / \rho_0}=\sqrt{\gamma}$, and $\Gamma = 0.08 \cdot 4 \pi \sqrt{\gamma}$, 
hence the rotation Mach number is $M=0.08$. With the above parameters the wave length of the sound waves 
$\lambda_s$ can be computed from \eqref{eqn.spl} as $\lambda_s = \pi c_0 / \omega \approx 39$. 
With the chosen grid resolution ($\Delta x = \Delta y = 4$ in the far field), the fourth order scheme resolves the 
acoustic waves with about 10 points per wavelength (PPW), while the second order scheme employs about 20 PPW.  
Simulations are performed until $t=500$, before the acoustic waves reach the corners of the outer border. The 
acoustic pressure field generated by the co-rotating vortex pair is shown in Fig. \ref{fig.vortex.p}. A comparison 
of our numerical simulations with the reference solution \eqref{eqn.spl} is depicted in Fig. \ref{fig.vortex.compare}. 
The reference solution \eqref{eqn.spl} is a time periodic solution. However, it is obvious for the present problem 
that when starting from an initially undisturbed pressure field, no sound signal can arrive at a given spatial point 
$\mathbf{x}=(x,y)$ before the time $|\mathbf{x}|/c_0$, hence the analytical reference solution is depicted in 
Fig. \ref{fig.vortex.compare} only for times larger than $|\mathbf{x}|/c_0$. Note further that in 
the present simulations the entire problem has been solved using the \textit{compressible} Euler equations from the 
near  field up to the very far field and that the singularity in the center of the potential vortices has been avoided 
by a Gaussian-type vorticity distribution inside the vortex core, as suggested in \cite{Mitchelletal,Diss-Thomas}. 
In contrast, the reference solution has been obtained with a matched asymptotic expansions technique for the radiated 
sound field of a pair of ideal incompressible potential vortices. Due to the modified flow field inside the
core radius with respect to the ideal potential vortex, we expect our sound wave amplitudes to be always lower than 
the ones of the reference solution, which is actually confirmed by the results shown in Fig. \ref{fig.vortex.compare}. 
With this said we can observe an overall good agreement with the analytical solution concerning phase and amplitude of 
the sound pressure signal for the fourth order ADER-WENO scheme. For comparison, a numerical solution obtained with 
a classical second order AMR scheme on a grid refined twice as much is also shown in Fig. \ref{fig.vortex.compare}. 
The mesh refinement for the second order scheme with respect to the fourth order method leads exactly to the same 
number of degrees of freedom used to represent the reconstructed solution $\mathbf{w}_h$. The second order scheme has to update four 
times more cells and due to the CFL condition also needs twice as many time steps compared to the fourth order scheme, 
which increases the number of zone updates by a factor of eight. Since the second order AMR scheme is 9.4 times cheaper 
per element update compared to the fourth order AMR method, see Table \ref{tab.efficiency}, the total CPU times of both 
simulations are comparable. 
However, due to the significantly higher numerical diffusion of the second order method even on the refined mesh, an  
\textit{unphysical} vortex merging appears, which causes the acoustic signal to cease completely after a certain time,  
since the merged vortices collapse into a single stationary vortex, which does not emit any sound waves. This is
clearly seen in the acoustic signals of Fig. \ref{fig.vortex.compare}, which also show that the second order scheme 
obtains much lower sound pressure amplitudes in the second point $\mathbf{x}_2=(200,0)$ even before the unphysical
vortex merging. The second observation point is about five propagated wavelengths away from the center of the vortex pair. 
We furthermore show the vortex configuration at the final time $t=500$ for both the fourth and the second order scheme in Fig. 
\ref{fig.vortex.merge}, as well as the time $t=281$ when the spurious numerical vortex merging takes place for the second 
order scheme. For a detailed study of \textit{physical} vortex mergers, see \cite{Melander,Waugh}. 

%
\subsection{Classical MHD equations}

In this section we consider a more complicated hyperbolic system than the Euler equations used in the 
previous section. We solve the classical, i.e. non--relativistic, equations of ideal magnetohydrodynamics (MHD) 
in three space dimensions. The MHD system introduces an additional difficulty for numerical schemes since the divergence 
of the magnetic field must remain zero for all times, i.e.  
\begin{equation}
\label{eqn.divfree}
  \frac{\partial B_x}{\partial x} +  
  \frac{\partial B_y}{\partial y} +  
  \frac{\partial B_z}{\partial z} = 0,  
\end{equation}
which for the continuous problem is always satisfied under the condition that the initial data of the magnetic field are 
divergence-free. From the discrete point of view this is not necessarily guaranteed and hence extra care is required in the 
discretization. In this article we use the hyperbolic version of the generalized Lagrangian multiplier (GLM) divergence cleaning 
approach proposed in \cite{Dedneretal}. It consists in adding an auxiliary variable $\Psi$ and one linear scalar PDE to the MHD 
system to transport divergence errors out of the computational domain with the artificial speed $c_h$. The augmented MHD system 
with hyperbolic GLM divergence cleaning has the state vector $\mathbf{u}$ given by 
\begin{equation}
\mathbf{u}^T=\left( \rho \,\, \rho \vec v^T  \,\, E \,\, \vec B^T \,\, \psi \right), 
\end{equation} 
and the flux tensor $\mathbf{F}=(\mathbf{f}, \mathbf{g}, \mathbf{h})$ is defined as: 
\begin{equation}
\mathbf{F} = \left( \begin{array}{c} 
    \rho \vec{v}^T \\ \rho \vec{v} \vec{v} + (p +  \frac{1}{8\pi} \vec B^2 )\, \mathbf{I} - \frac{1}{4\pi} \vec{B} \vec {B} \\
    \vec v^T ( E + p + \frac{1}{8\pi} \vec B^2) - \frac{1}{4\pi} \vec B^T (\vec v \cdot \vec B) \\
    \vec v \vec B - \vec B \vec v + \Psi \mathbf{I}, \\ c_h^2 \vec B^T \end{array} \right). 
\end{equation} 
with the velocity vector $\vec v = (v_x,v_y,v_z)^T$, the magnetic field vector $\vec B =(B_x,B_y,B_z)^T$ and the $3 \times 3$ identity 
matrix $\mathbf{I}$. The equation of state is the ideal gas law, hence 
\begin{equation}
   p = (\gamma-1)( E - \frac{1}{2} \rho \vec v^2 - \frac{\vec B^2}{8\pi} ).
\end{equation}

\begin{figure}
\begin{center}
\begin{tabular}{lcr}
\includegraphics[width=0.31\textwidth]{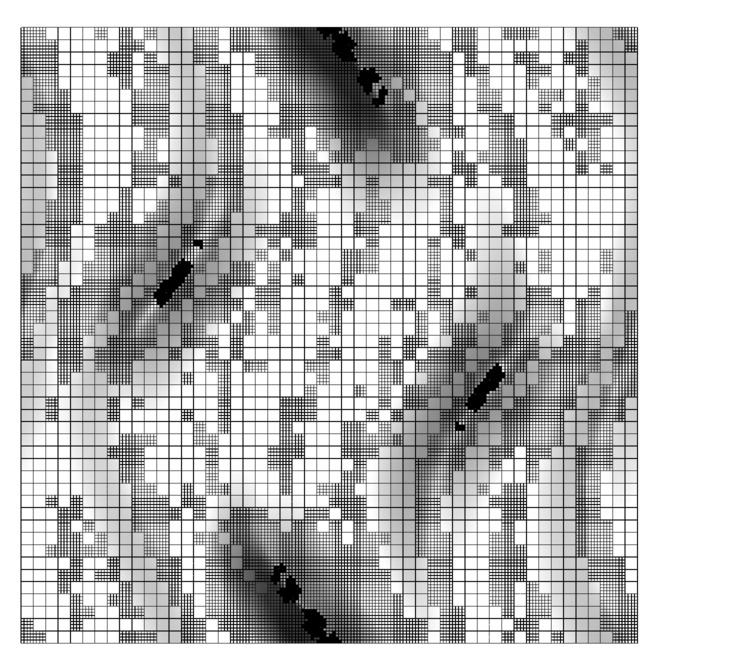}     &  
\includegraphics[width=0.31\textwidth]{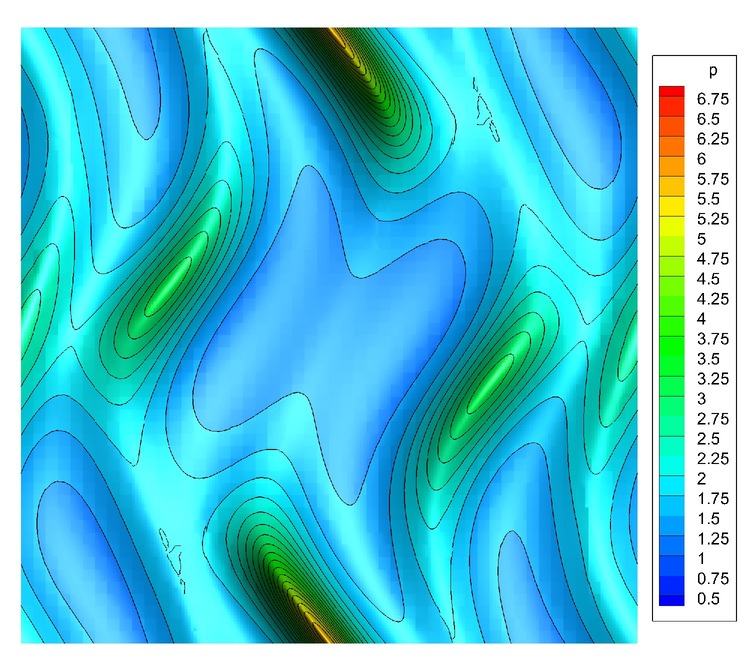}      &  
\includegraphics[width=0.31\textwidth]{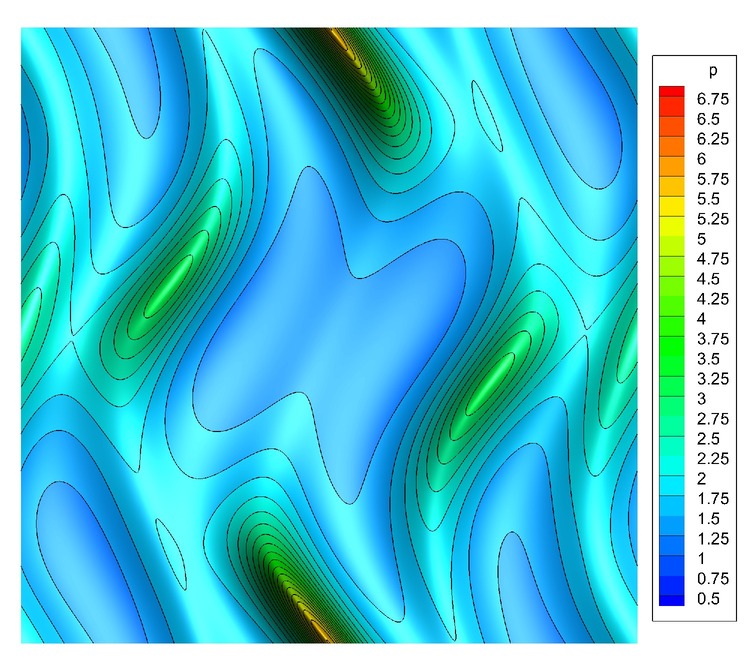}  \\   
\includegraphics[width=0.31\textwidth]{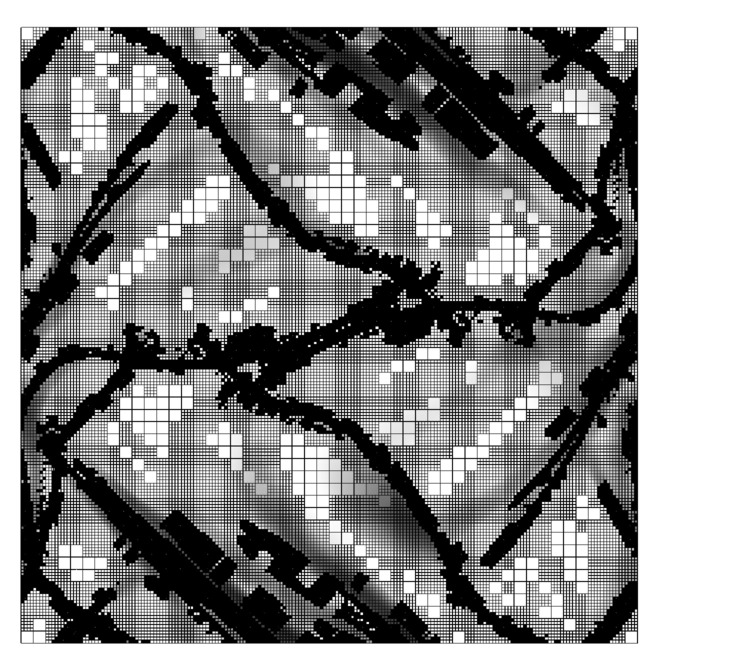}     &  
\includegraphics[width=0.31\textwidth]{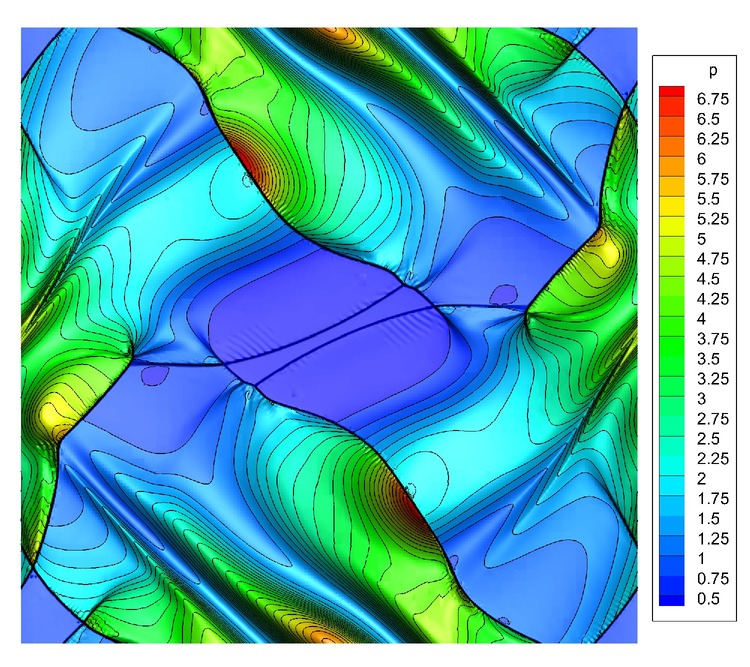}      &  
\includegraphics[width=0.31\textwidth]{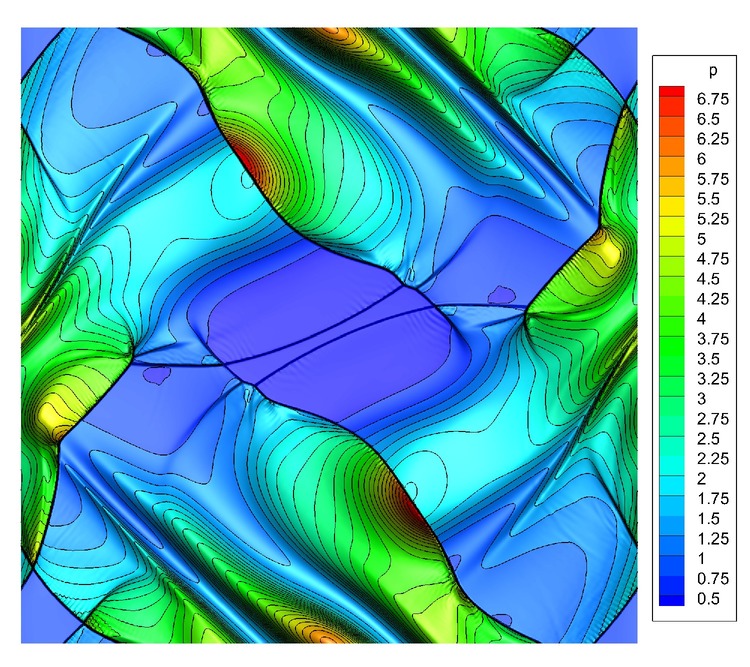}  \\   
\includegraphics[width=0.31\textwidth]{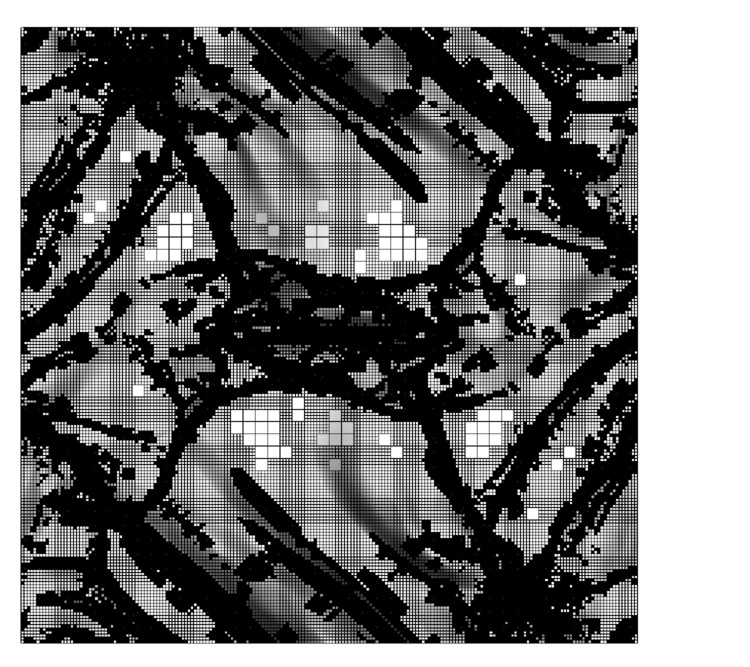}     &  
\includegraphics[width=0.31\textwidth]{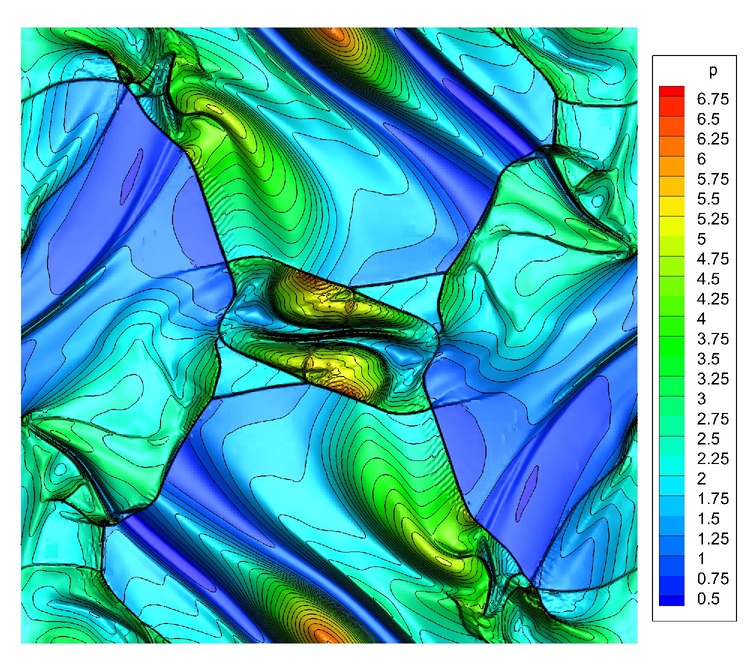}      &  
\includegraphics[width=0.31\textwidth]{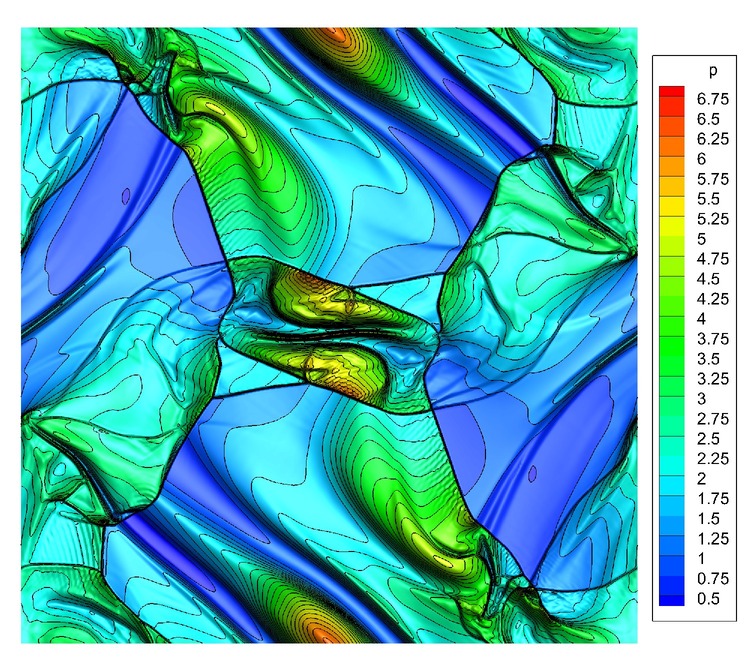}  \\   
\includegraphics[width=0.31\textwidth]{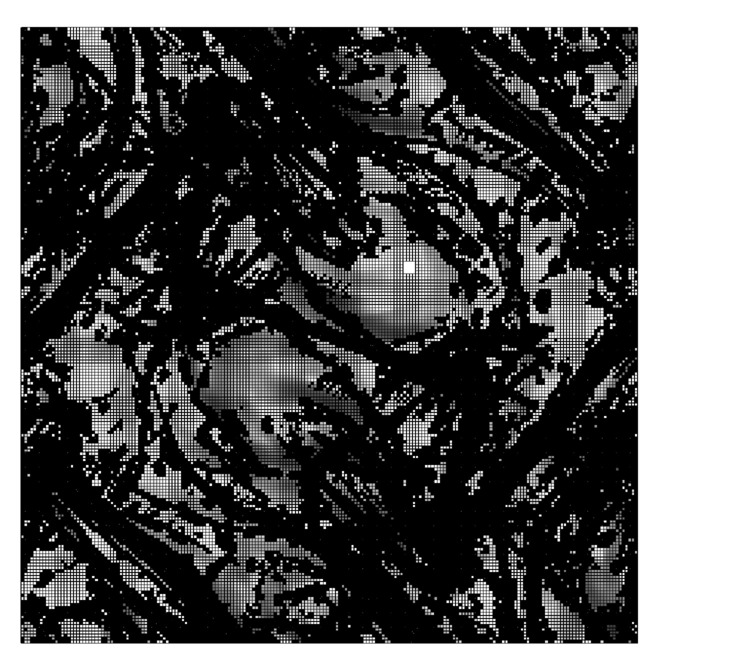}     &  
\includegraphics[width=0.31\textwidth]{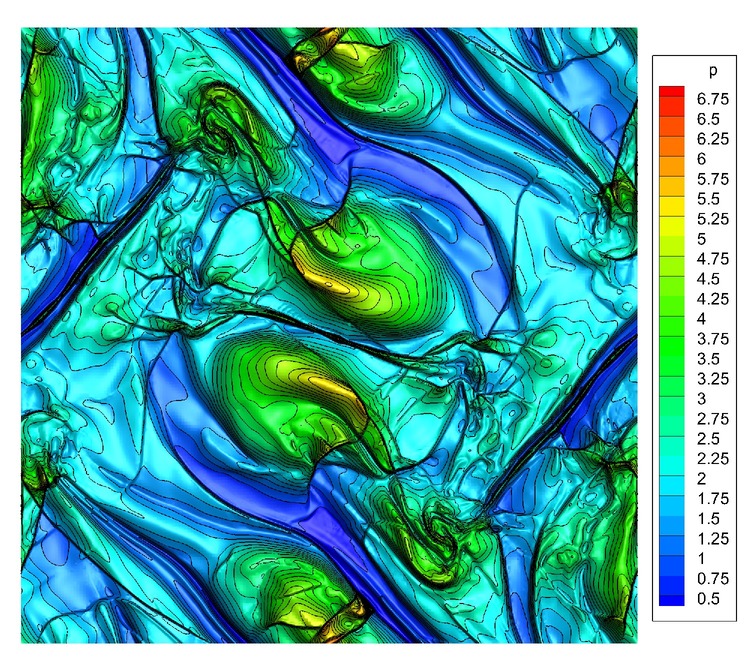}      &  
\includegraphics[width=0.31\textwidth]{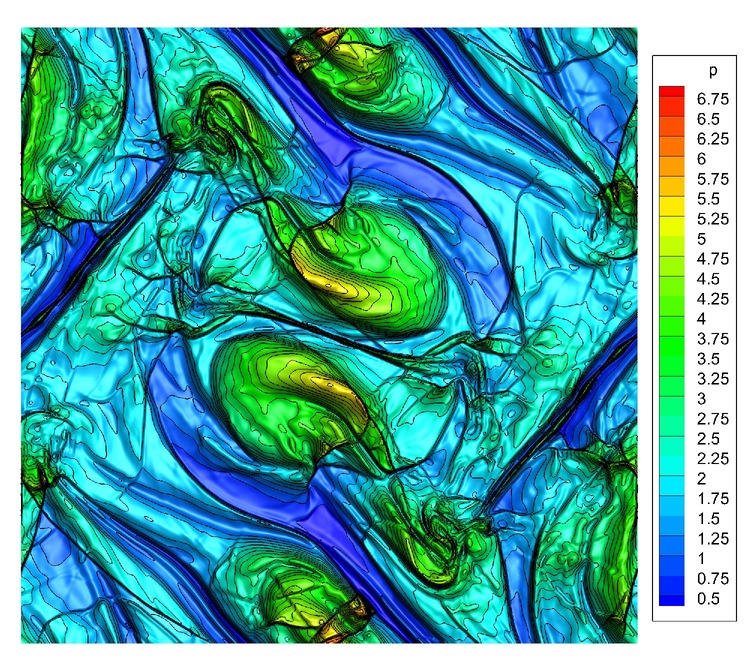}      
\end{tabular} 
\caption{ Orszag-Tang vortex system at times $t=0.5$, $t=2.0$, $t=3.0$ and $t=5.0$ from top to bottom. AMR grid (left), third order ADER-WENO solution obtained on the AMR grid (center) and on a fine uniform 
grid corresponding to the finest AMR grid level (right).} 
\label{fig.ot}
\end{center}
\end{figure}

\paragraph*{Orszag-Tang vortex system.} 
The first test case considered for the ideal MHD equations is the classical vortex system 
of Orszag and Tang \cite{OrszagTang} which was studied extensively in \cite{PiconeDahlburg} and \cite{DahlburgPicone}. The computational 
domain is $\Omega = \left[0;2\pi\right]^2$. We use the parameters of the computation of Jiang and Wu \cite{JiangWu}, scaling the magnetic 
field by $\sqrt{4\pi}$ due to the different normalization of the governing equations. The initial condition of the problem is given by  
\begin{equation}
  \left(\rho,u,v,p,B_x,B_y\right) = \left( \gamma^2, -\sin(y), \sin(x), \gamma, -\sqrt{4\pi} \sin(y), \sqrt{4\pi} \sin(2x)  \right), 
\end{equation}
with $w=B_z=0$ and $\gamma = \frac{5}{3}$. The divergence cleaning speed is set to $c_h=2.0$. The problem is solved up to $t=5.0$ using a 
third order ADER-WENO scheme with componentwise WENO reconstruction. The initial mesh on level zero is composed of $50 \times 50$ elements. 
We furthermore use $\mathfrak{r}=4$ and $\ell_{\max}=2$. This corresponds to an equivalent resolution on a uniform fine mesh of $800 \times 800$. 
To assess the accuracy and efficiency of our proposed AMR scheme, we also run a simulation on the uniform fine mesh for comparison. 
The results for pressure are shown in Fig. \ref{fig.ot} for $t=0.5$, $t=2.0$, $t=3.0$ and $t=5.0$, both, for the AMR grid as well as for the 
uniform grid corresponding to the finest AMR grid level. Our results are in agreement with the fifth order WENO finite difference solution computed by 
Jiang and Wu \cite{JiangWu} and with the unstructured third order WENO solution depicted in \cite{Dumbser2008} for the same problem. Furthermore, the 
AMR computations are in excellent agreement with the uniform fine grid reference solution. An efficiency comparison concerning memory requirements and 
CPU time is listed in Table \ref{tab.ot.compare}. The AMR method needs 322 time steps on the coarsest mesh, while the simulation on the uniform fine 
mesh needs 4148 time steps to reach the same final time. 
Even for this test problem, where most of the cells are refined and therefore only little gain is 
expected through the use of AMR techniques, we obtain still a speedup of a factor of 1.8 compared to the uniform fine mesh simulation. 

\begin{table}[!b]   
\caption{Memory and CPU time comparison of the third order ADER-WENO AMR method and ADER-WENO on a uniform fine grid for the Orszag--Tang problem.
         Memory consumption is measured in maximum number of elements and CPU time is normalized with respect to the simulation on the fine uniform mesh.} 
\begin{center} 
\renewcommand{\arraystretch}{1.0}
\begin{tabular}{cccc} 
\hline
     & AMR & Uniform &  ratio \\ 
\hline
Cells & 454525 & 640000  & 1.41 \\ 
CPU   & 0.547  & 1.0     & 1.83 \\ 
\hline
\end{tabular} 
\end{center}
\label{tab.ot.compare}
\end{table} 

\begin{figure}
\begin{center}
\begin{tabular}{lr}
\includegraphics[width=0.4\textwidth]{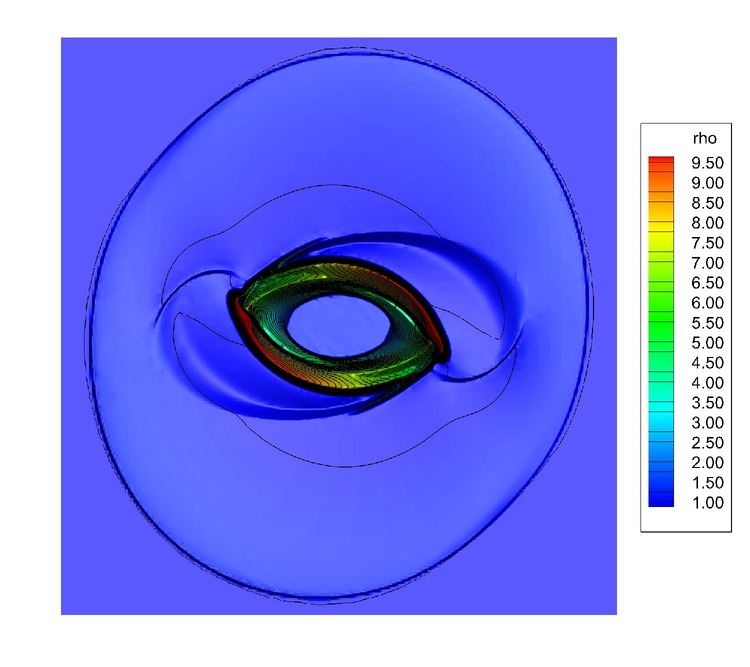}      &  
\includegraphics[width=0.4\textwidth]{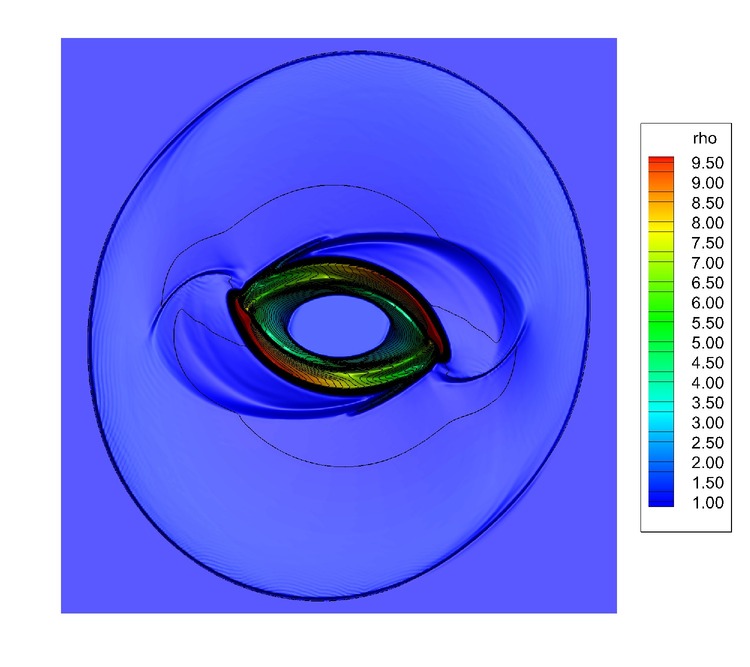}  \\   
\includegraphics[width=0.4\textwidth]{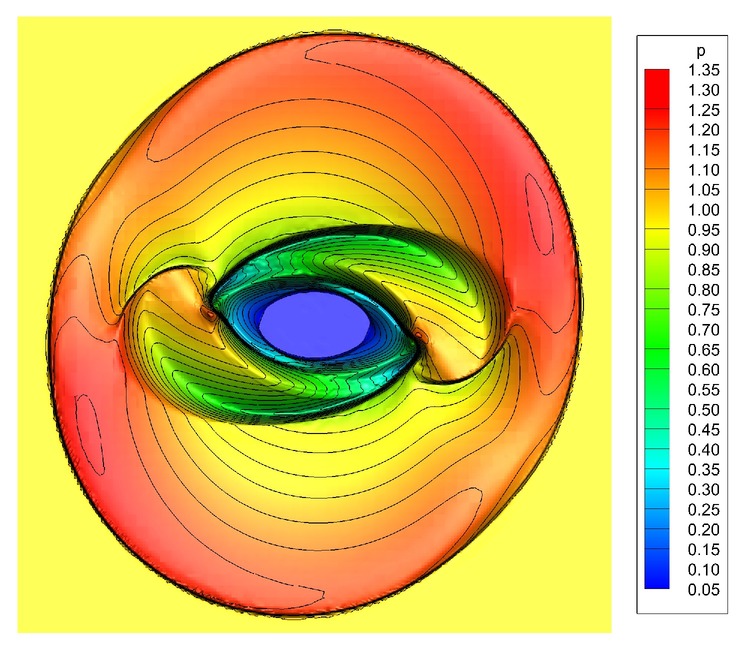}      &  
\includegraphics[width=0.4\textwidth]{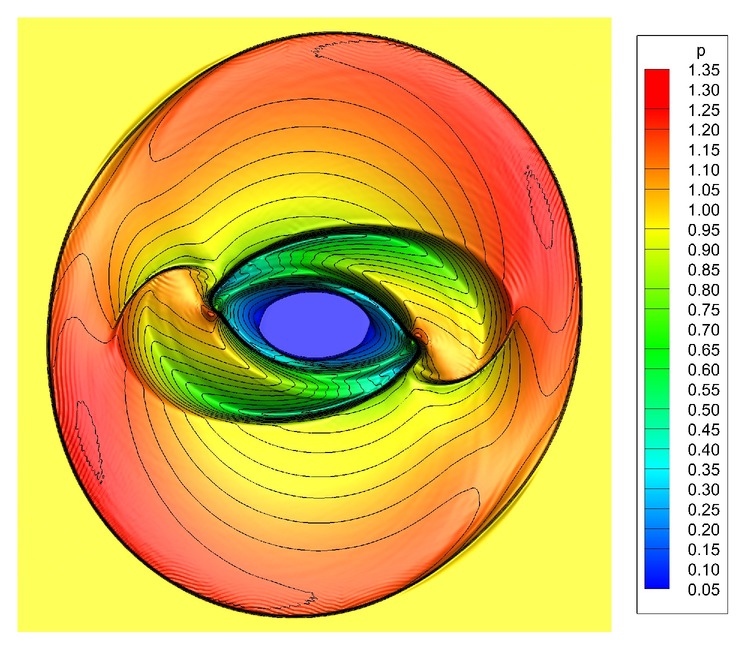}  \\   
\includegraphics[width=0.4\textwidth]{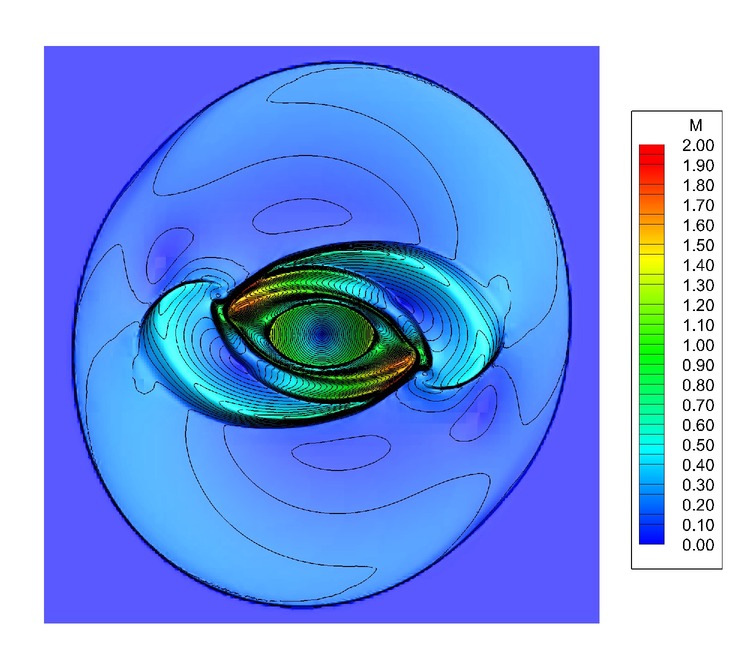}      &  
\includegraphics[width=0.4\textwidth]{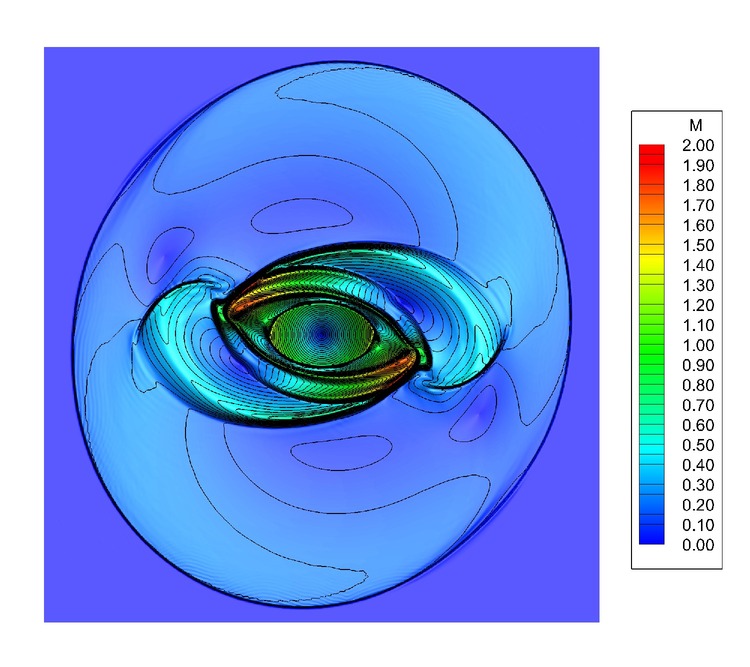}  \\   
\includegraphics[width=0.4\textwidth]{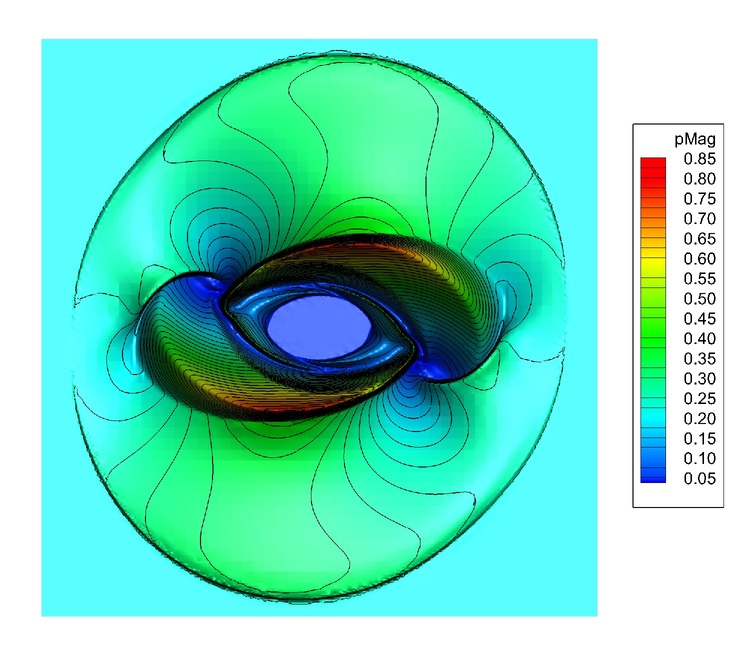}      &  
\includegraphics[width=0.4\textwidth]{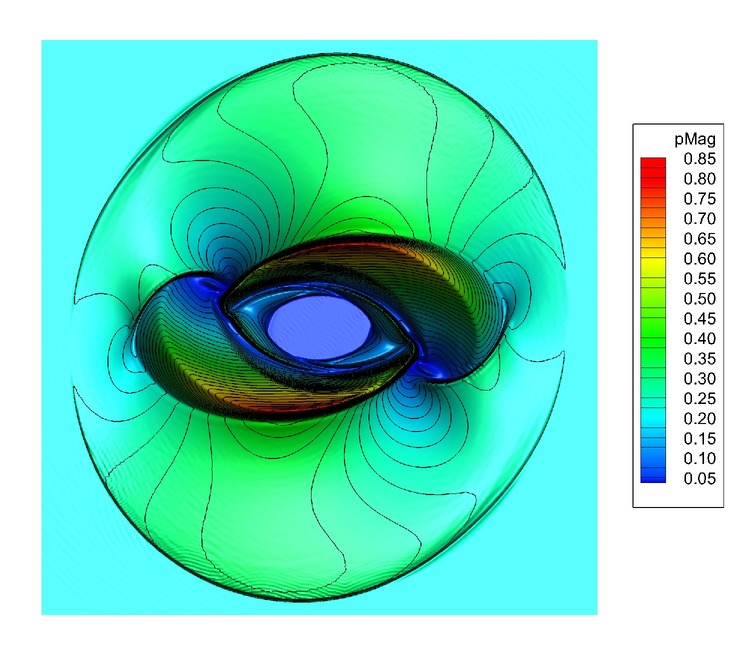}      
\end{tabular} 
\caption{ MHD rotor problem at time $t=0.25$. Third order ADER-WENO solution obtained on the AMR grid (left) and on a fine uniform 
grid corresponding to the finest AMR grid level (right).} 
\label{fig.mhdrotor}
\end{center}
\end{figure}

\begin{figure}
\begin{center}
\includegraphics[width=0.65\textwidth]{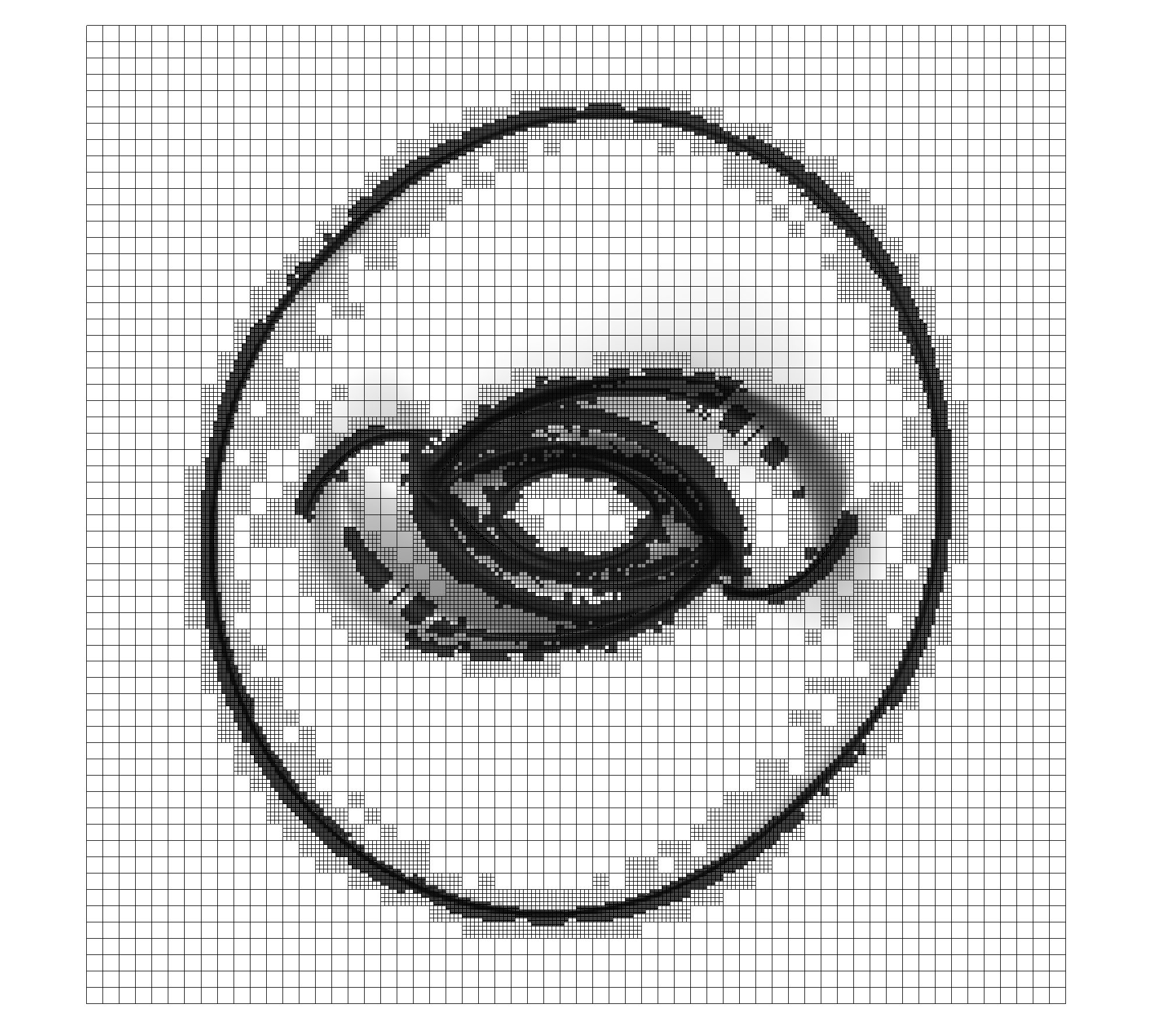}        
\caption{ AMR grid for the MHD rotor problem at time $t=0.25$.} 
\label{fig.mhdrotor.grid}
\end{center}
\end{figure}

\paragraph*{MHD rotor problem.} 
The second test case is the well-known MHD rotor problem proposed by Balsara and Spicer in \cite{BalsaraSpicer1999}. It consists of a rapidly rotating 
fluid of high density embedded in a fluid at rest with low density. Both fluids are subject to an initially constant magnetic field. The rotor 
causes torsional Alfv\'en waves to be launched into the fluid at rest. As a result the angular momentum of the rotor is 
diminished. The problem is set up on a computational domain $\Omega=[-0.6;0.6] \times [-0.6;0.6]$, using a third order ADER-WENO scheme. The AMR mesh 
on level zero contains $60 \times 60$ elements. With $\mathfrak{r}=4$ and $\ell_{\max}=2$ this simulation corresponds to a uniform fine mesh with  
$960 \times 960$ points resolution. As before, a uniform fine grid simulation is also performed to assess accuracy and efficiency of the proposed 
ADER-WENO scheme on AMR grids. The initial density of the rotor is $\rho=10$ for $0 \leq r \leq 0.1$ and $\rho=1$ for the ambient fluid. The 
rotor has a constant angular velocity $\omega$ that is determined in such a way to obtain a toroidal velocity of $v=\omega \cdot r = 1$ at $r=0.1$. 
The pressure is $p=1$ in the whole domain and the magnetic field vector is set to $\vec B = ( 2.5, 0, 0)^T$ in the entire domain. As proposed by 
Balsara and Spicer we apply a linear taper to the velocity and density field, however only in a very small range $0.1 \leq r \leq 0.105$ so 
that density and velocity match those of the ambient fluid at rest at a radius of $r=0.105$. The speed for the hyperbolic divergence cleaning 
is set to $c_h=2$ and $\gamma=1.4$ is used. Transmissive boundary conditions are applied at the outer boundaries. The final AMR mesh is depicted in 
Fig. \ref{fig.mhdrotor.grid}. The computational results on the AMR mesh are compared with those on the uniform fine mesh at time $t=0.25$ in Fig. 
\ref{fig.mhdrotor} for density, pressure, Mach number and magnetic pressure. One observes that both solutions agree very well with each other. Also 
compared to the results presented by Balsara and Spicer we note a very good agreement. We emphasize that thanks to the divergence cleaning, no 
spurious oscillations can be seen in the density field and in the magnetic pressure, as reported by Balsara and Spicer for Godunov schemes without 
divergence cleaning. 
The AMR method needs only 99 time steps on the coarsest mesh, while the simulation on the uniform fine mesh solution needs 1147 time steps to reach 
the same final time. 

In this test problem, the efficiency gain of AMR is particularly evident. The computation on a fine uniform mesh corresponding to the finest AMR level 
needs more than five times more elements and more than seven times more CPU time, see the detailed results reported in Table \ref{tab.rotor.compare}.

\begin{table}[!b]   
\caption{Memory and CPU time comparison of the third order ADER-WENO AMR method and ADER-WENO on a uniform fine grid for the MHD rotor problem.
         Memory consumption is measured in maximum number of elements and CPU time is normalized with respect to the total wallclock time on the 
         uniform mesh.} 
\begin{center} 
\renewcommand{\arraystretch}{1.0}
\begin{tabular}{cccc} 
\hline
     & AMR & Uniform &  ratio \\ 
\hline
Elements   & 179680 & 921600  & 5.13 \\ 
CPU time   & 0.140  & 1.0     & 7.14 \\ 
\hline
\end{tabular} 
\end{center}
\label{tab.rotor.compare}
\end{table}

\section{Conclusions}
\label{sec:concl}

In this article we have presented the first better than second order one--step ADER--WENO 
finite volume scheme on space--time adaptive AMR grids. The use of a high order one--step 
time stepping method, based here on a local space--time discontinuous Galerkin predictor, 
allows a straightforward implementation of time accurate local time stepping, where each 
AMR grid level runs on its own local time step. Furthermore, compared to the method of lines based
on Runge--Kutta time stepping, the use of a high order one--step scheme in time reduces 
the number of nonlinear WENO reconstructions and the number of necessary MPI communications.  
All these key features of our present scheme help to keep the overall overhead associated 
with the administration of the space--time adaptive mesh at a reasonable level, at most
$25\%$, as quantified in Table \ref{tab.efficiency}. 

We have carried out numerical convergence studies, confirming that the claimed higher order
in space and time is actually reached in practice. Furthermore, the scheme has been applied 
to a series of test problems in two and three space dimensions, solving the compressible Euler 
equations as well as the classical MHD equations. In our examples, it was clearly shown that 
also for better than second order schemes, the use of AMR is beneficial, compared to the use
of a uniform fine grid. We have also shown via numerical evidence, that even in the AMR 
context the use of higher order schemes is beneficial, in particular when small scale turbulent
structures, vortices and sound waves have to be resolved. All these physical phenomena require 
little numerical dissipation for their efficient simulation. 

Future research will concern the extension of the present method to the general family of 
$P_NP_M$ schemes introduced in \cite{Dumbser2008}, as well as the simulation of realistic  
problems in computational astrophysics. In further work we plan to extend the present scheme 
also to turbulent viscous flows, to chemically reacting multiphase flows as well as to 
nonconservative hyperbolic systems. 

\section*{Acknowledgments}

The research conducted here has been financed in parts by the European Research Council under 
the European Union's Seventh Framework Programme (FP7/2007-2013) in the frame of the research 
project \textit{STiMulUs}, ERC Grant agreement no. 278267. A.H. thanks Fundaci\'on Caja Madrid 
(Spain) for its financial support by a grant under the programme \textit{Becas de movilidad 
para profesores de las universidades p\'ublicas de Madrid}. 

\bibliography{AMR3D}
\bibliographystyle{plain}

\end{document}